\documentclass[leqno]{article}

\usepackage[utf8x]{inputenc}
\usepackage[english]{babel}

	\usepackage{stmaryrd}
	\usepackage{dsfont}
	\usepackage{mathrsfs}
	\usepackage{eucal}
	\usepackage{latexsym}
	\usepackage{yhmath}
	\usepackage{calc}
	\usepackage{enumitem}
	\usepackage{tocloft}
		
	\usepackage[toc,page]{appendix}
	\usepackage{amsthm, amsmath, amssymb, amsfonts}
	
	\usepackage{hyperref}

	\newcommand{\N}{\mathds{N}}
	\newcommand{\Z}{\mathds{Z}}
	\newcommand{\R}{\mathds{R}}
	
	\newcommand{\ds}{\displaystyle}

	\newcommand{\energyFlow}[2][]{\mathscr{E}#1 \left( #2 \right)}
	\newcommand{\energyDir}[2][]{\mathscr{D}#1 \left( #2 \right)}
	\newcommand{\effCon}[2][]{\mathscr{C}#1\left( #2 \right)}
	\newcommand{\effRes}[2][]{\mathscr{R}#1\left( #2 \right)}
	\newcommand{\WSF}{\mathsf{WSF}}
	\newcommand{\UST}{\mathsf{UST}}
	\newcommand{\USF}{\mathsf{USF}}
	\newcommand{\FSF}{\mathsf{FSF}}
	\newcommand{\LE}{\mathsf{LE}}
	
	\newcommand{\NRW}{\mathsf{NRW}}

	\newcommand{\strA}{\mathrm{A}}
	\newcommand{\strB}{\mathrm{B}}
	\newcommand{\strC}{\mathrm{C}}
	\newcommand{\strD}{\mathrm{D}}
	\newcommand{\strE}{\mathrm{E}}
	\newcommand{\strF}{\mathrm{F}}
	\newcommand{\strG}{\mathrm{G}}
	\newcommand{\strH}{\mathrm{H}}
	
	\newcommand{\strJ}{\mathrm{J}}
	\newcommand{\strK}{\mathrm{K}}
	\newcommand{\strL}{\mathrm{L}}
	
	\newcommand{\strN}{\mathrm{N}}
	
	\newcommand{\strP}{\mathrm{P}}
	\newcommand{\strQ}{\mathrm{Q}}
	\newcommand{\strR}{\mathrm{R}}
	\newcommand{\strS}{\mathrm{S}}
	\newcommand{\strT}{\mathrm{T}}
	\newcommand{\strU}{\mathrm{U}}
	\newcommand{\strV}{\mathrm{V}}
	\newcommand{\strW}{\mathrm{W}}

	\newcommand{\scrC}{\mathscr{C}}
	
	\newcommand{\scrE}{\mathscr{E}}
	\newcommand{\scrF}{\mathscr{F}}
	\newcommand{\scrG}{\mathscr{G}}

	\newcommand{\scrL}{\mathscr{L}}
	
	\newcommand{\scrN}{\mathscr{N}}
	
	\newcommand{\scrP}{\mathscr{P}}
	
	\newcommand{\scrR}{\mathscr{R}}
	\newcommand{\scrS}{\mathscr{S}}
	\newcommand{\scrT}{\mathscr{T}}

	\newcommand{\scrX}{\mathscr{X}}

		\newcommand{\card}[1]{\mathrm{card}\left(#1\right)}
		\newcommand{\pr}[1]{\mathrm{pr}_{#1}}
		
		\newcommand{\norm}[1]{\left\lVert #1 \right\rVert}
		\newcommand{\prodEsc}[2]{\left( #1\middle| #2 \right)}
		
		\newcommand{\Probability}[2][]{\mathds{P}{#1}\left(#2\right)}
		\newcommand{\Expectation}[2][]{\mathds{E}{#1}\left(#2\right)}
		
		\newcommand{\Covariance}[1]{\mathds{C}\mathrm{ov}\left(#1 \right)}
		\newcommand{\indic}[1]{\mathds{1}_{#1}}
		
		\newcommand{\integral}[1][]{\int\limits_{#1}}
		\newcommand{\bigO}[2][]{\mathrm{O}#1\left(#2\right)}
		\newcommand{\littleO}[2][]{\mathrm{o}#1\left(#2\right)}

		\newcommand{\japBracket}[2][]{\left\langle #2 \right\rangle#1}
		\newcommand{\spread}[2][]{\japBracket[_\eta]{#2}}

		\newcommand{\ball}[2][]{\mathrm{B}#1\left(#2\right)}
		
		\newcommand{\sphere}[2][]{\mathrm{S}#1\left(#2\right)}
		\newcommand{\frontier}[1]{\partial #1}
		
		\newcommand{\closure}[1]{\overline{#1}}
		
		\newcommand{\avec}[3]{{#1}_{#2}, \ldots, {#1}_{#3}}
		\newcommand{\supvec}[3]{{#1}^{#2}, \ldots, {#1}^{#3}}
		
		\newcommand{\distNorm}[2]{\mathsf{Norm}\left(#1;#2\right)}
		
\usepackage[headings]{fullpage}

\usepackage{fancyhdr}
	\newtheoremstyle{MiEstilo}	
		{\topsep}	
		{\topsep}	
		{\em}	
		{}				
		{\scshape}	
		{ )}			
		{ }			
		{}				
	\theoremstyle{MiEstilo}	
	\newcounter{auxCounter}
	\newcounter{auxSubCounter}
	\numberwithin{equation}{section}
	\newtheorem{Theorem}[equation]{Theorem (}
	
	\newtheorem{Proposition}[equation]{Proposition (}
	\newtheorem{SubCounter}{}[equation]
	\newtheorem{EqInside}[SubCounter]{(}
	\newtheorem{Lemma}[SubCounter]{Lemma (}
	
	\newtheorem{Corollary}[SubCounter]{Corollary (}
			\newtheoremstyle{MiEstiloSinNumero}	
				{\topsep}	
				{\topsep}	
				{\em}	
				{}				
				{\scshape}	
				{.}			
				{ }			
				{}				
			\theoremstyle{MiEstiloSinNumero}	
			\newtheorem*{Theorem*}{Theorem}
	\newtheoremstyle{DefRem}	
	{\topsep}		
	{\topsep}		
	{\normalfont}			
	{}				
	{\scshape}	
	{ )}				
	{ }				
	{}				
	\theoremstyle{DefRem}
	
	\newtheorem{Remark}[equation]{Remark (}
	\newtheoremstyle{Problem}	
		{\topsep}		
		{\topsep}		
		{\em}			
		{}				
		{\bfseries}	
		{.}				
		{ }				
		{}				
	\theoremstyle{Problem}

\usepackage{multicol}

\renewcommand{\labelenumi}{(\alph{enumi})}

\usepackage{graphicx}

\addtocontents{toc}{\protect\setcounter{tocdepth}{1}}

\begin{document}

	\title{\textbf{Uniform spanning forest on the integer lattice \\ with drift in one coordinate}}
	\author{\textsc{By Guillermo {Martinez Dibene}\thanks{e-mail: gmomtzcanada@gmail.com}}\\Mathematics Department \\ The University of British Columbia}
	\date{}
	\maketitle

	\begin{abstract}
		In this article we investigate the Uniform Spanning Forest $(\USF)$ in the nearest-neighbour integer lattice $\Z^{d + 1} = \Z \times \Z^d$ with an assignment of conductances that makes the underlying (Network) Random Walk ($\NRW$) drifted towards the right of the first coordinate. This assignment of conductances has exponential growth and decay; in particular, the measure of balls can be made arbitrarily close to zero or arbitrarily large. We establish upper and lower bounds for its Green's function. We show that in dimension $d = 1, 2$ the $\USF$ consists of a single tree while in $d \geq 3,$ there are infinitely many trees. We then show, by an intricate study of multiple $\NRW$s, that in every dimension the trees are one-ended; the technique for $d = 2$ is completely new, while the technique for $d \geq 3$ is a major makeover of the technique for the proof of the same result for the graph $\Z^d.$ We finally establish the probability that two or more vertices are $\USF$-connected and study the distance between different trees.

		\vspace{0.15cm}


		\vspace{0.15cm}

		\noindent \emph{Acknowledgements.} This research was partially supported by an NSERC (Canada) grant.
	\end{abstract}
		

	\tableofcontents

	\section{Introduction}
	\begin{center}
		\emph{All graphs are assumed directed, connected, without multiple edges or self-loops and denumerable. Given an edge $e,$ the direction of $e$ is given by a starting vertex $e^-$ and a ending vertex $e^+.$ All conductances are assumed strictly positive.}
	\end{center}
	This work is based on the PhD thesis of the author; should the reader need more details, consult \cite{MaDib:USFwithDrift}. A graph without cycles is called a \textbf{tree}; the tree is \textbf{spanning} if it contains all vertices of the graph. Given a finite connected graph $\strG,$ the number of subgraphs is also finite and thus, there are only finitely many spanning trees in that graph. The uniform spanning tree is, by definition, the probability measure on the set of spanning trees of $\strG$ satisfying $\mathds{P}(t) = \left( \# \text{ of spanning trees} \right)^{-1}.$	This can readily be generalised to the network case, $\mathds{P}(t) = k \cdot \mathsf{weight}(t),$ where $k$ is a normalising constant. Any of these measures (or the underlying random object) is called ``uniform spanning tree.''

	Suppose now that $(\strG, \mu)$ is an infinite network. The first paper to introduce $\UST$ on infinite graphs was Robin Pemantle's \cite{Pemantle:USFZd}. He considered $\Z^d$-balls corresponding to the norm $\ds \norm{x}_\infty = \max_{1 \leq i \leq d} |x_i|.$ He showed that the uniform spanning tree measure $\mu_n$ on the trees of $\strB_n$ converges weakly to some measure $\mu.$ This is the uniform spanning measure on $\Z^d.$ We will employ the notation $\USF$ to denote it.

	Later, Itai Benjamini, Russel Lyons, Yuval Peres and Oded Schramm, in \cite{BLPS:USF} and using tools unavailable to Pemantle, extended his result from $\Z^d$ to general networks. They started by considering an exhaustion of the graph created by considering increasing families of finite sets of vertices and then considering the networks induced on these sets. If $(\strG_n)$ is any such exhaustion, they showed that the uniform spanning forests measure of $\strG_n$ defines a sequence of measures that converge weakly to some measure $\mu^F.$ They named $\mu^F$ as ``free uniform spanning forest measure.'' The second way to construct the uniform spanning forest measure is by considering constraints in the boundary of $\strG_n.$ Indeed, the considered a sequence $(\strV_n)$ of families of finite sets of vertices and then constructed the network $\strG_n^W,$ which is the network induced on $\strV_n$ with the boundary of $\strV_n$ wired into a vertex. Again, they showed that the uniform spanning forest measures of $\strG_n^W$ converge to a measure $\mu^W$ and they called this measure the ``wired uniform spanning forest measure.'' Several properties of $\mu^F$ and $\mu ^W$ were studied, in particular, necessary and sufficient conditions were given as to when $\mu^F = \mu^W.$ We will employ $\FSF$ and $\WSF$ instead of $\mu^F$ and $\mu^W,$ respectively.

		\subsection{Basic definitions}
		Consider any network $(\strG, \mu).$ For a vertex $v,$ we define the \textbf{conductance} of $v$ by $\mu(v) = \sum\limits_{e^- = v} \mu(v).$ The network random walk ($\NRW$)  in $(\strG, \mu)$ is the Markov chain with state space $\strV$ (the vertex set of $\strG$) and transition probability $p(u, v) = \mu(u)^{-1} \left( \sum\limits_{e^- = u, e^+ = v} \mu(e) \right).$ If $(S_n)_{n \in \Z_+}$ is the $\NRW$ of $(\strG, \mu)$, we will write $\tau_\strH = \inf\{n \in \Z_+ \mid S_n \in \strH\}$ and $\tau_\strH^+ = \inf\{n \in \N \mid S_n \in \strH\}$ 	where $\strH$ is any set of vertices. For sake of simplicity, we write $\tau_a = \tau_{\{a\}},$ $\tau_a^+ = \tau_{\{a\}}^+.$ $\tau_\strH$ is the \textbf{hitting time} of $\strH$ and $\tau_a$ the \textbf{visit time} of $a;$ if $\strH = \strK^\complement,$ we will say \textbf{exit time} of $\strK.$

		Suppose $S = (S_n)_{n \in \Z_+}$ is the $\NRW$ on $\strG.$ Often  we will use $\Probability[^x]{S_n = y} = \Probability{S_n = y \mid S_0 = x}$ and $\Expectation[^x]{f(S_n)} = \Expectation{f(S_n) \mid S_0 = x}.$ The \textbf{Green's function} ``of the network $\strG$'' is, by definition $	G(x, y) = \sum\limits_{n \in \Z_+} \Probability[^x]{S_n = y} = \Expectation[^x]{\sum\limits_{n = 0}^\infty \indic{\{S_n = y\}}},$ for all $(x, y) \in \strV^2.$ Similarly, define the \textbf{Green's function restricted to leaving $\strV_0$} by $G_{\strV_0}(x, y) = \sum\limits_{m = 0}^\infty \Probability[^x]{S_m = y, m < \tau_{\strV_0}^\complement}.$
		
		Suppose that $S_0 = o$ is some fixed vertex of the network and suppose that $N = (N_t)_{t \geq 0}$ is a Poisson process in $\Z_+$ with unitary constant intensity which is independent of the process $S.$ We define the \textbf{continuous time $\NRW$} to be the process $\tilde{S} = \left( \tilde{S}_t \right)_{t \geq 0}$ given by $\tilde{S}_t = S_{N_t}.$ Define the \textbf{continuous time $\NRW$ transition probability} by $q_t^S(o, v) = \Probability[^o]{\tilde{S}_t = v},$ for $(t, z) \in [0, \infty) \times \strV.$

		Let $\gamma = (v_j)_{0 \leq j \leq m}$ be a path of vertices. The \textbf{loop-erasure} of $\gamma,$ denoted as $\LE(\gamma),$ is the path obtained from $\gamma$ by deleting cycles in chronological order. More specifically, it is defined as follows. Set $\gamma_0 = \gamma.$
		\begin{description}
			\item[$(\LE)$] Suppose we are at stage $k$ and we receive a path $\gamma_{k - 1} = \left( \avec{u}{1}{m_k} \right).$ If no vertex of this path is repeated, then return $\LE(\gamma) = \gamma_{k - 1}.$ Else, consider $\tau_k$ to be the \emph{first index} $j$ such that the vertex $u_j$ appears twice in $\gamma_{k - 1}.$ Then, let $\sigma_k$ to be the \emph{last index} $j$ such that $u_j = u_{\tau_k},$ consider now the path $\gamma_k = \left( \avec{u}{1}{\tau_k}, \avec{u}{\sigma_k + 1}{m_k} \right).$ Repeat.
		\end{description}

		A network $\Gamma$ is \textbf{transitive} if for every pair of vertices $x$ and $y$ there exists a network automorphism $\varphi$ of $\Gamma$ such that $\varphi(x) = y.$

		Let $\strA$ be a subset of vertices. The \textbf{interior} of $\strA$ is the set of all vertices of $\strA$ whose all neighbours are vertices of $\strA$ as well; we denote it with $\mathring{\strA}.$ The \textbf{closure} of $\strA,$ denoted $\closure{\strA},$ is the set of all vertices of $\strA$ or that are adjacent to $\strA.$ The \textbf{(exterior vertex) boundary} of $\strA,$ denoted $\partial \strA,$ is the set of all vertices \emph{outside} $\strA$ that are adjacent to $\strA.$

		Let $\strG$ be a network with vertex set $\strV.$ A \textbf{harmonic} function on the network $\strG$ is  any function $h:\strV \to \R$ such that $h(x) = \sum\limits_{y \sim x} p(x, y) h(y).$

		Let $(\strG, \mu)$ be a network with vertex set $\strV.$ A function $f:\strV \to \R$ is said to be a \textbf{Dirichlet} function if $\energyDir{f} = \sum\limits_{x, y} \mu(x, y) (f(x) - f(y))^2 < \infty.$ The finite or infinite number $\energyDir{f}$ is known as \textbf{Dirichlet energy} of $f.$

		\subsection{Wilson's algorithm}
		David Wilson published his algorithm in \cite{Wilson:GRST}. (This algorithm is for finite networks.)
		\begin{description}
			\item[$(\mathsf{WA})$] Order the vertices and set $t_0$ to be the tree consisting solely of the first vertex. Having defined $t_{k - 1},$ start an independent $\NRW$ from the first vertex not in $t_{k - 1}$ and run this random walk until it hits $t_{k - 1},$ consider the $\LE$ of this path and call $t_k$ the tree obtained as the union of this $\LE$-path with $t_{k - 1}.$ When all vertices have been searched, return $\mathfrak{T}$ the final tree constructed.
		\end{description}
		It is known that $\mathfrak{T}$ follows the $\UST$ law of the given (finite) network. A typical way of using Wilson's algorithm is via a more intricate structure; while Wilson introduced this structure, a more detailed exposition can be found in \cite{Barlow:LEWUST}. We refer to this version of Wilson's algorithm as ``Wilson's algorithm with stacks.'' See also \cite{LyPe:PTN}, Ch. 4. With the definitions of these references, we have the following result.
		\begin{Proposition}\label{Proposition: WA with stacks}
			Assume $(\strG, \mu)$ is a finite network. Fix a root $r \in \strV.$ With probability one, there are finitely many cycles that can ever be popped. Any two ordering of the searches of the vertices in $\strV \setminus \{r\}$ will pop all these cycles and the popping procedure will leave the same visible graph for the two orderings. This visible graph sitting on top of the stacks is a spanning tree and the distribution of this random tree is that of the $\UST$ of the given network.
		\end{Proposition}
		Observe that Wilson's algorithm run with stacks on a finite graph produces always the same tree provided the root is fixed in advance.

		We mention now the usual way to construct $\WSF$ on an infinite network $(\strG, \mu).$ Consider an exhaustion $\strV_q$ of the vertex set and denote by $(\strG_q, \mu)$ the network induced on $\strV_q$ with its boundary $\frontier \strV_q$ wired into a vertex. Run Wilson's algorithm with stacks taking $\frontier \strV_q$ as root. This produces a random tree $\mathfrak{T}_q.$ Then, $\mathfrak{T}_q \xrightarrow[q \to \infty]{\mathrm{weakly}} \mathfrak{F},$ where $\mathfrak{F}$ has distribution $\WSF,$ see \cite{LyPe:PTN}, sect. 10.1. There exists an alternative more dynamic (it provides a.s. convergence as opposed to weakly) method to sample the $\WSF$-distributed random object, which was first introduced in \cite{BLPS:USF}.
	
		Suppose $\strG$ is any network. Denote by $\strV$ the set of vertices of $\strG$ and by $\strE$ that of edges. Suppose $\xi:\N \to \strV$ is a bijection (an ``ordering'' of the vertex set). Let $(\Omega, \scrF, \mathds{P})$ be a probability space in which is possible to define an independent family $(S^v)_{v \in \strV}$ of $\NRW$s of this network, such that $S^v$ is a random walk started at $v.$ Define inductively the following random spanning subgraphs of $\strG:$
		\begin{description}
			\item[Start.] $L_1^\xi = \mathfrak{F}_1^\xi = \LE \left( S^{\xi(1)}_m \right)_{m \in \Z_+}.$
			\item[Inductive step.] Having defined $L_k^\xi$ and $\mathfrak{F}_k^\xi,$ define $L_{k + 1}^\xi = \LE \left( S^{\xi(k + 1)}_m \right)_{m = 0, \ldots, \tau_k}$ and $\mathfrak{F}_{k + 1} = \mathfrak{F}_k \cup L_{k + 1}^\xi,$ where $\tau_k = \tau_{\mathfrak{F}_k^\xi}\left( S^{\xi(k + 1)} \right),$ that is, $\tau_k$ is the hitting time of the forest currently present $\mathfrak{F}_k^\xi$ by the random walk $S^{\xi(k + 1)}.$
		\end{description}
		 Consider finally the random spanning subgraph of $\strG:$ $\mathfrak{F}^\xi = \bigcup\limits_{k \in \N} \mathfrak{F}_k^\xi.$ What Wilson's algorithm rooted at infinity \cite[Theorem 5.1]{BLPS:USF} states is that $\mathfrak{F}^\xi \sim \WSF_{\strG},$ in other words, the distribution of $\mathfrak{F}^\xi$ is that of the wired spanning forest of the network $\strG,$ and this happens independently of the choice of ordering $\xi.$

		\subsection{Definition of $\Gamma_d(\lambda)$ and summary of this paper}
		In this paper we will study $\Gamma_d(\lambda) = (\Z^{d + 1}, \mu),$ where $\mu$ is given  by
		\begin{equation}\label{Equation: Conductances of the network of interest}
			\mu((n, x), (n', x')) = e^{\lambda \max(n, n')}, \quad \text{ for all } (n, x) \sim (n', x').
		\end{equation}
		It can easily be seen that this assignment of conductances makes the underlying $\NRW$ to have a uniform drift to the right of the first axis. We will then investigate the main basic properties of the forest. In the remainder of this section we provide an overview of the results of this paper.

		In sect. \ref{Section: Greens function} we will prove the following theorem
		\begin{Theorem*}{\normalfont (Theorem (\ref{Theorem: Greens function bounds}))}
			There exist four constants $c_i$ ($i = 1, \ldots, 4$) such that for all vertices $z = (n, x) \neq 0$ of the network $\Gamma_d(\lambda),$
		\begin{displaymath}\textstyle
			G(0, z) \leq c_1 \begin{cases} e^{-c_2 \norm{z}} &\text{ if } \norm{x} > n, \\
			e^{-c_2 \frac{\norm{x}^2}{n}} \norm{z}^{-\frac{d}{2}} &\text{ if } \norm{x} \leq n,
			\end{cases}
		\end{displaymath}
		and
		\begin{displaymath}\textstyle
			G(0, z) \geq c_3 \begin{cases} e^{-c_4 \norm{z}} &\text{ if } \norm{x} > n, \\
			e^{-c_4 \frac{\norm{x}^2}{n}} \norm{z}^{-\frac{d}{2}} &\text{ if } \norm{x} \leq n,
			\end{cases}
		\end{displaymath}
		\end{Theorem*}
		To prove this bounds we will employ a decomposition of the Green's function using the continuous-time $\NRW$ and then use several ``well-known'' estimates on the latter.

		In sect. \ref{Section: Number of trees} we will establish is the following.
		\begin{Theorem*}{\normalfont (Theorem (\ref{Theorem: Equality between wired and free uniform spanning forest}))}
			On $\Gamma_d(\lambda),$ we have $\WSF = \FSF.$
		\end{Theorem*}
		This is done using the Liouville property and a coupling argument. This theorem then shows that $\USF$ is well defined on $\Gamma_d(\lambda).$ Next, using Fourier inversion theorem together with Plancherel's theorem, we will prove the following result, which we call ``bubble condition.''
		\begin{Theorem*}{\normalfont (Theorem (\ref{Theorem: The bubble condition}))}
			The Green's function  $z \mapsto G(0, z)$ of $\Gamma_d(\lambda)$ belongs to $\scrL^2$ if and only if $d \geq 3.$
		\end{Theorem*}
		 Shortly after, as a corollary of the bubble condition, we will also show that $\USF$ consists of a single tree when $d = 1, 2,$ and it consists of infinitely many infinite trees for $d \geq 3.$ 
		\begin{Theorem*}{\normalfont (Theorem (\ref{Theorem: Number of components in the forest}))}
			In $d = 1,2,$ $\WSF$ of $\Gamma_d(\lambda)$ is a.s. one tree; in $d \geq 3,$ there are a.s. infinitely many trees.
		\end{Theorem*}
		
		In sect. \ref{Section: Crossings of NRW} we estimate the frequency at which two independent $\NRW$ of $\Gamma_d(\lambda)$ will cross each other. We will prove the following result
		\begin{Theorem*}{\normalfont (Theorems (\ref{Theorem: In high dimensions two independent NRW may not cross at all}) and (\ref{Theorem: In dimension d = 1 the difference of two random walks is recurrent}), and Corollary (\ref{Corollary: In dimension d = 2 there are infinitely many crossings}))}
		Let $S$ and $S'$ be two independent network random walks in $\Gamma_d(\lambda)$ started at $0.$
		\begin{enumerate}
			\item If $d = 1,$ then a.s. there exists infinitely many $n$ such that $S_n = S_n' = 0.$
			\item If $d = 2,$ then a.s. there exist infinitely many pairs $(n, n')$ such that $S_n = S_{n'}'.$
			\item If $d \geq 3,$ there there exists a positive probability that for no pair $(n, n'),$ $S_n = S_{n'}'.$
		\end{enumerate}
	\end{Theorem*}
	In fact, in $d = 2,$ we will estimate the \emph{probability} that they will cross is roughly $\log n$ if they walk $n$ steps (a precise formulation of this is Theorem (\ref{Theorem: Logarithmic scale of crossings})). The proofs of dimensions $d \neq 2$ are relatively easy. The proof of dimension $d = 2$ is rather involved. We explain the gist of the idea now. We will consider a paraboloid given by the vertices $z = (n, x) \in \Z_+ \times \Z^d$ such that $\norm{x} \leq n^2$ and inside each of these, we will consider cylinders (\ref{Equation: Definition of the testing regions}) of appropriate dimensions. We will show that for certain regions (\ref{Equation: Region of restriction for the Green function}), the $\NRW$ started in the boundary of these regions will not go back to the testing cylinder (except with a small probability that can be controlled) and therefore, the Green's function virtually does not change value when we restrict it to leave these regions. We then employ a second-moment method to estimate the aforementioned probability. The calculations are quite tight and require a careful set up.

		In sect. \ref{Section: One endedness}, we prove
		\begin{Theorem*}{\normalfont (Theorems (\ref{Theorem: One end in dimension one}) and (\ref{Theorem: One end in dimension two}) and (\ref{Theorem: One end in high dimensions}))}
			In $\Gamma_d(\lambda),$ a.s. every tree in $\USF$ has one end.
		\end{Theorem*}	
		The proofs of dimensions $d = 1,$ $d = 2$ and $d \geq 3$ need to be done separately. In $d = 1$ we use planar duality and the rather easy-to-establish fact that $\Gamma_d(\lambda)$ is essentially its own dual; the idea is a makeover of the proof of the same result for the graph $\Z^2.$ In $d \geq 3,$ we followed closely the proof of one-endedness for the graphs $\Z^d$ given in \cite{LyPe:PTN}. We made adaptation according to our conveniences but the underlying idea is that should a tree have two ends, there will be strong conductivity in the boundary of boxes. The proof of case $d = 2$ is perhaps the most innovative contribution of this paper and was introduced from scratch. The first property we establish is that Wilson's algorithm rooted at infinity has some stochastic-stability (the law remains unaltered) if we fix the first steps and reorder future vertices (a stronger result for finite-graphs is proposition (\ref{Proposition: WA with stacks})). This will show that probabilities of future events depending on the past do not actually depend of the past. The next step is to show that it suffices to prove that the component of zero is one-ended. To this aim, we employ Wilson's algorithm rooted at infinity and take the first step of the algorithm to be the origin. Call $L_0$ this part of the forest. If the component at zero had at least two ends, then it would cross the boundary of each cylinder (centred at the origin) in at least two vertices. The aforementioned invariance of the law from reordering future searches having fixed the past will allow us to fix $L_0$ and then prove that, in a given a large cylinder $\strC = \left\{|n| \leq r, \norm{x} \leq r\right\},$ most vertices on $\partial \strC$ have infinitesimal probability to connect to zero, this will reduce the work to study a part of the left base of $\strC,$ namely $\strC_{r, 0} = \left\{n = -r, \norm{x} \leq r^{\frac{1}{2}} \right\}.$ We will then take a sparse subset here $(\strC_{r, 0}')$ and construct the branches from its vertices. The subset $\strC_{r, 0}'$ needs to be sparse enough so that the probability that the $\NRW$ started from its vertices create a second end from the origin is small. If now $z$ is any vertex in $\strC_{r, 0} \setminus \strC_{r, 0}',$ then we will show that the $\NRW$ started at $z$ has an overwhelming probability to hit the branch created by some $z' \in \strC_{r, 0}'$ and since this branch does not create a second end from the origin, $z$ will also not create a second end. The proof will conclude by gluing together all the estimates. Admittedly, the proof is elaborate and to ease the reading, we have provided well-detailed arguments.

		In sect. \ref{Section: Probability of same component} we will show
		\begin{Theorem*}{\normalfont (Theorems (\ref{Theorem: Probability of two points belonging in the same component}) and (\ref{Theorem: The probability of finitely many vertices in the same component}))}
			There exists a metric $\eta(z) = \max \left( |n|^{\frac{1}{2}}, \norm{x} \right)$ such
			\begin{displaymath}
				\Probability{z\text{ and }z'\text{ are in the same }\USF\text{-component}} \asymp \eta(z - z')^{-(d - 2)}.
			\end{displaymath}
		\end{Theorem*}
		To prove this, we follow similar techniques as those established in \cite{BKPS:GeomUSF}.  In a nutshell, we will need to investigate bounds of the convolution $\ds \sum_{z'' \in \Z^{d + 1}} G(z, z'') G(z', z'').$
	
		In sect. \ref{Section: Separation between components} we will establish the following
		\begin{Theorem*}{\normalfont (Theorem (\ref{Theorem: The separation between components}))}
			Let $D(z, z')$ be the minimal number of $\Z^{d + 1}$-edges outside the $\USF$-forest that connects the tree at $z$ and that at $z'.$ Almost surely, $\ds \max_{z, z'} D(z, z') = \left\lceil \dfrac{d - 2}{4} \right\rceil.$
		\end{Theorem*}
		The proof of this is done by adaptation of most of the methods of \cite{BKPS:GeomUSF} to the network $\Gamma_d(\lambda).$

		\subsubsection{Biased random walk}
		We make special mention of an article appearing in the arXiv.org in 2018. In \cite{SSSWX:BNRW} Zhan Shi, Vladas Sidoravicius, He Song, Longmin Wang, Kainan Xiang investigated properties of $\USF$ for a random walk they called ``biased random walk.'' They studied the network $\Z^d$ with assignment of conductances $c_\lambda(e) = \lambda^{-|e|},$ where $|e|$ is the graph-distance from the origin to the edge $e$ and $\lambda \in (0, 1)$ is a fixed parameter. The network considered in this paper and $\Gamma_d(\lambda)$ both have conductances which are far from uniform. However, the two assignments of conduces give rise to very different random walks, and very different properties of the uniform spanning forest.

		In \cite{SSSWX:BNRW} the $\NRW$ has a drift away from the origin, while in our case it has a uniform drift towards the right. Our assignments of conductances make the network transitive (indeed, translation-invariant), while theirs is not. They exploited this fact by noticing that the axes play a special r\^{o}le. This also leads to quite different properties of the uniform spanning forest. For example, the $\USF$ in \cite{SSSWX:BNRW} has $2^d$ trees if $d = 2,3,$ while in our case this number is one. Their methods and ours are also different: they used ``spectral radius'' and ``speed'' of the random walk, while we will not mention these. Another display of the difference in nature of the results in their paper when compared to this paper is when they count the number of intersections of random walk paths. They estimated these intersections by using a local limit theorem. Typically, local limit theorems have error terms of a polynomial decay in $\norm{z}.$ We will show that the Green's function in $\Gamma_d(\lambda)$ has exponential decay in some directions. Thus, a local limit theorem did not prove useful for us. Finally, another contrast in the two assignment of conductances. The fact that $c_\lambda(e)$ is bounded from below allowed them to prove one-endedness using the isoperimetric condition \cite[Theorem 10.43]{LyPe:PTN}. Our conductances decay exponentially, and so we cannot use the isoperimetric condition for one-endedness (or any other known technique so far); a new technique was derived from scratch to prove it.

	\section{Green's function bounds}\label{Section: Greens function}
	\subsection{Basic statistics of $\Gamma_d(\lambda)$}
	We will often write $z \in \Z^{d + 1} = \Z \times \Z^d$ as $z = (n, \avec{x}{1}{d}) = (n, x)$ and we will refer to $n$ as the ``drifted coordinate'' or the ``zeroth coordinate.'' A unit vector in $\Z^{d + 1}$ has one coordinate equal to 1 or -1 and all the rest equal to zero; for simplicity, $u$ and $u_i$ will denote unit vectors in $\Z^{d + 1}$ and $u_i$ is such that its $i$th coordinate is non zero ($0 \leq i \leq d$). The functions $\pr{i}:\Z^{d + 1} \to \Z,$ for $i = 0, 1, \ldots, d,$ are the projections and are defined in the obvious way. We will call $\Z u_i$ the ``$i$th factor'' as well as ``$i$th axis'' ($i = 0, \ldots, d$).  Let $\lambda > 0,$ fixed throughout the rest of this document.  Recall $\Gamma_d(\lambda)$ is the graph $\Z^{d + 1}$ with conductances (\ref{Equation: Conductances of the network of interest}). From these conductances, one can get the transition kernel for network random walk. The probability that the next step should be $u \in \Z^{d + 1},$ a unit vector with our conventions, is given by
	\begin{equation}\textstyle\label{Equation: Transition density for network random walk}
		p(u) = \left\{\begin{array}{lcl}
			(2d + 1 + e^\lambda)^{-1} &\text{ if } &u = (-1, 0) \text{ or } u = u_i, i = 1, \ldots, d,\\
			e^\lambda (2d + 1 + e^\lambda)^{-1} &\text{ if } &u = (1, 0).
		\end{array}\right.
	\end{equation}
	The $\frac{1}{2}$-lazy version of the above transition density is $p^L(0) = \frac{1}{2},$ $p^L(u) = \frac{1}{2} p(u)$ for unit $u.$

	Denote by $\zeta \cdot \varepsilon_t$ the measure on $\Z^d$ such that it has total mass equal to $\zeta$ at the vertex $t.$ Notice that if $Y = \left( \supvec{Y}{(0)}{(d)} \right)$ is a random element of $\Z^{d + 1}$ with distribution (\ref{Equation: Transition density for network random walk}), then ($i = 1, \ldots, d$)
	\begin{equation}\textstyle\label{Equation: Distribution of the step lengths of the network}
		Y^{(0)} \sim (2d + 1 + e^\lambda)^{-1} \Big( \varepsilon_{-1} + 2d \cdot \varepsilon_0 + e^\lambda \cdot \varepsilon_1 \Big), \quad Y^{(i)} \sim (2d + 1 + e^\lambda)^{-1} \Big( \varepsilon_{-1} + (2d - 1 + e^\lambda) \cdot \varepsilon_0 + \varepsilon_1),
	\end{equation}
	which immediately allows one to conclude that
	\begin{equation}\textstyle\label{Equation: First and second moment of the network}
		\Expectation{Y} = (a, 0) = \left( \frac{e^\lambda - 1}{2d + 1 + e^\lambda}, 0 \right),
		\quad
		\Covariance{Y} =
		\left[\begin{matrix}
			\sigma_0^2 &0 \\0 &\sigma^2 I_d
		\end{matrix}\right] =
		\left[\begin{matrix}
			\frac{e^\lambda + 1}{2d + 1 + e^\lambda}&0\\
			0 &\frac{2}{2d + 1 + e^\lambda} I_d
		\end{matrix}\right],
	\end{equation}
	where $I_d$ is the identity on $\R^d.$ From the Central Limit Theorem, it follows immediately that if $(Y_n)$ is a family of independent random elements following (\ref{Equation: Transition density for network random walk}),
	\begin{displaymath}\textstyle
		\frac{1}{\sqrt{n}} \sum\limits_{j = 1}^n (Y_j - (a, 0)) \xrightarrow[n \to \infty]{\text{weakly}} \distNorm{0}{ \left[\begin{matrix} \sigma_0^2 &0\\0 &\sigma^2 I_d\end{matrix}\right]}.
	\end{displaymath}
	The Fourier transform of $Y$ is
	\begin{equation}\textstyle\label{Equation: Fourier transform of the step-lengths}\textstyle
		\varphi_Y(h) = \Expectation{e^{i \prodEsc{h}{Y}}} = \frac{e^\lambda e^{i h_0} + e^{-i h_0}}{2d + 1 + e^\lambda} + \frac{2}{2d + 1 + e^\lambda} \sum\limits_{j = 1}^d \cos h_j,
	\end{equation}
	where $h = (\avec{h}{0}{d}) \in [-\pi, \pi]^{d + 1}.$ Observe also $\varphi_Y(h) = 1$ if and only if $1 = e^{i h_0} = \cos h_j$ for $j = 0, \ldots, d,$ and this is the same as saying that $h_0 = \ldots = h_d = 0.$

	\subsection{Continuous time network random walk}\label{Section: CTNRW}
	Recall the definition of the continuous time network random walk. In $\Gamma_d(\lambda),$ we can write
	\begin{equation}\textstyle\label{Equation: Definition of the CTNRW}
		\tilde{S}_t = \left( \tilde{B}_t, \tilde{X}_t \right), \quad (t \geq 0).
	\end{equation}
	The following properties come from well-known facts regarding Poisson processes or are otherwise given a reference. Set $p = \Probability{Y = u} = (2d + 1 + e^\lambda)^{-1}.$
	\begin{enumerate}[leftmargin=1.15\parindent]
		\item $\tilde{B}$ and $\tilde{X}$ are independent continuous-time Markov processes.
	
		\item $\tilde{X}_t = X_{M_t},$ where $X$ is a standard random walk in $\Z^d$ and $M = (M_t)_{t \geq 0}$ is a Poisson process in $\Z_+$ with constant intensity equal to $2dp = \frac{2d}{2d + 1 + e^\lambda}.$
	
		\item $\tilde{B}_t = B_{L_t},$ where $B$ is a biased random walk in $\Z$ with jumps to $\pm 1,$ $L = (L_t)_{t \geq 0}$ is a Poisson process in $\Z_+$ with constant intensity equal to $\sigma_0^2 = (1 + e^\lambda)p = \frac{1 + e^\lambda}{2d + 1 + e^\lambda},$ and the probability of $B$ having a positive jump is $\frac{e^\lambda}{1 + e^\lambda}.$
	
		\item Let $q^X$ be the \emph{transition probability} of $\tilde{X},$ that is, $q^X(t, y) = q_t^X(y) = \Probability[^0]{\tilde{X}_t = y}.$ Then \cite[Theorem 5.17]{Barlow:RWHK} provides upper bounds
	\end{enumerate}
	\begin{equation}\textstyle\label{Equation: Upper exponential bounds for the CTSRW kernel}
		q^X(t, y) \leq
		\begin{cases}
			\frac{2}{d} \exp\left( -\frac{\norm{y}^2}{2dp e^2 t} \right) &\text{ if } \norm{y} \leq 2dpe t, \\
			\frac{1}{2d} \exp\left( -2dpt - \norm{y} \log\left( \frac{ \norm{y} }{ 2dpe t } \right) \right) &\text{ if } \norm{y} \geq 2dpe t;
		\end{cases}
	\end{equation}
	\hspace{\parindent}and lower bounds, valid for $\norm{y} \geq 2dp t > 0,$
	\begin{equation}\textstyle\label{Equation: Lower exponential bounds for the CTSRW kernel}
		q^X(t, y) \geq c \exp\left(-c' \norm{y} \left[1 + \log \left( \frac{\norm{y}}{2dp t} \right) \right] \right).
	\end{equation}
	
	\begin{enumerate}[leftmargin=1.15\parindent]\setcounter{enumi}{4}
		\item By Theorem 6.19 bearing in mind Definition 3.28, Corollary 3.30 and Definition 5.19 of \cite{Barlow:RWHK} there exists four constants $c_i > 0$ ($i = 1, \ldots, 4$) such that for all $(t, y)$ for which $2dp t > \norm{y},$
	\end{enumerate}
	\begin{equation}\textstyle\label{Equation: Gaussian bounds for the CTSRW kernel}
		\begin{split}
			q^X_t(y) &\leq c_1 t^{-\frac{d}{2}} e^{-c_2 \frac{\norm{y}^2}{t}}, \\
			q^X_t(y) &\geq c_3 t^{-\frac{d}{2}} e^{-c_4 \frac{\norm{y}^2}{t}}.
		\end{split}
	\end{equation}

	\begin{enumerate}[leftmargin=1.15\parindent]\setcounter{enumi}{5}
		\item\label{Item: Result on Martins book} There exists a pair of universal constants $c, c' > 0$ such that for all $r, T \geq 1,$ $\Probability{\sup\limits_{0 \leq s \leq T} \norm{\tilde{X}_s} \geq r} \leq c e^{-c' \frac{r^2}{T}}.$ This follows from Theorem 4.33 of \cite{Barlow:RWHK} applied to the standard random walk (so $\beta = 2$), you will need to recall Definitions 4.14 and 4.18, and use Lemma 4.20, which holds with $\alpha = d$ and $\beta = 2$ for the standard random walk on $\Z^d.$
	\end{enumerate}

	We conclude this section with the following easy remark.
	\begin{Remark}\label{Remark: The norm is proportional to the maximum entry}
		Let $c > 0.$
		\begin{enumerate}
			\item If $z = (n, x)$ satisfies $|n| \leq c \norm{x},$ then $\norm{x} \leq \norm{z} \leq \sqrt{1 + c^2} \norm{x}.$
			\item If $z = (n, x)$ satisfies $\norm{x} \leq c |n|,$ then $|n| \leq \norm{z} \leq \sqrt{1 + c^2} |n|.$
			\item There exists $c' > 0$ such that for all $z \in \Z^{d + 1},$ $e^{-c \norm{z}} \leq c' \norm{z}^{-\frac{d}{2}}.$
		\end{enumerate}
	\end{Remark}

	\subsection{Some estimates of sums}\label{Section: Estimates of sums}
	We collect together here some estimates of sums which will be used later. All proofs can be found in \cite{MaDib:USFwithDrift}, sect. 2.3.

	\begin{Proposition}\label{Proposition: Asymptotic estimate on integrals of radial functions with controlled monotonicity}
		Let $L \in \N$ and $f:[L, \infty) \to \R_+^*$ be a monotone function such for a pair $a, b > 0,$ we have $a f(j + 1) \leq f(j) \leq b f(j + 1)$ for integers $j \geq L.$ Let $r > L + 1.$ Then, $\sum\limits_{L \leq \norm{x} \leq r} f(\norm{x}) \asymp \integral[L \leq \norm{x} \leq r] dx\ f(\norm{x}) = \sigma_d \integral[L]^r dt\ f(t) t^{d - 1},$ where $\sigma_d$ is the ``surface area'' of the $d$-dimensional sphere and any implicit constant may depend solely on dimension and on the pair $(a, b).$
	\end{Proposition}

	\begin{Proposition}\label{Proposition: Asymptotic estimates on Eulerian sums}
		Let $a \geq 0,$ $b, \gamma > 0$ and $r > 2.$
		\begin{enumerate}
			\item If $a > 0,$ $\sum\limits_{\norm{x} \leq r} \norm{x}^a e^{-\gamma \norm{x}^b} \asymp \min \left( r^b, \gamma^{-1} \right)^{\frac{a + d}{b}}.$
			\item If $a  = 0,$ $\sum\limits_{0 < \norm{x} \leq r} e^{-\gamma \norm{x}^b} \asymp \min\left( r^b, \gamma^{-1} \right)^{\frac{d}{b}};$ thus, if $r^{-b} \leq \gamma \leq 1,$ $\sum\limits_{\norm{x} \leq r} e^{-\gamma \norm{x}^b} \asymp \gamma^{-\frac{d}{b}}.$
		\end{enumerate}
		(Any implicit constant may depend on $a, b$ and dimension but not on $\gamma$ or on $r.$)
	\end{Proposition}

	\begin{Proposition}\label{Proposition: Sums of Greens function, ver 1}
		Let $a, b > 0.$ There exists a pair of constants $c_1, c_2 > 0$ such that for all $t \in \Z^d$ and all $m \in \N,$
		\begin{enumerate}
			\item if $\norm{t} > bm,$ $e^{-a \norm{t}} \leq \sum\limits_{\norm{x} \leq bm} e^{-a \left( \frac{\norm{x}^2}{m} + \norm{x - t} \right)} \leq	c_1 e^{-c_2 \norm{t}};$

			\item if $\norm{t} \leq bm,$ $	e^{-a \frac{\norm{t}^2}{m}} \leq \sum\limits_{\norm{x} \leq bm} e^{-a \left( \frac{\norm{x}^2}{m} + \norm{x - t} \right)} \leq c_1 e^{-c_2 \frac{\norm{t}^2}{m}}.$
		\end{enumerate}
	\end{Proposition}

	\begin{Proposition}\label{Proposition: Estimates in the p series}
		For every $\alpha > 1$ and all integers $m \geq 2,$ $\frac{1}{\alpha - 1} m^{-(\alpha - 1)} \leq \sum\limits_{k \geq m} k^{-\alpha} \leq \frac{2^{\alpha - 1}}{\alpha - 1} m^{-(\alpha - 1)}.$
	\end{Proposition}

	\begin{Proposition}\label{Proposition: Estimates in the p series with an exponential factor}
		For every $a \in \R$ and $b \geq 1,$ there exists a constant $c = c(a, b) > 0$ such that for all $\gamma > 1,$ $	\sum\limits_{1 \leq j \leq b\gamma} j^a e^{-\frac{\gamma}{j}} \leq c \gamma^{a + 1}.$
	\end{Proposition}

	\begin{Proposition}\label{Proposition: Sums of Greens function, ver 2}
		There exists a universal constant $c > 0$ such that for all systems $(a, p, q, r, t)$ with $r > 2,$ $a, p, q > 0$ and $t \in \Z^d,$
		\begin{displaymath}\textstyle
			\sum\limits_{\norm{x} \leq r} e^{-a \left( \frac{\norm{x}^2}{p} + \frac{\norm{x - t}^2}{q} \right)}
			\leq \begin{cases}
				\ds c \min \left( r^2, \frac{p}{a} \right) ^\frac{d}{2}e^{-\frac{a}{4} \frac{\norm{t}^2}{q}}&\text{ if } \norm{t} > 2r \\
				\ds c \left[ \min \left( r^2, \frac{p}{a} \right)^\frac{d}{2} e^{-\frac{a}{4} \frac{\norm{t}^2}{q}} + \min \left( \frac{\norm{t}^2}{4}, \frac{q}{a} \right)^{\frac{d}{2}} e^{-\frac{a}{4} \frac{\norm{t}^2}{p}} \right] &\text{ if } \norm{t} \leq 2r.
			\end{cases}
		\end{displaymath}
	\end{Proposition}

	\subsection{Some properties of (network) random walk}
	These properties are well-known and are included for completeness of the exposition. We denote by $S = (S_n)$ the network random walk of a network $\Gamma$  and $G$ will denote its Green's function.

	\begin{Proposition}\label{Proposition: The Greens function is a constant multiple the probability of reaching a point}
		Let $\Gamma$ be any network. For any two vertices $o$ and $v$ of $\Gamma,$ $G(o, v) = \integral[0]^\infty dt\ q_t^S(o, v).$ If $\Gamma$ is transitive and $v$ is any vertex, $G(o, v) = G(o, o) \Probability[^o]{\tau_v < \infty},$ where $\tau_v = \tau_v \left( \tilde{S} \right) = \inf\left\{t \geq 0;\ \tilde{S}_t = v \right\}.$
	\end{Proposition}

	\begin{Proposition}\label{Proposition: Probability for the existence of splitting levels}
		Consider a random walk $S$ on $\Z$ with step-lengths distributed as $q \cdot \varepsilon_{-1} + r \cdot \varepsilon_0 + p \cdot \varepsilon_1.$ Assume $1 > p + q > p > q > 0.$ Consider the events $\strC_u = \{\text{for exactly one } n, S_n = u\} = \bigcup\limits_{n \in \Z_+} \bigcap\limits_{m \neq n} \{S_n = u, S_m \neq u\}$ and $\strE_u = \{\text{for at least one } n, S_n = u\} = \bigcup\limits_{n \in \Z_+} \{S_n = u\}.$ Then, for $u \geq 0$ and $v < 0,$
		\begin{enumerate}
			\item $\Probability[^0]{\strC_u} = \Probability[^0]{\tau_0^+ = \infty} = p - q,$
			\item $\Probability[^0]{\bigcup\limits_{u = 0}^{n^2} \strC_u} \geq 1 - (1 - (p - q))^n - n \left(\frac{q}{p} \right)^n,$
			\item $\Probability[^0]{\strE_v} = \left( \frac{p}{q} \right)^v.$
		\end{enumerate}
	\end{Proposition}
	\begin{proof}
		Items (a) and (c) are well-known. For convenience of the reader, we provide details of item (b).

		Consider $\strD_u$ the complement of $\strC_u.$ It suffices to show $\Probability[^0]{\bigcap\limits_{u = 0}^n \strD_{un}} \leq (1 - (p - q))^n + n \left(\frac{q}{p} \right)^n.$ With respect to $\mathds{P}^0,$ the events $\strD_u$ ($u \geq 0$) are a.s. the events where the state $u$ is visited twice or more. Introduce $\theta_t$ defined by ``the time of first return to $t$'': $\theta_t = \inf\{k >\tau_t \mid S_k = t\} \quad (t \in \Z_+).$ It is clear that $(\tau_{un})_{u = 1, \ldots, n}$ is strictly increasing, that $\tau_t < \theta_t$ for all $t \in \Z_+,$ and also that the complement of the event $\bigcap\limits_{j = 1}^n \{\theta_{(j - 1) n} < \tau_{j n}\}$ is contained in the event $\bigcup\limits_{j = 1}^n \{\tau_{jn} < \theta_{(j - 1)n}\}.$ Therefore, $\Probability[^0]{\bigcap\limits_{u = 0}^n \strD_{un}} \leq \Probability[^0]{\theta_0 < \tau_n < \theta_n < \ldots < \tau_{n^2} < \theta_{n^2} < \infty} + \sum\limits_{j = 1}^n \Probability[^0]{\tau_{jn} < \theta_{(j - 1)n} < \infty}.$ By the strong Markov property, $\Probability[^0]{\tau_{jn} < \theta_{(j - 1)n} < \infty}=\Probability[^0]{\strE_{-n}} = \left( \frac{q}{p} \right)^n.$ Set $\scrE_n = \scrF_{\tau_n}^S = \sigma(S_k; 0 \leq k \leq \tau_n)$ and
		\begin{displaymath}\textstyle
			\gamma_k \mathop{=}\limits^{\mathrm{def.}} \Probability[^0]{\theta_0 < \tau_n < \theta_n < \ldots < \tau_{k n} < \theta_{k n} < \infty}.
		\end{displaymath}
		Observe that the strong Markov property implies
		\begin{align*}
			\gamma_n &= \Probability[^0]{\theta_0 < \tau_n < \theta_n < \ldots < \tau_{n^2} < \theta_{n^2} < \infty} \\
			&= \Expectation[^0]{\indic{\{\theta_0 < \tau_n < \infty\}} \Expectation[^0]{ \indic{\{\theta_n < \ldots < \tau_{n^2} < \theta_{n^2} < \infty\}}\middle|\ \scrE_n}} \\
			&= \Probability[^0]{\theta_0 < \tau_n < \infty} \gamma_{n - 1} \leq \Probability[^0]{\tau_0^+ < \infty} \gamma_{n - 1},
		\end{align*}
		which leads at once to $\gamma_n \leq (1 - (p - q))^n$ since $\gamma_1 \leq 1 - (p - q).$
	\end{proof}

	\subsection{A hitting time for the continuous time network random walk}
	Again, we compile some properties of certain random times associated with the continuous time random walk (\ref{Equation: Definition of the CTNRW}). These properties are arguably well-known and are included for completeness of the exposition. In what follows, we use (\ref{Equation: First and second moment of the network}) $p = (2d + 1 + e^\lambda)^{-1},$ $a = (e^\lambda - 1)p$ and $\sigma_0^2 = (e^\lambda + 1)p.$
	Introduce
	\begin{displaymath}\textstyle
		\tau_n^0 = \inf \left\{ t \geq 0; \tilde{B}_t = n \right\} \quad (n \in \N).
	\end{displaymath}
	Then, $\tau_n^0$ is an almost surely finite stopping time, it is independent of the process $\tilde{X}$  and by successive applications of the strong Markov property at the times $\tau_j^0$ ($1 \leq j \leq n - 1$), it is immediately seen that the law of $\tau_n^0$ is the $n$th convolution power of the law of $\tau_1^0.$ Let $\zeta$ denote the Laplace transform of $\tau_1^0,$ that is, $\zeta(t) = \Expectation[^0]{e^{t \tau_1^0}}.$ The following proposition follows from elementary results and its proof is therefore omitted, see \cite{MaDib:USFwithDrift}, propositions (2.5.1) and (2.5.2).
	\begin{Proposition}
		\begin{enumerate}
			\item We have $\Expectation[^0]{\tau_1^0} = \frac{1}{a}.$

			\item Then $\zeta(t) < \infty$ for $t < (e^{\frac{\lambda}{2}} - 1)^2 p$ and $\zeta(t) = \infty$ for $t > (e^{\frac{\lambda}{2}} - 1)^2 p.$ We have the following bounds
			\begin{displaymath}\textstyle
				\zeta(t) \leq
				\begin{cases}
					\frac{2e^\lambda p}{2e^\lambda p - t} &\text{ for } -\infty < t \leq 0, \vspace{0.1cm}\\
					\frac{2e^\lambda p}{\sigma_0^2 - t} &\text{ for } 0 < t < (e^{\frac{\lambda}{2}} - 1)^2 p.
				\end{cases}
			\end{displaymath}
		\end{enumerate}
	\end{Proposition}
	
	\begin{Proposition}\label{Proposition: Exponential decay of the tails of the hitting time of a hyperplane}
		\begin{enumerate}
			\item For $\alpha \leq \frac{n}{2e^\lambda p},$ $\Probability[^0]{\tau_n^0 \leq \alpha} \leq \exp\left( -2e^\lambda p \alpha + n \log \left( \frac{2 e^{\lambda + 1} p \alpha}{n} \right) \right).$

			\item For $0 < \beta < \frac{1}{2e^\lambda p}$ there exists a constant $u(\beta) > 0$ such that for all $n,$ $\Probability[^0]{\tau_n^0 \leq \beta n} \leq e^{-u(\beta) n}.$

			\item There exists two constants $T_0, K_0 > 0$ such that $\Probability[^0]{\tau_n^0 \geq \alpha} \leq e^{-T_0(\alpha - \frac{n}{a})}$ valid for $\alpha \geq (K_0 + \frac{1}{a}) n.$
		\end{enumerate}
	\end{Proposition}
	\begin{proof}
		The bound in (c) is a Corollary of Theorem 1 and Lemma 5, in Chapter III of \cite{Petrov:SIRV}, and observe we can apply these results since the previous proposition shows the Laplace transform of $\tau_1^0$ to be absolutely convergent on an interval around zero. We now prove the bound for $\Probability[^0]{\tau_n^0 \leq \alpha}.$ Markov's inequality shows, for $t < 0,$ $\Probability[^0]{\tau_n^0 \leq \alpha} = \Probability[^0]{e^{t \tau_n^0} \geq e^{t \alpha}} \leq e^{-t \alpha} \zeta(t)^n,$ where $\zeta$ is the Laplace transform of $\tau_1^0.$ From the previous proposition, we know $\zeta(t) \leq \frac{2 e^\lambda p}{2 e^\lambda p - t}$ for $t \leq 0,$ and so $\Probability[^0]{\tau_n^0 \leq \alpha} \leq e^{-t \alpha + n \log \left( \frac{2 e^\lambda p}{2 e^\lambda p - t} \right) }$ for $t < 0.$ Minimising gives the minimiser $t = 2 e^\lambda p - \frac{n}{\alpha},$ which is valid only when $t < 0,$ that is, only when $\alpha < \frac{n}{2e^\lambda p},$ the bound follows readily upon evaluation at the minimiser. When $\alpha = \frac{n}{2e^\lambda p},$ we reach the trivial bound $\Probability[^0]{\tau_n^0 \leq \alpha} \leq 1,$ and so it is also valid in this case. Finally, the substitution $\alpha = \beta n$ implies $\Probability[^0]{\tau_n^0 \leq \beta n} \leq e^{-u(\beta) n},$ where $u(\beta) = 2e^\lambda p \beta - \log(2e^{\lambda + 1} p \beta).$ A standard minimisation shows $u(\beta) > u\left( \frac{1}{2e^\lambda p} \right) = 0$ for $0 < \beta < \frac{1}{2e^\lambda p}.$
	\end{proof}

	\begin{Proposition}\label{Proposition: Exponential bound on the visit time of a point in the vertical hyperplane}
		Recall from (\ref{Proposition: The Greens function is a constant multiple the probability of reaching a point}) the definition of $\tau_z$ for $z = (n, x) \in \Z^{d + 1}.$ There exists $c_1, c_2 > 0$ such that, for all $x \in \Z^d,$ $\Probability[^0]{\tau_{(0, x)} < \infty} \leq c_1 e^{-c_2 \norm{x}}.$
	\end{Proposition}
	\begin{proof}
		Let $r = r_x$ be the smallest integer $k \geq \frac{\norm{x}}{2}.$ Then $\Probability[^0]{\tau_{(0, x)} < \infty} \leq \Probability[^0]{\tau_{(0, x)} \leq \tau_r^0 \leq c r} + \Probability[^0]{\tau_{(0, x)} \leq \tau_r^0, \tau_r^0 > c r} + \Probability[^0]{\tau_r^0 < \tau_{(0, x)} < \infty},$ where $c \geq K_0 + \frac{1}{a}$ is chosen so that $c r_x \geq (K_0 + \frac{1}{a}) r_x$ for all $x \in \Z^d,$ $x \neq 0$ and $K_0$ is as in (\ref{Proposition: Exponential decay of the tails of the hitting time of a hyperplane}). Then, recalling $\tilde{S} = \left( \tilde{B}, \tilde{X} \right)$ and item (\ref{Item: Result on Martins book}) in \S \ref{Section: CTNRW},
		\begin{displaymath}\textstyle
			\Probability[^0]{\tau_{(0, x)} \leq \tau_r^0 \leq c r} \leq \Probability[^0]{\sup\limits_{0 \leq s \leq c r} \norm{\tilde{X}_s} \geq \norm{x}} \leq c' e^{-c'' \frac{\norm{x}^2}{c \norm{x}}} = c' e^{-\frac{c''}{c} \norm{x}},
		\end{displaymath}
		for two constants $c', c'' > 0.$ From (\ref{Proposition: Exponential decay of the tails of the hitting time of a hyperplane}), $\Probability[^0]{\tau_{(0, x)} \leq \tau_r^0, \tau_r^0 > c r} \leq e^{-c' r},$ for $c' = T_0(c - \frac{1}{a}).$ Applying the strong Markov property at the time $\tau_r^0,$ we reach that
		\begin{align*}
			\Probability[^0]{\tau_r^0 < \tau_{(0, x)} < \infty} &\leq \Expectation[^0]{\Probability[^{\tilde{S}_{\tau_r^0}}]{\tau_{(0, x)} < \infty}} \leq \Probability[^0]{\tilde{B} \text{ ever reaches } -r} \\
			&= \Probability[^0]{B \text{ ever reaches } -r} \leq e^{-c' \norm{x}},
		\end{align*}
		for some $c' > 0.$
	\end{proof}

	\begin{Proposition}\label{Proposition: The Greens decomposed in its two main coordinates as a sum of an integral}
		Recall from (\ref{Proposition: The Greens function is a constant multiple the probability of reaching a point}) the definition of $\tau_z$ for $z \in \Z^{d + 1}.$ If $z = (n, x),$ we have $G(0, z) = G(0, 0) \sum\limits_{y \in \Z^d} \Probability[^0]{\tau_{(0, x - y)} < \infty} \integral[0]^\infty q^X_t(y) \Probability{\tau_n^0 \in dt}.$
	\end{Proposition}
	\begin{proof}
		By (\ref{Proposition: The Greens function is a constant multiple the probability of reaching a point}), we have $G(0, z) = G(0, 0) \Probability[^0]{\tau_z < \infty}.$ Apply the strong Markov property at the time $\tau_n^0$ and use that $\tau_n^0$ and $\tilde{X}$ are mutually independent,
		\begin{align*}\textstyle
			\Probability[^0]{\tau_z < \infty} &\textstyle= \Expectation[^0]{ \Probability[^{ \tilde{S}_{\tau_n^0}}]{\tau_z < \infty} } = \sum\limits_{y \in \Z^d} \Probability[^{(n, y)}]{\tau_{(n, x)} < \infty} \Probability[^0]{ \tilde{S}_{\tau_n^0} = (n, y)} \\
			&\textstyle= \sum\limits_{y \in \Z^d} \Probability[^0]{\tau_{(0, x - y)} < \infty} \integral[0]^\infty \Probability[^0]{ \tilde{X}_t = y} \Probability{\tau_n^0 \in dt}.
		\end{align*}
		Since $q^X_t(y) = \Probability[^0]{ \tilde{X}_t = y},$ we reach the desired conclusion.
	\end{proof}

	\subsection{Statement and proof of the Theorem}
	Consider the Green's function $G$ of $\Gamma_d(\lambda).$ By transience, the probability $p_0 = \Probability[^0]{\tau_0 < \infty}$ that the network random walk started at zero returns to zero is (strictly) less than one. Thus, the number of visits to zero, started at zero, has finite expectation since this number follows a geometric law with probability of success $p_0,$ in other words, $G(0, 0)$ is some finite number, depending only on $\lambda$ (as such, regarded as universal). The next theorem deals with the bounds $G(0, z)$ for $z \neq 0.$
	\begin{Theorem}\label{Theorem: Greens function bounds}
		There exist four constants $c_i$ ($i = 1, \ldots, 4$) such that for all vertices $z = (n, x) \neq 0$ of the network $\Gamma_d(\lambda),$
		\begin{displaymath}\textstyle
			G(0, z) \leq c_1 \begin{cases} e^{-c_2 \norm{z}} &\text{ if } \norm{x} > n, \\
			e^{-c_2 \frac{\norm{x}^2}{n}} \norm{z}^{-\frac{d}{2}} &\text{ if } \norm{x} \leq n,
			\end{cases}
		\end{displaymath}
		and
		\begin{displaymath}\textstyle
			G(0, z) \geq c_3 \begin{cases} e^{-c_4 \norm{z}} &\text{ if } \norm{x} > n, \\
			e^{-c_4 \frac{\norm{x}^2}{n}} \norm{z}^{-\frac{d}{2}} &\text{ if } \norm{x} \leq n,
			\end{cases}
		\end{displaymath}
	\end{Theorem}
	\begin{proof}
		The proof is rather lengthy and we divide it into several steps.
	
		\vspace{0.15cm}
		\noindent \textsc{Proof of the upper bounds.}
		
		We begin with proving the upper bound when $z \in \Z_- \times \Z^d. $ We start with the following simple observation: bearing in mind Remark (\ref{Remark: The norm is proportional to the maximum entry}), if $z = (n, x) \in \Z_- \times \Z^d,$ the upper bound will be established if we prove $G(0, z) \leq c e^{-c' \max(|n|, \norm{x})}$ for a pair of constants $c, c' > 0$ independent of $z.$ 
	
		We analyse separately the cases $\norm{x} < |n|$ and $\norm{x} \geq |n|.$
	
		Assume first $z \in \Z_- \times \Z^d$ is such that $\norm{x} < |n|,$ notice $n \leq 0.$ We claim there exists a pair of constants $c, c' > 0$ (independent of $z = (n, x)$) such that $G(0, z) \leq c e^{-c' |n|}.$ To see this, recall (\ref{Proposition: The Greens function is a constant multiple the probability of reaching a point}) and and (\ref{Proposition: Probability for the existence of splitting levels}), so that
		\begin{displaymath}\textstyle
			G(0, z) = G(0, 0) \Probability[^0]{\tau_z < \infty} \leq G(0, 0) \Probability[^0]{B \text{ ever reaches level } n} = G(0, 0) e^{c' n},
		\end{displaymath}
		where $\tau_z = \tau_z\left( \tilde{S} \right)$ is the visit time of the point $z$ by the continuous time random walk $\tilde{S}$ (\ref{Equation: Definition of the CTNRW}), proving the claim.
	
		Assume now $z \in \Z_- \times \Z^d$ satisfies $\norm{x} \geq |n|,$ notice $n \leq 0.$ We have that $G(0, z) = \frac{\mu(z)}{\mu(0)} G(0, -z),$ and since $\mu(z) \asymp e^{\lambda n},$ we have the following preliminary result.
		\begin{Lemma}\label{Lemma: Green function upper bounds for negative and small n}
			If the Green's function bounds hold on the region of $z \in \Z_+ \times \Z^d$ such that $0 \leq n \leq \norm{x},$ then they also hold on the region of points $z \in \Z_- \times \Z^d$ such that $|n| \leq \norm{x}.$
		\end{Lemma}
		Observe that (\ref{Proposition: Exponential bound on the visit time of a point in the vertical hyperplane}) shows that the Green's function bounds hold on the set $\{0\} \times \Z^d.$ Thus, with the previous lemma in mind, the case $z \in \Z_- \times \Z^d$ has been established.
	
		Let us focus now on the case $z = (n, x) \in \Z_+^* \times \Z^d.$ By virtue of (\ref{Proposition: Exponential bound on the visit time of a point in the vertical hyperplane}) and (\ref{Proposition: The Greens decomposed in its two main coordinates as a sum of an integral}), we have
		\begin{equation}\textstyle
			G(0, z) \leq c \sum\limits_{y \in \Z^d} e^{-c' \norm{x - y}} \integral[0]^\infty q_t^X(y) \Probability{\tau_n^0 \in dt}, \tag{$\ast$}
		\end{equation}
		where $\tau_n^0$ is the visit time of $n$ by the walk $\tilde{B}$ and $q_t^X$ is the continuous time transition probability of $\tilde{X}.$
		\begin{Lemma}
			There exist two constants $c_1, c_2 > 0,$ independent of $z,$ such that
			\begin{displaymath}\textstyle
				\integral[0]^\infty q^X_t(y) \Probability{\tau_n^0 \in dt} \leq \begin{cases}
					c_1 e^{-c_2n} + n^{-\frac{d}{2}} e^{-c_2 \frac{\norm{y}^2}{n} } &\text{ if } \norm{y} \leq \frac{d}{e^\lambda + 1} n \vspace{0.15cm}\\
					c_1 e^{-c_2 \norm{y}} &\text{ if } \norm{y} > \frac{d}{e^\lambda + 1} n.
				\end{cases}
			\end{displaymath}
		\end{Lemma}
		\begin{proof}[Proof of lemma.]
			Recall $p = (2d + 1 + e^\lambda)^{-1}.$ We have $	\integral[0]^\infty q^X_t(y) \Probability{\tau_n^0 \in dt} = \integral[0]^\infty q^X_{\frac{t}{2dp}}(y) \Probability{2dp\ \tau_n^0 \in dt}.$ We will tackle the two cases.
		
			\textbf{Case 1:} $\norm{y} \leq \frac{d}{e^\lambda + 1} n.$ We split the integral as follows
			\begin{displaymath}\textstyle
				\integral[0]^\infty q^X_{\frac{t}{2dp}}(y) \Probability{2dp\ \tau_n^0 \in dt} = \left[ \integral[0]^{\frac{\norm{y}}{e}} + \integral[\frac{\norm{y}}{e}]^{\frac{d}{e^\lambda + 1} n} + \integral[\frac{d}{e^\lambda + 1} n]^{ 2dp(K_0 + \frac{1}{a}) n } + \integral[ 2dp(K_0 + \frac{1}{a}) n]^\infty \right] q^X_{\frac{t}{2dp}}(y) \Probability{2dp\ \tau_n^0 \in dt},
			\end{displaymath}
			where $K_0$ is as in (\ref{Proposition: Exponential decay of the tails of the hitting time of a hyperplane}). In what follows, we will use the bounds for $\Probability{\tau_n^0 \leq \alpha}$ and $\Probability{\tau_n^0 \leq \beta n}$ of this proposition without further notice. The first two integrals will be bounded using the bounds (\ref{Equation: Upper exponential bounds for the CTSRW kernel}). In the region $0 \leq t \leq \frac{\norm{y}}{e},$ we have $q^X_{\frac{t}{2dp}}(y) \leq \frac{1}{2d} e^{-t - \norm{y} \log \left(\frac{\norm{y}}{et} \right)}$ and the function $t \mapsto -t-\norm{y} \log(\frac{\norm{y}}{et})$ is easily seen to be increasing on $(0, \norm{y}),$ hence
			\begin{align*}\textstyle
				\integral[0]^{\frac{\norm{y}}{e}} q^X_{\frac{t}{2dp}}(y) \Probability{2dp\ \tau_n^0 \in dt} &\textstyle\leq \frac{1}{2d} e^{-\frac{\norm{y}}{e}} \Probability{\tau_n^0 \leq \frac{\norm{y}}{2d p e}} \leq \frac{1}{2d} e^{-\frac{\norm{y}}{e}} \Probability{\tau_n^0 \leq \frac{n}{2 e^{\lambda + 1} p}} \leq e^{-\frac{\norm{y} + n}{e}} \leq e^{-\frac{n}{e}}.
			\end{align*}
			Similarly, in the region $\frac{\norm{y}}{e} \leq t < \infty,$ we have the bound $q^X_{\frac{t}{2dp}}(y) \leq \frac{2}{d} e^{-\frac{\norm{y}^2}{e^2 t}},$ and the term on the right is an increasing function of $t > 0.$ Hence,
			\begin{align*}\textstyle
				\integral[\frac{\norm{y}}{e}]^{\frac{d}{e^\lambda + 1} n} q^X_{\frac{t}{2dp}}(y) \Probability{2dp\ \tau_n^0 \in dt} &\textstyle\leq \frac{2}{d} e^{- \left( \frac{e^{\lambda - 2}}{d} \right) \frac{\norm{y}^2}{n}} \Probability{\tau_n^0 \leq \frac{n}{2p(e^\lambda + 1)}} \leq e^{- \left( \frac{e^{\lambda - 2}}{d} \right) \frac{\norm{y}^2}{n} - cn} \leq e^{-cn}
			\end{align*}
			for some constant $c = u\left( \frac{1}{2p(e^\lambda + 1)} \right) > 0.$ We use now the bounds (\ref{Equation: Gaussian bounds for the CTSRW kernel}). Consider the interval $\frac{d}{e ^\lambda + 1} n \leq t \leq 2dp \left( K_0 + \frac{1}{a} \right) n,$ we have a Gaussian decay $q^X_{\frac{t}{2dp}}(y) \leq c_1 t^{-\frac{d}{2}} e^{-c_2 \frac{\norm{y}^2}{t}}.$ On this interval, $t \asymp n$ and so, there is a pair of constants $c, c' > 0$ such that $\integral[\frac{d}{e^\lambda + 1} n]^{ 2dp (K_0 + \frac{1}{a}) n } q^X_{\frac{t}{2dp}}(y) \Probability{2dp\ \tau_n^0 \in dt} \leq c n^{-\frac{d}{2}} e^{-c' \frac{\norm{y}^2}{n}}.$ Next, for $t \geq 2dp \left( K_0 + \frac{1}{a} \right) n,$ we use the factor $t^{-\frac{d}{2}}$ of the Gaussian estimate for $q^X_{\frac{t}{2dp}}(y)$ to conclude there exists $c > 0$ such that $q^X_{\frac{t}{2dp}}(y) \leq c_1 t^{-\frac{d}{2}} \leq c n^{-\frac{d}{2}}.$ Thus,
			\begin{displaymath}\textstyle
				\integral[ 2dp( K_0 + \frac{1}{a}) n]^\infty q^X_{\frac{t}{2dp}}(y) \Probability{2dp\ \tau_n^0 \in dt} = c n^{-\frac{d}{2}} \Probability{\tau_n^0 \geq \left( K_0 + \frac{1}{a} \right) n} = c n^{-\frac{d}{2}} e^{-T_0 K_0 n} \leq e^{-T_0K_0 n},
			\end{displaymath}
			again by (\ref{Proposition: Exponential decay of the tails of the hitting time of a hyperplane}). Putting all the bounds found, this concludes the case $\norm{y} \leq \frac{d}{e^\lambda + 1} n.$

			\textbf{Case 2:} $\norm{y} > \frac{d}{e ^\lambda + 1} n.$ Here we divide
			\begin{displaymath}\textstyle
				\integral[0]^\infty q^X_{\frac{t}{2dp}}(y) \Probability{2dp\ \tau_n^0 \in dt} = \left[ \integral[0]^{\frac{\norm{y}}{e}} + \integral[\frac{\norm{y}}{e}]^{ 2 \sigma_0^2 (K_0 + \frac{1}{a}) \norm{y} } + \integral[ 2 \sigma_0^2 (K_0 + \frac{1}{a}) \norm{y}]^\infty \right] q^X_{\frac{t}{2dp}}(y) \Probability{2dp\ \tau_n^0 \in dt}.
			\end{displaymath}
			Recall the function $t \mapsto -t - \norm{y} \log \left( \frac{\norm{y}}{e t} \right)$ is increasing on $(0, \norm{y}).$ As in \textbf{Case 1}, we can bound
			\begin{displaymath}\textstyle
				q^X_{\frac{t}{2dp}}(y) \leq \begin{cases}
					e^{-\frac{\norm{y}}{e}} &\text{ if } t \leq \frac{\norm{y}}{e}, \vspace{0.05cm}\\
					\frac{2}{d} e^{-c \norm{y}} &\text{ if } \frac{\norm{y}}{e} \leq t \leq 2 \sigma_0^2 (K_0 + \frac{1}{a}) \norm{y}, \vspace{0.05cm}\\
					c_1 &\text{ if } t \geq 2 \sigma_0^2 (K_0 + \frac{1}{a}) \norm{y},
				\end{cases}
			\end{displaymath}
			where $c_1$ is the constant appearing on the Gaussian bound of $q_t^X$ and $c = \big( 2 \sigma_0^2 e^2 (K_0 + \frac{1}{a}) \big)^{-1}.$ We may assume, should the need arise, $K_0$ is so large so that $2 \sigma_0^2 (K_0 + \frac{1}{a}) \geq 1.$ It is obvious then that the first two integrals decay exponentially in $\norm{y},$ as for the third one, we employ (\ref{Proposition: Exponential decay of the tails of the hitting time of a hyperplane}) to obtain
			\begin{displaymath}\textstyle
				\integral[ 2 \sigma_0^2 (K_0 + \frac{1}{a}) \norm{y}]^\infty q^X_{\frac{t}{2dp}}(y) \Probability{2dp\ \tau_n^0 \in dt} \leq c_1 \Probability{\tau_n^0 \geq \frac{\sigma_0^2 (K_0 + \frac{1}{a}) \norm{y}}{dp}},
			\end{displaymath}
			the condition $\frac{d}{e^\lambda + 1} n < \norm{y}$ guarantees $\Probability{\tau_n^0 \geq \frac{\sigma_0^2 (K_0 + \frac{1}{a}) \norm{y}}{dp}} \leq e^{-T_0 K_0 \frac{e^\lambda + 1}{d} \norm{y}}.$ This completes \textbf{Case 2} and the lemma.
		\end{proof}
	
		With the lemma proven, let us return to the proof of the upper bounds for the Green's function. The reader is reminded that there only remains to prove the bounds for $z = (n, x)$ in $\Z_+^* \times \Z^d.$ The previous lemma and ($\ast$) imply that there exist two constants $c, c' > 0$ such that
		\begin{displaymath}\textstyle
			G(0, z) \leq c\left( \sum\limits_{\norm{y} \leq b n} e^{-c'(\norm{x - y} + n)} + \sum\limits_{\norm{y} \leq b n} n^{-\frac{d}{2}} e^{-c' \left( \norm{x - y} + \frac{\norm{y}^2}{n} \right)} + \sum\limits_{\norm{y} > bn} e^{-c'(\norm{x - y} + \norm{y})} \right),
		\end{displaymath}
		where $b = \frac{d}{e ^\lambda + 1}.$
	
		Bearing in mind (\ref{Remark: The norm is proportional to the maximum entry}), we only need to prove there exists a pair of constants $c, c' > 0$ such that $G(0, z) \leq c e^{-c' \norm{x}}$ for all $\norm{x} > bn$ and that there exists a pair of constants $c, c' > 0$ such that $G(0, z) \leq c n^{-\frac{d}{2}} e^{-c' \frac{\norm{x}^2}{n}}$ for all $\norm{x} \leq bn.$ We consider first a sum of the form $\sum\limits_{\norm{y} \leq bn} e^{-c(\norm{x - y} + \alpha n)}.$ Set $c' = \min \left(c, \frac{c \alpha}{2b} \right),$ then
		\begin{align*}\textstyle
			\sum\limits_{\norm{y} \leq bn} e^{-c(\norm{x - y} + \alpha n)} &\textstyle\leq \sum\limits_{\norm{y} \leq bn} e^{-c'\norm{x - y} - c \alpha n} \leq e^{-c'\norm{x}} \sum\limits_{\norm{y} \leq bn} e^{\frac{c \alpha}{2b} \norm{y} - c \alpha n} \\
			&\textstyle\leq e^{-c' \norm{x}} L n^d e^{-\frac{c \alpha }{2} n} \leq c'' e^{-c' \norm{x} - \frac{c \alpha }{4} n},
		\end{align*}
		in which $L$ is a constant, depending only in dimension, satisfying that $\card{\ball[_{\Z^d}]{0; r}} \leq L r^d$ for all $r \geq 1$ and $\ds c'' = \sup_{n \in \N} Ln^d e^{-\frac{c \alpha}{4} n} < \infty.$ Thus, we have proven the following lemma.
		\begin{Lemma}
			Let $\alpha, \beta, \gamma > 0.$ There exist constants $c, c' > 0$ such that $\sum\limits_{\norm{y} \leq \beta n} e^{-\gamma (\norm{x - y} + \alpha n)} \leq c e^{-c' \norm{(n, x)}},$ for all $(n, x) \in \Z \times \Z^d.$
		\end{Lemma}

		Next, we consider the sum $\sum\limits_{\norm{y} > bn} e^{-c'(\norm{x - y} + \norm{y})}.$ In the case that $\norm{x} \leq bn,$ then we may bound $e^{-c'(\norm{x - y} + \norm{y})} \leq e^{-c'b n} e^{-c' \norm{x - y}},$ thus $\sum\limits_{\norm{y} > bn} e^{-c'(\norm{x - y} + \norm{y})} \leq e^{-c'b n} \sum\limits_{y \in \Z^d} e^{-c' \norm{y}}.$ In the case that $\norm{x} > b|n|,$ we may bound $\norm{x - y} \geq \norm{x} - \norm{y},$ thus $e^{-c'(\norm{x - y} + \norm{y})} \leq e^{-\frac{c'}{2} \norm{x - y} - c'\norm{y}} \leq e^{-\frac{c'}{2} \norm{x} - \frac{c'}{2} \norm{y}}.$ Therefore $\sum\limits_{\norm{y} > bn} e^{-c'(\norm{x - y} + \norm{y})} \leq e^{-\frac{c'}{2} \norm{x}} \sum\limits_{y \in \Z^d} e^{-\frac{c'}{2} \norm{y}}.$ We have established the following lemma.
		\begin{Lemma}
			Let $\alpha, \beta > 0.$ There exists a pair of constants $c, c' > 0$ such that $\sum\limits_{\norm{y} > \beta n} e^{-\alpha(\norm{x - y} + \norm{y})} \leq c e^{-c' \norm{(n, x)}}$ for  every $(n, x) \in \Z \times \Z^d.$
		\end{Lemma}
		
		By virtue of the previous lemmas, it remains to show that
		\begin{displaymath}\textstyle
			\sum\limits_{\norm{y} \leq b n} n^{-\frac{d}{2}} e^{-c' \left( \norm{x - y} + \frac{\norm{y}^2}{n} \right)} \leq
			\begin{cases}
				c e^{-c' \norm{x}} &\text{ if } \norm{x} > bn, \\
				c e^{-c' \frac{\norm{x}^2}{n}} n^{-\frac{d}{2}} &\text{ if } \norm{x} \leq bn.
			\end{cases}
		\end{displaymath}
		
		Assume first that $\norm{x} > bn > 0.$ The estimate (\ref{Proposition: Sums of Greens function, ver 1}) shows at once $\sum\limits_{\norm{y} \leq bn} e^{-c \left(\norm{x - y} + \frac{\norm{y}^2}{n} \right)} \leq c e^{-c' \norm{x}}$ for a pair of constants $c, c' > 0.$ This concludes the case $\norm{x} > bn > 0.$ Assume now that $\norm{x} \leq bn.$ By (\ref{Proposition: Sums of Greens function, ver 1}), we also have $\sum\limits_{\norm{y} \leq b n} n^{-\frac{d}{2}} e^{-c' \left( \norm{x - y} + \frac{\norm{y}^2}{n} \right)} \leq c n^{-\frac{d}{2}} e^{-c' \frac{\norm{x}^2}{n}}$ for (possibly) another pair $c, c' > 0.$ This concludes the case $\norm{x} \leq bn,$ and thus, the proof of the upper bounds by virtue of Lemma (\ref{Lemma: Green function upper bounds for negative and small n}).

		\newpage
		\noindent \textsc{Proof of lower bounds.}
	
		These are much simpler. We have $G(0, z) = \frac{\mu(z)}{\mu(0)} G(0, - z).$ If $z = -(n, x)$ with $n > 0,$ then $\mu(z) \asymp e^{-\lambda n}$ and so $G(0, z) \geq c e^{-\lambda n} G(0, -z) \geq ce^{-\lambda \norm{z}} G(0, -z).$ Therefore, if we prove lower bounds for $n > 0,$ we will reach the lower bounds for $n < 0.$ For the case $n = 0$ observe the following: by (\ref{Proposition: The Greens function is a constant multiple the probability of reaching a point}),
		\begin{displaymath}\textstyle
			G(0, (0, x)) = G(0, 0) \Probability[^0]{\tau_{(0, x)} < \infty} \geq G(0, 0) \Probability[^0]{\tau_{(0, x)} < \infty, \tilde{S}_{\sigma_1} = (1, 0)},
		\end{displaymath}
		where $\sigma_1$ is the time of the first jump of the Poisson process $N$ in the definition of $\tilde{S}.$ Now,
		\begin{displaymath}\textstyle
			\Probability[^0]{\tau_{(0, x)} < \infty, \tilde{S}_{N_{\sigma_1}} = (1, 0)} = e^\lambda p \integral[0]^\infty dt\ e^{-t} \Probability[^{(1, 0)}]{\tau_{(0, x)} < \infty} = e^\lambda p \Probability[^0]{\tau_{(-1, x)} < \infty}.
		\end{displaymath}
		Thus, if the correct bounds hold for $n < 0,$ then $G(0, (0, x)) \geq e^\lambda p G(0, (-1, x)) \geq c e^{-c' \norm{(1, x)}} \norm{(1, x)}^{-\frac{d}{2}},$ and it is obvious that $\norm{(1, x)} \asymp \norm{(0, x)}$ (\ref{Remark: The norm is proportional to the maximum entry}). Whence, the proof of the lower bounds reduces to the case $n > 0.$
	
		When $n > 0,$ we can apply (\ref{Proposition: The Greens decomposed in its two main coordinates as a sum of an integral}). Now,
		\begin{displaymath}\textstyle
			\integral[0]^\infty q_t^X(y) \Probability{\tau_n^0 \in dt} \geq \integral[\frac{2d}{e^\lambda} n]^{4(K_0 + \frac{1}{a})dp\ n} q_{\frac{t}{2dp}}^X(y) \Probability{2dp\ \tau_n^0 \in dt}.
		\end{displaymath}
		Observe that the lower bounds of (\ref{Equation: Lower exponential bounds for the CTSRW kernel}) and (\ref{Equation: Gaussian bounds for the CTSRW kernel}) are of the same type (with possibly different constants) whenever $\norm{y} \asymp t.$ On the interval of integration, $t \asymp n,$ therefore
		\begin{displaymath}\textstyle
			q_{\frac{t}{2dp}}^X(y) \geq \begin{cases}
				c n^{-\frac{d}{2}} e^{-c' \frac{\norm{y}^2}{n}} &\text{ if } n \geq \norm{y}, \\
				c \norm{y}^{-\frac{d}{2}} e^{-c' \norm{y}} &\text{ if } n < \norm{y}.
				\end{cases}
		\end{displaymath}
		Also, (\ref{Proposition: Exponential decay of the tails of the hitting time of a hyperplane}) gives $\Probability{\frac{2d}{e^\lambda} n \leq 2dp\ \tau_n^0 \leq 4(K_0 + \frac{1}{a})dp\ n} \geq c,$ for a constant $c > 0$ independent of $n.$ Therefore,
		\begin{displaymath}\textstyle
			\frac{G(0, z)}{G(0, 0)} \geq c \left( \sum\limits_{n \geq \norm{y}} \Probability[^0]{\tau_{(0, x - y)} < \infty} n^{-\frac{d}{2}} e^{-c' \frac{\norm{y}^2}{n}} + \sum\limits_{n < \norm{y}} \Probability[^0]{\tau_{(0, x - y)} < \infty} \norm{y}^{-\frac{d}{2}} e^{-c' \norm{y}} \right).
		\end{displaymath}
		We discard all the terms except $y = x,$ and reach
		\begin{displaymath}\textstyle
			G(0, z) \geq \begin{cases}
				c G(0, 0) n^{-\frac{d}{2}} e^{-c' \frac{\norm{x}^2}{n}} &\text{ if } n \geq \norm{x} \\
				c G(0, 0) \norm{x}^{-\frac{d}{2}} e^{-c' \norm{x}} &\text{ if } n < \norm{x}.
			\end{cases}
		\end{displaymath}
		This terminates the proof of Theorem (\ref{Theorem: Greens function bounds}).
	\end{proof}	

	\section{The number of trees in $\USF$}\label{Section: Number of trees}
In this section we will establish a number of basic properties of the uniform spanning forest. $\Gamma_d(\lambda)$ has the Liouville property, this implies both measures $\WSF$ and $\FSF$ coincide in this network. Finally, establish when the samples are trees or forests.

	\subsection{Coupling of random walks}\label{Section: Coupling of random walks}
	Consider two vertices $x$ and $y$ in the network $(\strG, \mu).$ We will say that the network random walk started at $x$ can be \textbf{classically coupled} with the network random walk started at $y$ if the following condition holds:
	\begin{description}
		\item[($\mathsf{CC}$)] There exists a probability space $(\Omega, \scrF, \mathds{P}),$ a filtration $(\scrF_n)_{n \in \Z_+}$ in this probability space, two Markov processes with respect to this filtration and defined on this probability space, denoted by $(S_n, \scrF_n)$ and $\left( S_n', \scrF_n \right),$ both of them having for transition probability that defined by the conductances of the network, $S_0 = x,$ $S_0' = y$ and there exists a stopping time $\tau$ relative to the filtration $(\scrF_n)_{n \in \Z_+},$ such that $\tau < \infty$ $\mathds{P}$-a.e. and for those $\omega \in \Omega$ for which $\tau(\omega) < \infty, $ $S_{\tau(\omega)+n}(\omega) = S_{\tau(\omega)+n}'(\omega)$ for any $n \in \N.$
	\end{description}
	When this happens, we will call $\tau$ the \textbf{coupling time} and the condition $\tau < \infty$ $\mathds{P}$-a.e. is expressed by saying that ``$\tau$ is successful.'' Obviously, for a classical coupling to exist, it is necessary that $x$ and $y$ should be connected in $(\strG, \mu),$ but this is far from sufficient.

	\begin{Remark}
		Observe that by definition, a classical coupling is a coupling $T = \left( S, S' \right)$ of the network random walk $S$ started at $x$ and the network random walk $S'$ started at $y$ satisfying strong conditions and it does not necessarily exists; for instance, the standard lattice $\Z^1$ with $x = 0$ and $y= 1,$ any two standard random walks will always be at distance at least one and so, no classical coupling can exist (the ``parity problem'').
	\end{Remark}

	\begin{Proposition}\label{Proposition: Decomposition of RW}
		Denote by $u_i$ ($1 \leq i \leq d$) the canonical unit vectors of $\Z^d.$ Let $\rho$ be a probability measure on $\Z^d$ such that
		\begin{enumerate}
			\item the relation $\rho(x) > 0$ implies $x$ is an integer multiple of some $u_i;$
			\item $\rho(\Z^* u_i) > 0$ for each $i = 1, \ldots, d.$
		\end{enumerate}
		Introduce two probability measures, first $\varphi(i) = \rho(\Z^* u_i) + \frac{\rho(0)}{d}$ ($1 \leq i \leq d$) and then $\nu_i(k) = \frac{\rho(k u_i)}{\varphi(i)}$ for $k \in \Z^*$ and $\nu_i(0) = \frac{\rho(0)}{d\varphi(i)}.$ Suppose $F_j \sim \varphi$ and $X_j^{(i)} \sim \nu_i$ are independent random elements ($1 \leq i \leq d$ and $j \in \Z_+$). Set $Y_j = X_j^{(F_j)} = \sum\limits_{i = 1}^d X_j^{(i)} \indic{\{F_j = i\}}.$ Then, for any $S_0$ independent of the families $X_j^{(i)}$ and $F_j,$ $S_n = S_0 + Y_1 + \ldots + Y_n$ is a random walk on $\Z^d$ started at $S_0$ with transition probability $\rho:$ $\Probability{Y_j = x} = \rho(x)$ ($x \in \Z^d$).
	\end{Proposition}

	\begin{Remark}\label{Remark: Decomposition of RW for lazy version in Gamma lambda}
		When $\rho$ of (\ref{Proposition: Decomposition of RW}) is given by the $\frac{1}{2}$-lazy version of (\ref{Equation: Transition density for network random walk}), the appropriate distributions are $\varphi(0) = (1 + e^\lambda)(2(2d + 1 + e^\lambda))^{-1} + (2(d + 1))^{-1},$ and for $i = 1, \ldots, d,$ $\varphi(i) = (2d + 1 + e^\lambda)^{-1} + (2(d + 1))^{-1}.$ While the $\nu_i$ are given as $	\nu_0(-1) = (2(2d + 1 + e^\lambda) \varphi(0))^{-1},$ $\nu_0(0) = (2(d + 1) \varphi(0))^{-1},$ $\nu_0(1) = e^\lambda (2(2d + 1 + e^\lambda) \varphi(0))^{-1}$ and, for $i = 1, \ldots, d,$ $\nu_i(1) = \nu_i(-1) = (2(2d + 1 + e^\lambda) \varphi(i))^{-1},$ $\nu_i(0) = (2(d + 1) \varphi(i))^{-1}.$
	\end{Remark}

	The following theorem is arguably well-known while a source where it is proven is elusive to find. The reader can find a complete proof in \cite{MaDib:USFwithDrift}, theorem (3.2.3). The proof-strategy follows by classically coupling each coordinate at a time.

	\begin{Theorem}\label{Theorem: Coupling of two random walks on Zd}
		Let $\rho$ be a probability measure on $\Z^d$ such that
		\begin{enumerate}
			\item the relation $\rho(x) > 0$ implies $x = k u_i$ for some $k \in \Z$ and some canonical unitary vector $u_i;$
			\item $\rho(\Z^\ast u_i) > 0$ for each $i = 1, \ldots, d;$
			\item for every $i = 1, \ldots, d,$ the probability measures $\nu_i,$ defined in (\ref{Proposition: Decomposition of RW}), and $\bar{\nu_i}$ defined by $\bar{\nu}_i:k \mapsto \nu_i(-k)$ ($k \in \Z$) are such that their convolution $\nu_i \ast \bar \nu_i$ defines a transition probability kernel for a Markov chain on $\Z$ that is irreducible and recurrent.
	\end{enumerate}
		Then, for any pair of vertices $x$ and $y$ in $\Z^d,$ there exists a classical coupling of random walks started at $x$ and $y,$ respectively.
	\end{Theorem}

	\subsection{The Liouville property}\label{Section: Liouville property}
	We will say that a network satisfies
	\begin{description}[leftmargin=!,labelwidth=\widthof{($\mathsf{SLP}$)}]
		\item[($\mathsf{LP}$)] The \textbf{Liouville property}  if every \emph{bounded} harmonic function on the network is constant.
		\item[($\mathsf{SLP}$)] The \textbf{strong Liouville property}  if every \emph{positive} harmonic function on the network is constant.
		\item[($\mathsf{DLP}$)] The \textbf{Dirichlet Liouville property} if every \emph{Dirichlet} harmonic function is constant.
	\end{description}

	\begin{Remark}
		\begin{enumerate}
			\item For any network, we have the following implications from the above properties
			\begin{displaymath}\textstyle
				(\mathsf{SLP}) \implies (\mathsf{LP}) \implies (\mathsf{DLP}).
			\end{displaymath}
			The first of these implications is trivial for if $h$ is a bounded harmonic function, then $h + c$ is a positive harmonic function for large enough constant $c,$ then condition ($\mathsf{SLP}$) implies $h + c$ is constant, and so is $h.$ The second of this implications is harder and it fully proved in \cite{MaDib:USFwithDrift}, appendix D (the author followed closely the steps of Soardi \cite{Soardi:PTIN}).

			\item It is a well-known fact that the graph $\Z^d$ satisfies the strong Liouville property. One proof is given in \cite{Barlow:RWHK}. Here, M. Barlow followed an already established strategy: showing that $\Z^d$ satisfies the so-called \textbf{elliptic Harnack inequality} ($\mathsf{EHI}$) and also showing that for any network $(\mathsf{EHI}) \implies (\mathsf{SLP}).$

			\item The network $\Gamma_d(\lambda)$ does not satisfy ($\mathsf{SLP}$). Consider the function $h(z) = e^{-\lambda n}$ (with $z = (n, x)$ as usual). Then, with $p = (2d + 1 + e^\lambda)^{-1},$
			\begin{displaymath}\textstyle
				\sum\limits_{z' \sim z} p(z, z') h(z') = p e^{-\lambda (n - 1)} + 2d p e^{-\lambda n} + e^\lambda p e^{-\lambda (n + 1)} = e^{-\lambda n} = h(z),
			\end{displaymath}
			so that $h$ is harmonic, and proving the claim.

			\item If on a given network and for any pair of vertices, two lazy network random walks, one started at each vertex, can be classically coupled, then this network satisfies the Liouville property. In other words,
			\begin{displaymath}\textstyle
				(\mathsf{CC}) \implies (\mathsf{LP}).		
			\end{displaymath}
			First notice that a function $h:\strV \to \R$ is harmonic relative to the $\NRW$ if and only if it is harmonic relative to the $\frac{1}{2}$-lazy $\NRW.$ Next, consider a bounded harmonic function $h,$ any two vertices $x, y$ and the two $\frac{1}{2}$-lazy $\NRW$s $S$ and $S',$ started at $x$ and $y,$ respectively, with successful coupling time $\tau.$ Since $h(S)$ and $h(S')$ are bounded martingales, we obtain $h(x) = \Expectation{h(S_\tau)} = \Expectation{h(S_\tau')} = h(y).$
		\end{enumerate}
	\end{Remark}

	\begin{Theorem}
		The network $\Gamma_d(\lambda)$ has the Liouville property.
	\end{Theorem}
	\begin{proof}
		By virtue of the previous proposition, it suffices to prove that any two lazy network random walks (started at any two vertices) of $\Gamma_d(\lambda)$ can be classically coupled. In this case, the hypothesis of theorem (\ref{Theorem: Coupling of two random walks on Zd}) are satisfied: the measure $\nu_i \ast \bar \nu_i,$ where $\nu_i$ is as in (\ref{Remark: Decomposition of RW for lazy version in Gamma lambda}) is clearly symmetric and bounded, it is aperiodic by lazyness, thus $\nu_i \ast \bar \nu_i$ is also recurrent and irreducible.
	\end{proof}

	\subsection{Liouville property and the uniform spanning forest measure}
	The Liouville property is relevant to the context of uniform spanning forests due to the following theorem.
	\begin{Theorem}
		Let $\Gamma$ be any network. A sufficient condition for the measures $\WSF$ and $\FSF,$ of the wired and free uniform spanning forest on $\strG,$ to be the same, is that this network should satisfy the Liouville property.
	\end{Theorem}
	For the proof, see \cite[Theorem 7.3]{BLPS:USF}. As a corollary of the foregoing.
	\begin{Theorem}\label{Theorem: Equality between wired and free uniform spanning forest}
		$\WSF = \FSF$ on $\Gamma_d(\lambda).$
	\end{Theorem}

	By virtue of the equality between $\WSF$ and $\FSF$ on $\Gamma_d(\lambda),$ we will denote this measure as $\USF$ from now onwards.

	\subsection{The number of components}
	To count the number of components, the reader should recall the definition of transitive network. We know $\Gamma_d(\lambda)$ is transitive by considering translations.

	We need the following results:
	\begin{description}
		\item[\cite{LyPe:PTN}, Theorem 10.24.] In any transient transitive network $(\strG, \mu),$ with vertex set $\strV$ and Green's function $G,$ the number of intersections of any two independent network random walks is almost surely infinity if $\sum\limits_{z \in \strV} G(0, z)^2 = \infty$ and is almost surely finite otherwise.
		 \item[\cite{BLPS:USF}, Theorem 9.2.] On any network, for the number of components in $\WSF$ to be one, it is necessary and sufficient that, two independent network random walk started at any two vertices, should intersect (a.s.).
		 \item[\cite{BLPS:USF}, Theorem 9.4.] If there exists a vertex $o$ such that, for almost every realisation, the number of intersections of two independent copies of the network random walk started at $o$ is finite, then the number of components in $\WSF$ is infinite.
	\end{description}

	\begin{Theorem}\label{Theorem: Number of components in the forest}
		Consider the network $\Gamma_d(\lambda).$ The $\USF$ here (\ref{Theorem: Equality between wired and free uniform spanning forest}) consists of a single tree when $d = 1, 2$ and it has an infinite number of components when $d \geq 3.$ From now onwards, we use $\UST$ in lieu of $\USF$ for $d = 1, 2.$
	\end{Theorem}
	\begin{proof}
		By the theorems stated above, it all reduces to show that $\sum\limits_{z \in \Z^{d + 1}} G(0, z)^2$ diverges for $d = 1,2$ and that it converges otherwise. In the sake of clarity, we separate the calculations in the next proposition, with them the theorem is proved.
	\end{proof}

	As mentioned in the proof above, the next theorem states when the Green's function belongs to $\mathscr{L}^2.$

	\begin{Theorem}\label{Theorem: The bubble condition}{\normalfont (The ``bubble condition.'')}
		\begin{enumerate}
			\item If $d = 1, 2,$ $\sum\limits_{z \in \Z^{d + 1}} G(0, z)^2 = \infty.$
			\item If $d \geq 3,$ $\sum\limits_{z \in \Z^{d + 1}} G(0, z)^2 < \infty.$
		\end{enumerate}
	\end{Theorem}
	\begin{proof}
		Recall the Fourier transform of the step-lengths of the network random (\ref{Equation: Fourier transform of the step-lengths}). Now, if $(S_n)$ is the network random walk started at zero and $\varphi_n$ is its Fourier transform of $S_n,$ then $\varphi_n(h) = \varphi_Y(h)^n.$ Hence, for all $h \neq 0,$ $\widehat{G}(h) = \sum\limits_{n = 0}^\infty \varphi_n(h) = (1 - \varphi_Y(h))^{-1},$ where $\widehat{G}$ denotes the Fourier transform of the Green's function $G.$ By Parseval's relation,
		\begin{equation}\textstyle\label{Equation: Parsevals relation for the Green function}
			\sum\limits_{z \in \Z^{d + 1}} G(0, z)^2 = \frac{1}{(2\pi)^{d + 1}} \integral[{[-\pi, \pi]^{d + 1}}] dh\ \left| 1 - \varphi_Y(h) \right|^{-2},
		\end{equation}
		where the two integrals will converge to the same value, or they will both diverge.
	
		Using second order expansion in (\ref{Equation: Fourier transform of the step-lengths}) and with $a$ as in (\ref{Equation: First and second moment of the network}) and $p = (2d + 1 + e^\lambda)^{-1},$ we reach
		\begin{displaymath}\textstyle
			\varphi_Y(h) = 1 + i a h_0 - \frac{a}{2} h_0^2 - p \sum\limits_{j = 0}^d h_j^2 + \littleO{h_0^2} + \bigO{\norm{h}^4}, \quad h \to 0.
		\end{displaymath}
		Denote $\alpha_j = p$ for $j = 1, \ldots, d$ and $\alpha_0 = p + \frac{a}{2},$ $\varphi_Y(h) = 1 + i a h_0 - \sum\limits_{j = 0}^d \alpha_j h_j^2 + \littleO{h_0^2} + \bigO{\norm{h}^4}.$ Then, 
		\begin{align}\textstyle\tag{$\ast$}
			\left| 1 - \varphi_Y(h) \right|^2 &\textstyle= a^2 h_0^2 + i a h_0 \sum\limits_{j = 0}^d \alpha_j h_j^2 + \littleO{h_0^2} + \bigO{h_0 \norm{h}^4} \\\nonumber
			&\textstyle- i a h_0 \sum\limits_{j = 0}^d \alpha_j h_j^2 + \left( \sum\limits_{j = 0}^d \alpha_j h_j^2 \right)^2 + \bigO{\norm{h}^6} \\\nonumber
			&\textstyle= a^2 h_0^2 + \left( \sum\limits_{j = 0}^d \alpha_j h_j^2 \right)^2 + \littleO{h_0^2} + \bigO{h_0 \norm{h}^4} + \littleO{\norm{h}^4}.
		\end{align}
		Divide the right-hand side of (\ref{Equation: Parsevals relation for the Green function}) as the sum of $I_\varepsilon + J_\varepsilon,$ where $I_\varepsilon$ has domain of integration $\strK_\varepsilon = [-\varepsilon, \varepsilon] \times \ball[_{\R^d}]{0; \varepsilon}$ and $J_\varepsilon$ has domain of integration $[-\pi, \pi]^{d + 1} \setminus \strK_\varepsilon.$ It follows at once $J_\varepsilon < \infty.$ We will show the existence of a positive $\varepsilon$ such that $I_\varepsilon = \infty$ for $d = 1, 2$ and $I_\varepsilon < \infty$ for $d \geq 3.$
	
		Assume first $d = 1, 2.$ Then by ($\ast$) for $\varepsilon$ small enough, $|1 - \varphi_Y(h)|^2 \leq c^{-1}(h_0^2 + \norm{h'}^4)$ for some $c > 0.$ Then,
		\begin{displaymath}\textstyle
			I_\varepsilon \geq c \integral[{\ball[_{\R^d}]{0; \varepsilon}}] dh'\ \integral[-\varepsilon]^\varepsilon dh_0\ \frac{1}{h_0^2 + \norm{h'}^4} = 2c \integral[{\ball[_{\R^d}]{0; \varepsilon}}] dh'\ \integral[0]^\varepsilon dh_0\ \frac{1}{h_0^2 + \norm{h'}^4}.
		\end{displaymath}
		Changing measures to surface measure on the spheres of $\R^d,$ there is a constant $c > 0$ such that
		\begin{displaymath}\textstyle
			I_\varepsilon \geq c \integral[0]^\varepsilon d \rho\ \rho^{d - 1} \integral[0]^\varepsilon dh_0\ \frac{1}{h_0^2 + \rho^4} = c \integral[0]^\varepsilon d\rho\ \rho^{d - 3} \arctan \left( \frac{\varepsilon}{\rho^2} \right) = \infty.
		\end{displaymath}
	
		Similarly from ($\ast$), we can assume $\varepsilon > 0$ is small enough so that for some constant $c > 0,$ $\left| 1 - \varphi_Y(h) \right|^2 \geq c^{-1} (h_0^2 + \norm{h'}^4).$ As before, we change to surface measure of spheres and use that $\rho^{d - 1} \leq \rho^2$ since $\varepsilon$ can be assumed smaller than unity and $d \geq 3$, then
		\begin{displaymath}\textstyle
			I_\varepsilon \leq c \integral[\strK_\varepsilon] \frac{1}{h_0^2 + \norm{h'}^4} d(h_0, h') \leq c' \integral[0]^\varepsilon d \rho\ \rho^2 \cdot \frac{1}{\rho^2} \arctan \left( \frac{h_0}{\rho^2} \right) \Bigg|_{h_0 = 0}^{h_0 = \varepsilon} < \infty.
		\end{displaymath}
		Therefore, we reach that the sum (\ref{Equation: Parsevals relation for the Green function}) is finite for all $d \geq 3.$
	\end{proof}

	\section{Crossings of $\NRW$}\label{Section: Crossings of NRW}
	\subsection{Dimensions $d = 1, 3, 4, \ldots$}
	\begin{Theorem}\label{Theorem: In high dimensions two independent NRW may not cross at all}
		Let $d \geq 3.$ Consider two independent network random walks $S$ and $S'$ on $\Gamma_d(\lambda).$ Suppose they have arbitrary (random) initial points $S_0$ and $S_0',$ respectively. Then, almost surely, the paths of $S$ and $S'$ will cross finitely many times.
	\end{Theorem}
	\begin{proof}
		By the bubble condition (\ref{Theorem: The bubble condition}), $\sum\limits_{z \in \Z^{d + 1}} G(0, z)^2 < \infty$ in this range of dimensions. The theorem then follows from \cite[Theorem 10.24]{LyPe:PTN}.
	\end{proof}

	The previous theorem is not true in dimension $d = 1,2.$ However, there is an important distinction between the cases $d = 1$ and $d = 2.$ The gist of the distinction lies in how frequently the two random paths cross. In this section we consider the case $d = 1.$
	\begin{Theorem}\label{Theorem: In dimension d = 1 the difference of two random walks is recurrent}
		Let $d = 1.$ Consider two independent network random walks $S$ and $S'$ on $\Gamma_1(\lambda)$ with arbitrary starting points, $S_0$ and $S_0'.$ Then, $S - S'$ is a recurrent random walk; in particular, for almost every realisation, there will be infinitely many $n$ for which $S_n - S_n' = S_0 - S_0'.$ Thus, if both $S$ and $S'$ start at the same point, then $S$ and $S'$ will collide infinitely often.
	\end{Theorem}
	\begin{proof}
		By means of a translation by $-(S_0 - S_0'),$ we can assume $S_0 - S_0' = 0.$ The difference $S - S'$ is a random walk, started at zero and has step-length distribution
		\begin{align*}
			\mu &= p^2 e^\lambda \Big( \varepsilon_{(-2, 0)} + \varepsilon_{(2, 0)} \Big) + p^2 \Big( \varepsilon_{(0, 2)} + \varepsilon_{(0, -2)} \Big) \\
			&+ p^2(e^\lambda + 1) \Big( \varepsilon_{(-1, -1)} + \varepsilon_{(1, 1)} + \varepsilon_{(-1, 1)} + \varepsilon_{(1, -1)} \Big) + p^2 (e^{2 \lambda} + 3) \varepsilon_{(0, 0)},
		\end{align*}
		where $\varepsilon_{(h, k)}$ denotes the measures of unitary total mass at $(h, k) \in \Z^2$ and $p = \frac{1}{3 + e^\lambda}.$ Observe that $S - S'$ is a symmetric, aperiodic and irreducible random walk with bounded step lengths on the lattice generated by $u_1 = (1, 1)$ and $u_2 = (-1, 1).$ By \cite[Theorem 4.1.1]{LaLim:RWMI} $S - S'$ is recurrent.
	\end{proof}

	The remainder of this chapter is devoted into \emph{estimating how frequently} there will be a crossing between two network random walk paths on $\Gamma_2(\lambda).$ Whence, in the upcoming sections we will develop the necessary tools to tackle dimension $d = 2.$

	\subsection{Two elementary results}
	For the sake of reference, we write the following. The ``second moment inequality" states that for a non-negative random variable
	\begin{equation}\textstyle\label{Equation: Second moment lower bound}
		\Probability{Z > \frac{1}{2} \Expectation{Z}} \geq \frac{\Expectation{Z}^2}{4 \Expectation{Z^2}}.
	\end{equation}

	Consider now an increasing sequence of $\sigma$-fields $(\mathscr{E}_n)_{n \in \Z_+},$ with $\mathscr{E}_0$ the $\sigma$-field generated by $\varnothing,$ on some probability space $(\strE, \mathscr{E}, \mathds{P}).$ Suppose that $\strA_n$ is an event measurable up to time $n,$ that is $\strA_n \in \mathscr{E}_n.$ Set $q_n = \Probability{\strA_n \mid \mathscr{E}_{n - 1}}$ (so that $q_1 = \Probability{\strA_1}$). We have the following result.
	\begin{Proposition}\label{Proposition: Levys Borel Cantelli lemma}
		The event $\strN = \{\strA_n \mathrm{\ i.o.}\} \triangle \left\{\sum\limits_{n = 1}^\infty q_n = \infty \right\}$ is $\mathds{P}$-null (``Levy's extension of Borel-Cantelli lemma'').
	\end{Proposition}
	The proof of this result is a standard application of martingale theory to $X_n = \sum\limits_{j = 1}^n (\indic{\strA_j} - q_j).$ See section 7.8 of \cite{Ash:RAP} (apply theorem 7.8.4).

	\subsection{Setup for $\Gamma_2(\lambda)$}
	As before, we will write $z = (n, x);$ any affix, such as a sub or superscript, will be given to the three symbols: if we write $z'$ at some point, it will be assumed that $z' = (n', x'),$ similarly, if we talk about $n_1$ and $x_1,$ we will be assuming $z_1 = (n_1, x_1),$ etc.

	Divide the space $\Z \times \Z^2$ into the two sets $\strS_1 = \{z; n \leq \norm{x}\}$ and $\strS_2 = \{z; n > \norm{x}\}.$ Notice $0 \in \strS_1$ and $\strS_1 \cap \strS_2 = \varnothing.$ To simplify notation, we define (for any $y \in \Z^r, r \geq 1$)
	\begin{equation}\textstyle\label{Equation: Japanese bracket}
		\japBracket{y} = \max(\norm{y}, 1),
	\end{equation}
	and notice $\japBracket{y} = \norm{y}$ unless $y = 0;$ also $\japBracket{y} \asymp \norm{y} + 1 \asymp \sqrt{\norm{y}^2 + 1}.$

	Consider the following tessellation ($p, q \in \N$ and $i = 1,2$) of the ``positive half-space'' $\Z_+ \times \Z^2$:
	\begin{equation}\textstyle\label{Equation: Division of the space into a tessellation}
		\begin{split}
			\strR_1(p, q) &= \{z : 9^p  < n \leq 2 \cdot 9^p, 3^{q - 1} < \norm{x} \leq 3^q \}, \\
			\strR_2(p, q) &= \{z : 2 \cdot 9^p  < n \leq 9^{p + 1}, 3^{q - 1} < \norm{x} \leq 3^q \}\\
			\strR_1(p, 0) &= \{z : 9^p  < n \leq 2 \cdot 9^p, \norm{x} \leq 1 \},\\
			\strR_2(p, 0) &= \{z : 2 \cdot 9^p  < n \leq 2 \cdot 9^{p + 1}, \norm{x} \leq 1 \},\\
			\strR_1(0, q) &= \{z : 0 \leq n \leq 2, 3^{q - 1} < \norm{x} \leq 3^q \}\\
			\strR_2(0, q) &= \{z : 3 \leq n \leq 9, 3^{q - 1} < \norm{x} \leq 3^q \}\\
			\strR_1(0, 0) &= \{z : 0 \leq n \leq 2, \norm{x} \leq 1\} \\
			\strR_2(0, 0) &= \{z : 3 \leq n \leq 9, \norm{x} \leq 1\}
		\end{split}
	\end{equation}
	For $p \in \N,$ we shall consider the following region
	\begin{equation}\textstyle\label{Equation: Definition of the testing regions}
		\strD_p = \bigcup\limits_{q = 0}^p \strR_1(p, q) = \{z : 9^p < n \leq 2 \cdot 9^p, \norm{x} \leq 3^p\}.
	\end{equation}
	The following easy remark contains all the technical estimates that will be employed later in this section. Consider the Green's function bounds and fix the constants $c_i$ ($1 \leq i \leq 4$) of (\ref{Theorem: Greens function bounds}).
	\begin{Remark}\label{Remark: Principal calculations for the number of crossings in dimension d equals 2}
		For $i = 1, 2$ and $(p, q) \in \Z_+^2,$ we have the following.
		\begin{enumerate}
			\item $\card{\strR_i(p, q)} \asymp 9^{p + q}.$
	
			\item $\strR_{i_1}(p_1, q_1)$ and $\strR_{i_2}(p_2, q_2)$ are disjoint unless $(i_1, p_1, q_1) = (i_2, p_2, q_2).$
	
			\item $\card{\strD_p} \asymp 9^p \sum\limits_{q = 0}^p 9^q \asymp 9^{2p}.$

			\item We have $\strR_i(p, q) \subset \strS_1$ if $2p + 3 \leq q$ and $\strR_i(p, q) \subset \strS_2$ if $q \leq 2p.$
	
			\item If $z \in \Z_+ \times \Z^2$ is in $\strS_1$ and $c$ is a positive constant, then $0 \leq n \leq \norm{x},$ this implies $\norm{x} \leq \norm{z} \leq \sqrt{2} \norm{x},$ and then $\frac{1}{\sqrt{2}} e^{-\sqrt{2}c\norm{x}} \japBracket{x}^{-1} \leq e^{-c \norm{z}} \japBracket{z}^{-1} \leq e^{-c \norm{x}} \japBracket{x}^{-1}.$ Therefore, if $3^{q - 1} \leq \norm{x} \leq 3^q,$ $\frac{1}{\sqrt{2}} e^{-\sqrt{2} c 3^q} 3^{-q} \leq e^{-c \norm{z}} \japBracket{z}^{-1} \leq 3 e^{-\frac{c}{3} 3^q} 3^{-q}.$
	
			\item If $z \in \Z_+ \times \Z^2$ is in $\strS_2$ and $c$ is a positive constant, then $\norm{x} < n,$ this implies $0 < n \leq \norm{z} \leq \sqrt{2} n,$ and then $\frac{1}{\sqrt{2}} e^{-c \frac{\norm{x}^2}{n}} n^{-1} \leq e^{-c \frac{\norm{x}^2}{n}} \norm{z}^{-1} \leq e^{-c \frac{\norm{x}^2}{n}} n^{-1}.$ Therefore, if $9^p \leq n \leq 9^{p + 1}$ and $3^{q - 1} \leq \norm{x} \leq 3^q$ $\frac{1}{9\sqrt{2}} e^{-c 9^{q - p}} 9^{-p} \leq e^{-c \frac{\norm{x}^2}{n}} \norm{z}^{-1} \leq e^{-\frac{c}{81} 9^{q - p}} 9^{-p}.$
	
			\item Define the function $\varphi_\strU:\Z \times \Z^2 \to \R$ given by $\varphi_\strU(z) = c_1 e^{-c_2 \norm{z}} \japBracket{z}^{-1}$ if $z$ lies in region $\strS_1$ and $\varphi_\strU(z) = c_1 e^{-c_2 \frac{\norm{x}^2}{n}} \norm{z}^{-1}$ for $z$ within $\strS_2.$ Then $G(0, z) \leq \varphi_\strU(z).$ If $z \in \strR_1(p, q) \cup \strR_2(p, q),$
			\begin{displaymath}
				\varphi_\strU(z) \leq \begin{cases}
					3c_1 e^{-\frac{c_2}{3} 3^q} 3^{-q} &\text{ if } 2p + 3 \leq q, \vspace{0.2cm} \\
					c_1 e^{-\frac{c_2}{81} 9^{q - p}} 9^{-p} &\text{ if } q \leq 2p + 2.
				\end{cases}
			\end{displaymath}
	
			\item Define $\varphi_\strL(z) = c_3 e^{-c_4 \norm{z}} \japBracket{z}^{-1}$ if $z \in \strS_1$ and $\varphi_\strL(z) = c_3 e^{-c_4 \frac{\norm{x}^2}{n}} \norm{z}^{-1}$ if $z \in \strS_2.$ Then, $G(0, z) \geq \varphi_\strL(z).$ If $z \in \strR_i(p, q)$ ($i = 1,2$),
			\begin{displaymath}
				\varphi_\strL(z) \geq \begin{cases}
					\frac{c_3}{\sqrt{2}} e^{-\sqrt{2} c_4 3^q} 3^{-q} &\text{ if } 2p + 3 \leq q, \vspace{0.25cm} \\
					\frac{c_3}{9 \sqrt{2}} e^{-c_4 9^{q - p}} 9^{-p} &\text{ if } q \leq 2p + 2.
				\end{cases}
			\end{displaymath}

			\item For $z = (n, x) \in \Z \times \Z^2$ define $\dot{z} = (-n, x) \in \Z \times \Z^2$ (the ``reflection'' of $z$ through the ``plane'' $\{0\} \times \Z^2$). Then, $z \mapsto \dot{z}$ is a bijection of $\Z_\pm \times \Z^2$ onto $\Z_\mp \times \Z^2,$ it is the identity when restricted to $\{0\} \times \Z^2.$ For every $z \in \Z_- \times \Z^2,$ $\varphi_\strU(z) \leq \varphi_\strU \left( \dot{z} \right).$

			\item There are constants $d_1 = 3c_1,$ $d_2 = \frac{c_2}{81},$ $d_3 = \frac{c_3}{18}$ and $d_4 = 2c_4,$ such that the relation $z \in \strR_i(p, q)$ ($p, q \in \Z_+,$ $i = 1,2$) imply
			\begin{displaymath}
				G(0, z)  \leq \varphi_\strU(z) \leq \begin{cases}
					d_1 e^{-d_2 3^q} 3^{-q} &\text{ if } 2p + 1 \leq q, \\
					d_1 e^{-d_2 9^{q - p}} 9^{-p} &\text{ if } q \leq 2p.
				\end{cases}
			\end{displaymath}
			and
			\begin{displaymath}
				G(0, z) \geq \varphi_\strL(z) \geq \begin{cases}
					d_3 e^{-d_4 3^q} 3^{-q} &\text{ if } 2p + 1 \leq q, \\
					d_3 e^{-d_4 9^{q - p}} 9^{-p} &\text{ if } q \leq 2p.
				\end{cases}
			\end{displaymath}

			\item $G(0, z) \asymp 9^{-p}$ with any constant being uniform for $z \in \strD_p$ and $p \in \N.$
		\end{enumerate}
	\end{Remark}

		\subsubsection{Some preliminary sums of the Green's function of $\Gamma_2(\lambda).$}
		In what follows, the estimates in Remark (\ref{Remark: Principal calculations for the number of crossings in dimension d equals 2}) will be used, often without referencing to this remark.
		\begin{Proposition}\label{Proposition: Sum of the squares of the Green function}
			Let $S$ and $S'$ be two independent network random walks started at $0$ in $\Gamma_2(\lambda).$ Then, with any implicit constants being independent of $p \in \N,$
			\begin{enumerate}
				\item $\sum\limits_{z \in \strD_p} G(0, z)^2 \asymp 1;$ and
				\item $\sum\limits_{z \in \strD_p - \strD_p} G(0, z) G(0, -z) \asymp 1.$
				\item Also, there exists a constant $c > 0$ such that $\sum\limits_{z \in \strD_p - \strD_p} G(0, z)^2 \leq c p,$ for all $p \in \N.$
			\end{enumerate}
		\end{Proposition}
		\begin{proof}
			We know $G(0, z)^2 \asymp 9^{-2p}$ and $\card{\strD_p} \asymp 9^{2p}$ for $z \in \strD_p$ and the implicit constants independent of $p;$ this proves the first assertion $\sum\limits_{z \in \strD_p} G(0, z)^2 \asymp 1.$ There remains to show the other two estimates. First, it is clear that $\sum\limits_{z \in \strD_p - \strD_p} G(0, z) G(0, -z) \geq G(0,0)^2.$ Hence, there remains to prove the existence of constants $c > 0$ and $c' > 0$ such that $\sum\limits_{z \in \strD_p - \strD_p} G(0, z) G(0, -z) \leq c$ and $\sum\limits_{z \in \strD_p - \strD_p} G(0, z)^2 \leq c' p,$ for all $p \in \N.$ Set $\strQ^{(p)} = \strD_p - \strD_p$ and $\strQ_+^{(p)} = \strQ^{(p)} \cap (\Z_+ \times \Z^2).$ Consider the function $\varphi_\strU$ introduced in (\ref{Remark: Principal calculations for the number of crossings in dimension d equals 2}), we regard the constants $c_k$ and $d_k$ ($k = 1, \ldots, 4$) of (\ref{Remark: Principal calculations for the number of crossings in dimension d equals 2}) to be fixed during the argument. Then $\sum\limits_{z \in \strQ^{(p)}} G(0, z) G(0, -z) \leq 2 \sum\limits_{z \in \strQ_+^{(p)}} G(0, z) G(0, -z) \leq 2 \sum\limits_{z \in \strQ_+^{(p)}} \varphi_\strU(z) \varphi_\strU(-z),$ and $\sum\limits_{z \in \strQ^{(p)}} G(0, z)^2 \leq \sum\limits_{z \in \strQ^{(p)}} \varphi_\strU(z)^2 \leq 2 \sum\limits_{z \in \strQ_+^{(p)}} \varphi_\strU(z)^2$ where the second inequality follows from considering the reflection $z \mapsto \dot{z}$ (\ref{Remark: Principal calculations for the number of crossings in dimension d equals 2})(i). To simplify notation, introduce $\strR(s, t)$ to be the union $\strR_1(s, t) \cup \strR_2(s, t),$ these  sets being defined in (\ref{Equation: Division of the space into a tessellation}). We know that there exists a constant $L > 0$ such that for all $s, t \in \Z_+,$ $\card{\strR(s, t)} \leq L 9^{s + t}.$ Then, one can write $\strQ_+^{(p)} \subset \bigcup\limits_{(s, t)} \strR(s, t),$ where the indices run over all choices of $0 \leq s, t \leq p + 1.$ Now, the sets $\strR(s, t)$ are pairwise disjoint and, therefore, we split $\sum\limits_{z \in \strQ_+^{(p)}} \varphi_\strU(z) \varphi_\strU(-z) \leq \sum\limits_{(s, t)} \sum\limits_{z \in \strR(s, t)} \varphi_\strU(z) \varphi_\strU(-z) = P_1 + P_2,$ 	and, in a similar fashion, $\sum\limits_{z \in \strQ_+^{(p)}} \varphi_\strU(z)^2 \leq T_1 + T_2,$ where $P_1$ and $T_1$ are the sums corresponding to $t < s$ and $P_2$ and $T_2,$ the sums for $t \geq s.$ We shall show that $P_1$ and $P_2$ are bounded by constants independent of $p,$ and $T_1 \leq c p$ and $T_2 \leq c' p$ for a pair of constants $c, c' > 0$ independent of $p;$ having this, we will conclude $\sum\limits_{z \in \strD_p - \strD_p} G(0, z) G(0, -z) \leq c$ and $\sum\limits_{z \in \strQ^{(p)}} G(0, z)^2 \leq c' p$ for some constants $c > 0$ and $c' > 0,$ as desired.
		
			We now prove $P_1$ and $P_2$ are bounded by constants independent of $p.$ Observe from (\ref{Remark: Principal calculations for the number of crossings in dimension d equals 2})(g) that for all $z \in \strR(s, t),$ either $	\varphi_\strU(z) \leq d_1 e^{-d_2 3^t} 3^{-t}$ or $\varphi_\strU(z) \leq d_1 e^{-d_2 9^{t - s}} 9^{-s}.$ Also, for $z \in \strQ_+^{(p)}$ we have $\varphi_\strU(-z) = c_1 e^{-c_2 \norm{z}} \norm{z}^{-1}.$ Therefore (use (\ref{Remark: Principal calculations for the number of crossings in dimension d equals 2})(a)),
			\begin{align*}\textstyle
				P_2 &\textstyle\leq \sum\limits_{s = 0}^{p + 1} \sum\limits_{t = s}^{p + 1} \sum\limits_{z \in \strR(s, t)} \varphi_\strU(z) \varphi_\strU(-z) \leq 3c_1 d_1 L \sum\limits_{s = 0}^{p + 1} \sum\limits_{t = s}^{p + 1} \big( e^{-d_2 3^t} 9^s + e^{-d_2 9^{t - s}} 3^t \big) e^{-\frac{c_2}{3} 3^t} \\
				&\textstyle\leq 3c_1 d_1 L  \sum\limits_{l = 0}^\infty e^{-\frac{c_2}{3} 3^l} \left( \left[ \sum\limits_{k = 0}^\infty e^{-d_2 3^k} \right] 9^l + \left[ \sum\limits_{k = 0}^\infty e^{-d_2 9^k} 3^k \right] 3^l \right) < \infty.
			\end{align*}
			The sum $P_1$ is easier,
			\begin{align*}\textstyle
				P_1 &\textstyle\leq \sum\limits_{s = 0}^{p + 1} \sum\limits_{t = 0}^{s - 1} \sum\limits_{z \in \strR(s, t)} \varphi_\strU(z) \varphi_\strU(-z) \leq c_1 d_1 L \sum\limits_{s = 0}^{p + 1} \sum\limits_{t = 0}^{s - 1} 9^{s + t} e^{-d_2 9^{t - s}} 9^{-s} e^{-c_2 9^s} 9^{-s} \\
				&\textstyle\leq c_1 d_1 L \sum\limits_{s = 0}^{p + 1} e^{-c_2 9^s} 9^{-s} \sum\limits_{t = 0}^{s - 1} 9^t = \frac{c_1 d_1 L}{8} \sum\limits_{k = 0}^\infty e^{-c_2 9^k} < \infty.
			\end{align*}
		
			We continue with the proof of the bounds of $T_1$ and $T_2,$ beginning with $T_2 \leq c p,$ for some constant $c > 0.$  By definition of $T_2,$ $T_2 = \sum\limits_{s = 0}^{p + 1} \sum\limits_{t = s}^{p + 1} \sum\limits_{z \in \strR(s, t)} \varphi_\strU(z)^2.$ As before, we have $\varphi_\strU(z)^2 \leq d_1^2 \Big(  e^{-2d_2 3^t} 9^{-t} + e^{-2d_2 9^{t - s}} 9^{-2s} \Big).$ Therefore, with $C = d_1^2 L,$
			\begin{align*}\textstyle
				T_2 &\textstyle\leq C \sum\limits_{s = 0}^{p + 1} \sum\limits_{t = s}^{p + 1} \left( e^{-2 d_2 3^t} 9^{-t} + e^{-2 d_2 9^{t - s}} 9^{-2s} \right) 9^{s + t} = C \sum\limits_{s = 0}^{p + 1} \sum\limits_{t = s}^{p + 1} \left( e^{-2 d_2 3^t} 9^s + e^{-2 d_2 9^{t - s}} 9^{t - s} \right) \\
				&\textstyle\leq C \sum\limits_{k = 0}^\infty e^{-d_2 3^k} 9^k \sum\limits_{k = 0}^\infty e^{-d_2 3^k} + C (p + 2) \sum\limits_{k = 0}^\infty e^{-2 d_2 9^k} 9^k \leq C' p,
			\end{align*}
			where $C' = C \sum\limits_{k = 0}^\infty e^{-d_2 3^k} 9^k \sum\limits_{k = 0}^\infty e^{-d_2 3^k} + 3C \sum\limits_{k = 0}^\infty e^{-2 d_2 9^k} 9^k < \infty$ 	is independent of $p.$ Similarly, notice one can express $T_1 = \sum\limits_{s = 0}^{p + 1} \sum\limits_{t = 0}^{s - 1} \sum\limits_{z \in \strR(s, t)} \varphi_\strU(z)^2,$ bearing this in mind and that $\varphi_\strU(z) \leq d_1 9^{-s}$ for $0 \leq t \leq s$ it follows immediately that, for some positive constants $c$ and $c',$ $	T_1 \leq c \sum\limits_{s = 0}^{p + 1} \sum\limits_{t = 0}^s 9^{s + t} 9^{-2 s} \leq c' p,$ as desired.
		\end{proof}
	
		\subsubsection{Estimates on the Green's function on some regions of $\Gamma_2(\lambda).$}
		Denote with $S$ and $S'$ two independent network random walks of $\Gamma_2(\lambda).$

		Consider the following region ($p \in \N$)
		\begin{equation}\textstyle\label{Equation: Region of restriction for the Green function}
			\strU_p = \left\{ (n, x) \in \Z^3 : |n| \leq 4 \cdot 9^{p k_0}, \norm{x} \leq 3^{(p + 1)k_0} \right\},
		\end{equation}
		where $k_0$ is a positive integer to be determined later on. Define the 	``separating cylinders'' (creating a barrier between $\strD_{p k_0}$ and $\strD_{(p + 1) k_0}$ for $p \in \N$) to be the (vertex) boundary $\strF_p = \frontier \strU_p,$ and set $\strF_0 = \{0\}.$ Having these sets, define ($p \in \Z_+$)
		\begin{equation}\textstyle\label{Equation: Hitting time of the frontier Fp}
			T_p = \tau_{\strF_p}(S)
		\end{equation}
		the hitting time of $\strF_p$ by $S.$ If $z \in \Z^3,$ let
		\begin{equation}\textstyle\label{Equation: The first index such that S0 belong to Up}
			p_z = \min\{p \mid z \in \strU_p\}.
		\end{equation}
		Having defined $p_z,$ observe that when $S_0 = o,$ then a.s. the sequence of random variables $(T_p)_{p \geq p_o}$ is increasing and each one of them is finite (by transience).

		Next, by definition $\strF_p \cap \strU_p = \varnothing.$ It will be useful to divide the set $\strF_p$ into several regions. Geometrically, the set $\strF_p$ can be thought as a cylinder and we need to consider the body and the two bases separately. Analytically, if $(n, x) \in \strF_p$ then $|n| > 4 \cdot 9^{p k_0}$ or else $\norm{x} > 3^{(p + 1)k_0}.$ We can specialise further for either $|n| = 4 \cdot 9^{p k_0} + 1$ and $\norm{x} \leq 3^{(p + 1)k_0}$ or else, $|n| \leq 4 \cdot 9^{p k_0}$ and $3^{(p + 1)k_0} \leq \norm{x} \leq 3^{(p + 1)k_0} + 1.$ Thus, the exponential scales being considered permit us to assume that $\norm{z},$ for $z \in \strF_p,$ satisfies :
		\begin{enumerate}
			\item $|n| = 4 \cdot 9^{p k_0}$ and $\norm{x} \leq 3^{(p + 1)k_0},$ or
			\item $|n| < 4 \cdot 9^{p k_0}$ and $\norm{x} = 3^{(p + 1)k_0}.$
		\end{enumerate}
		With this assumption, we introduce the following sets. The ``left-hand base of the cylinder $\strF_p$'' is the set of $z \in \strF_p$ such that $n = -4 \cdot 9^{p k_0};$ the ``right-hand base of the cylinder $\strF_p$'' is the set of $z \in \strF_p$ such that $n = 4 \cdot 9^{p k_0}.$ The points in $\strF_p$ that belong to neither the left nor the right-hand base of the cylinder will be referred to as the ``body of the cylinder $\strF_p,$'' in other words, the body of the cylinder is the set of points $z \in \strF_p$ such that $|n| < 4 \cdot 9^{p k_0}$ and $\norm{x} = 3^{(p + 1)k_0}.$
	
		Define $(T_p')$ as above changing the process $S$ for the process $S'.$ Observe that $(T_p)$ and $(T_p')$ are independent families.
	
		Having the times $(T_p)$ and $(T_p')$ we now consider ``sections'' of the network random walks and the number of intersections in $\strD_{p k_0}$ of two such sections. Thus, define ($p \in \N$)
		\begin{equation}\textstyle\label{Equation: Definition of the number of crossing in the sections}
			M_p = \sum\limits_{z \in \strD_{p k_0}} \sum\limits_{m = T_{p - 1}}^{T_p} \sum\limits_{m' = T_{p - 1}'}^{T_p'} \indic{\{S_m = z\}} \indic{\{S_{m'}' = z\}}.
		\end{equation}
		Observe that $M_p$ counts the number of intersections of the two random walks inside the region $\strD_{p k_0}$ starting in the separating cylinder $\strF_{p - 1}.$ We remark that, if $S_0 = o,$ then $S$ cannot be inside $\strD_{p k_0}$ before the time $T_{p - 1}$ provided $p - 1 \geq p_o;$ similarly, $S'$ cannot be inside $\strD_{p k_o}$ before time $T_{p - 1}'$ provided $S_0' = o'$ and $p - 1 \geq p_{o'}.$

		We will study the Green's function restricted to the region $\strU_p:$ $G_{\strU_p}(0, z) =  \Expectation[^0]{\sum\limits_{m = 0}^{T_p - 1} \indic{\{S_m = z\}}},$ for $z \in \Z^3.$ For consistency we denote $\strU_0 = \{0\}.$ For the sake of simplicity, we will write $G_p = G_{\strU_p}.$
		\begin{Proposition}\label{Proposition: Asymptotic bound between the Green function and the restricted Green function on d equals 2}
			For every $p \in \N,$ $\strD_{p k_0} - \strU_{p - 1} \subset \strU_p.$ Moreover, for any $\varepsilon > 0$ there exists an index $k_0 = k_0(\varepsilon)$ such that for some $p_0 \in \N$ the relations $p \geq p_0$ and $z \in \strD_{p k_0} - \strU_{p - 1}$ imply $G(0, z) \leq (1 + \varepsilon) G_p(0, z).$
		\end{Proposition}
		\begin{proof}
			The proof contains several steps. We begin by noticing that for $z_1 \in \strD_{p k_0}$ and $z_2 \in \strU_{p - 1}$ one has by (\ref{Equation: Definition of the testing regions}) and (\ref{Equation: Region of restriction for the Green function}) $9^{p k_0} < n_1 \leq 2 \cdot 9^{p k_0}$ and $\norm{x_1} \leq 3^{p k_0},$ and $-4 \cdot 9^{(p - 1) k_0} \leq n_2 \leq 4 \cdot 9^{(p - 1) k_0}$ and $\norm{x_2} \leq 3^{p k_0},$ which implies
			\begin{EqInside}\label{Equation: Size of the coordinates in Dpk0 - Up-1}
				Any $(n, x) \in \strD_{p k_0} - \strU_{p - 1}$ satisfies $\frac{1}{2} \cdot 9^{p k_0} \leq n \leq \frac{5}{2} 9^{p k_0}, \norm{x} \leq 2 \cdot 3^{p k_0}.$
			\end{EqInside}
			\noindent In particular, $z_1 - z_2$ belongs to the region $\strU_p$ so this establishes the first assertion in (\ref{Proposition: Asymptotic bound between the Green function and the restricted Green function on d equals 2}).
		
				We now prove the second assertion of (\ref{Proposition: Asymptotic bound between the Green function and the restricted Green function on d equals 2}). First, $G(0, z) = G_p(0, z) + \Expectation[^0]{G(S_{T_p}, z)}.$ Observe that $S_{T_p}$ is in $\strF_p$ and, therefore, we need to estimate $G(z_3, z)$ for $z_3 \in \strF_p.$ This will be done in the upcoming three lemmas.
	
			We will pause momentarily the proof of the proposition in order to establish these lemmas and resume it once they have been proven.
		
			\begin{Lemma}\label{Lemma: Asymptotic bound for Green function restricted to Up}
				\begin{enumerate}
					\item For every $\varepsilon > 0,$ there exists a $p_0 = p_0(\varepsilon)$ such that for all $p \geq p_0,$ all $k_0 \in \N,$ and all $z_1 \in \strD_{p k_0}, z_2 \in \strU_{p - 1}$ and $z_3 \in \strF_p$ for which $n_1 - n_2 - n_3 \leq 2 \cdot 3^{(p + 1) k_0}$ and $z_3$ does not lie in the left-hand base of $\strF_p,$ we have $G(z_3, z_1 - z_2) \leq \varepsilon G(0, z_1 - z_2).$
	
					\item There exists a constant $c > 0$ such that for $p \in \N$ and all $z_1 \in \strD_{p k_0}, z_2 \in \strU_{p - 1}$ and $z_3 \in \strF_p$ for which $n_1 - n_2 - n_3 \geq 2 \cdot 3^{(p + 1) k_0}$ and $z_3$ does not lie in the left-hand base of $\strF_p$ then $G(z_3, z_1 - z_2) \leq c 9^{-k_0} G(0, z_1 - z_2).$
				\end{enumerate}
			\end{Lemma}
			\begin{proof}[Proof of lemma.]
				Suppose the constants $c_k$ ($1 \leq k \leq 4$) of (\ref{Theorem: Greens function bounds}) are given. By (\ref{Equation: Size of the coordinates in Dpk0 - Up-1}) and the estimates on the Green's function, we have $G(0, z_1 - z_2) \asymp 9^{-p k_0}.$
	
				We prove part (a) first. Assume $z_1, z_2, z_3$ are as stated.  Consider first the case $n_3 < 4 \cdot 9^{p k_0},$ that is to say, consider first when $z_3$ is not in the right-hand base of $\strF_p.$ Then $\norm{x_3} = 3^{(p + 1) k_0},$ which in turn yields $\frac{1}{3} \cdot 3^{(p + 1) k_0} \leq \norm{x_1 - x_2 - x_3} \leq 2 \cdot 3^{(p + 1) k_0}.$ Suppose first that $n_1 - n_2 - n_3 \leq \frac{1}{3} \cdot 3^{(p + 1) k_0},$ then $z_1 - z_2 - z_3$ is in $\strS_1,$ implying by (\ref{Theorem: Greens function bounds})
				\begin{align*}
					G(z_3, z_1 - z_2) &\leq c_1 e^{-c_2 \norm{z_1 - z_2 - z_3}} \norm{z_1 - z_2 - z_3}^{-1} \leq 3c_1 e^{-\frac{c_2}{3} 3^{(p + 1) k_0}} 3^{-(p + 1) k_0} \\
					&= \left[ 3c_1 e^{ -\frac{c_2}{3} 3^{(p + 1) k_0} } 3^{(p - 1) k_0} \right] 9^{-p k_0} \leq \left[ 3c_1 e^{ -\frac{c_2}{3} 3^{p k_0} } 3^{p k_0} \right] 9^{-p k_0}.
				\end{align*}
				Define $u_1(t) = 3c_1 e^{-\frac{c_2}{3} 3^t} 3^t,$ so that $G(z_3, z_1 - z_2) \leq u_1(p k_0) 9^{-p k_0},$ for $z_1,$ $z_2$ and $z_3$ as stated. Next, suppose that $\frac{1}{3} \cdot 3^{(p + 1) k_0} \leq n_1 - n_2 - n_3 \leq 2 \cdot 3^{(p + 1) k_0}.$ Then,
				\begin{align*}
					G(z_3, z_1 - z_2) &\leq c_1 \left(e^{-c_2 \norm{z_1 - z_2 - z_3}} + e^{-c_2 \frac{\norm{x_1 - x_2 - x_3}^2}{n_1 - n_2 - n_3}} \right) \norm{z_1 - z_2 - z_3}^{-1} \\
					&\leq 3c_1  \left(e^{-\frac{c_2}{3} 3^{(p + 1) k_0}} + e^{-\frac{c_2}{18} 3^{(p + 1) k_0}} \right) 3^{-(p + 1) k_0} \\
					&= \left[ 3c_1  \left(e^{-\frac{c_2}{3} 3^{(p + 1) k_0}} + e^{-\frac{c_2}{18} 3^{(p + 1) k_0}} \right) 3^{(p - 1) k_0} \right] 9^{-p k_0} \\
					&\leq \left[ 3c_1  \left(e^{-\frac{c_2}{3} 3^{p k_0}} + e^{-\frac{c_2}{18} 3^{p k_0}} \right) 3^{p k_0} \right] 9^{-p k_0} = u_2(p k_0) 9^{-p k_0},
				\end{align*}
				where $u_2(t)$ is defined to be $3c_1  \left(e^{-\frac{c_2}{3} 3^t} + e^{-\frac{c_2}{18} 3^t} \right) 3^t.$ Consider now the case where $z_3$ is in the right-hand base of the cylinder, that is to say, assume $n_3 = 4 \cdot 9^{p k_0}$ and $\norm{x_3} \leq 3^{(p + 1) k_0}$ but in this case, $-\frac{7}{2} \cdot 9^{p k_0} \leq n_1 - n_2 - n_3 \leq -\frac{3}{2} \cdot 9^{p k_0}$ and again, by (\ref{Theorem: Greens function bounds}), $G(z_3, z_1 - z_2) \leq c_1 e^{-c_2 |n_1 - n_2 - n_3|} |n_1 - n_2 - n_3|^{-1} \leq \frac{2 c_1}{3} e^{-\frac{2 c_2}{3} 9^{p k_0}} 9^{-p k_0}.$ Here we set $u_3(t) = \frac{2 c_1}{3} e^{-\frac{2 c_2}{3} 9^t},$ so that in this case $G(z_3, z_1 - z_2) \leq u_3(p k_0)  9^{-p k_0}.$ Notice that in all three cases above, we bounded $G(z_3, z_1 - z_2) \leq u_i(p k_0) 9^{-p k_0},$ and it is clear that $u_i(t) \to 0$ as $t \to \infty$ in each case. Then, there exists a $t_0$ such that $\sup\limits_{t \geq t_0} \max\limits_{i = 1, 2, 3} u_i(t) \leq \varepsilon.$ We can define $p_0$ to be the least integer greater than or equal to $t_0;$ this $p_0$ works for all $k_0$ since $p k_0 \geq p.$ This finishes part (a).
			
				We now prove part (b). Suppose that $z_1, z_2, z_3$ are as stated in (b). This implies the relation $n_1 - n_2 - n_3 \geq 2 \cdot 3^{(p + 1) k_0}.$ By (\ref{Equation: Size of the coordinates in Dpk0 - Up-1}), we have $n_1 - n_2 - n_3 \leq 7 \cdot 9^{p k_0}$ and, therefore $n_1 - n_2 -n_3$ lies in an interval of the form $\strJ_\alpha = \left[ (\alpha - 1) 3^{(p + 1) k_0}, \alpha 3^{(p + 1) k_0} \right]$ where $\alpha = 3, \ldots, 7 \cdot 3^{(p - 1) k_0}.$ Regardless the value of $\alpha,$ the relation $n_1 - n_2 - n_3 \in \strJ_\alpha$ implies $z_1 - z_2 - z_3 \in \strS_2$ and these relations imply,
				\begin{displaymath}
					G(z_3, z_1 - z_2) \leq c_1 e^{-c_2 \frac{\norm{x_1 - x_2 - x_3}^2}{n_1 - n_2 - n_3}} (n_1 - n_2 - n_3)^{-1} \leq 2 c_1 e^{-\frac{c_2}{9} \alpha^{-1} 3^{(p + 1) k_0}} \alpha^{-1} 3^{-(p + 1) k_0}.
				\end{displaymath}
				Observe that for any constant $c > 0,$ the function $u_c(t) = t^{-1} e^{-c t^{-1}}$ for $t > 0$ has an absolute and global maximum at $t = c.$ Its maximum equals $u_c(c) = e^{-1} c^{-1}.$ We deduce $G(z_3, z_1 - z_2) \leq 2c_1 e^{-1} \frac{9}{c_2} 9^{-(p + 1) k_0} = c' 9^{-k_0} 9^{-p k_0},$ the desired conclusion since $G(0, z_1 - z_2) \asymp 9^{-p k_0},$ with any implicit constant independent of $p$ and of $k_0.$
			\end{proof}
	
			\begin{Lemma}
				For any $\varepsilon > 0,$ there exist indices $k_0 = k_0(\varepsilon) \in \N$ and $p_0  = p_0(\varepsilon) \in \N$ such that the relations $p \geq p_0$ and $z \in \strD_{p k_0} - \strU_{p - 1}$ imply $\sup\limits_{\substack{z_3 \in \strF_p \\ n_3 > -4 \cdot 9^{p k_0}}} G(z_3, z) \leq \varepsilon G(0, z).$
			\end{Lemma}
			\begin{proof}[Proof of lemma.]
				By (\ref{Lemma: Asymptotic bound for Green function restricted to Up}), given $\varepsilon > 0$ there is a $p_0 \in \N$ such that if $p \geq p_0,$ $k_0 \in \N$ and for all $z_1 \in \strD_{p k_0},$ $z_2 \in \strU_{p - 1}$ and $z_3 \in \strF_p$ with $n_3 > -4 \cdot 9^p,$ we either have $G(z_3, z_1 - z_2) \leq \varepsilon G(0, z_1 - z_2)$ or $G(z_3, z_1 - z_2) \leq c 9^{-k_0} G(0, z_1 - z_2),$ where the constant $c$ is universal. Choose $k_0 \in \N$ such that $c 9^{-k_0} \leq \varepsilon.$ Hence, $G(z_3, z_1 - z_2) \leq \varepsilon G(0, z_1 - z_2)$ uniformly in $z_3 \in \strF_p$ with $n_3 > -4 \cdot 9^p$ for all $p \geq p_0,$ which is the conclusion to be reached.
			\end{proof}
		
			\begin{Lemma}\label{Lemma: Probability STp lies in the left hand base}
				For any $\varepsilon > 0,$ there exist indices $k_0 = k_0(\varepsilon) \in \N$ and $p_0 = p_0(\varepsilon) \in \N$ such that the relation $p \geq p_0$ implies $\Probability[^0]{S_{T_p} \in \{z \in \strF_p : n = -4 \cdot 9^p\}} \leq \varepsilon.$ Also, there is a universal constant $c > 0$ such that for all $p \in \N$ and all $z \in \strD_{p k_0} - \strU_{p - 1},$ $\sup\limits_{\substack{z_3 \in \strF_p \\ n_3 = -4 \cdot 9^p}} G(z_3, z) \leq cG(0, z).$
			\end{Lemma}
			\begin{proof}[Proof of lemma.]
				Set $\strA_p = \{z \in \strF_p : n = -4 \cdot 9^{p k_0}\}.$ Consider the network random walk $S$ in $\Gamma_2(\lambda)$ started at zero and define $B = \pr{0}(S).$ By (\ref{Proposition: Probability for the existence of splitting levels}), $\Probability{S_{T_p} \in \strA_p} \leq \Probability{B \text{ ever reaches } -4 \cdot 9^{p k_0}} = e^{- 4 \lambda 9^{p k_0}}.$ This proved the first claim. To prove the second claim, use (\ref{Equation: Size of the coordinates in Dpk0 - Up-1}), if $z_3 \in \strA_p$ and $z \in \strD_{p k_0} - \strU_{p - 1}$ then $G(z_3, z) \asymp 9^{-p k_0}$ with any implicit constant being universal in $p, k_0, z, z_3.$ Similarly, $G(0, z) \asymp 9^{-p k_0}$ for $z \in \strD_{p k_0} - \strU_{p - 1}$ with any implicit constant universal. The second claim follows readily.
			\end{proof}

			We are ready to continue with the proof of Proposition (\ref{Proposition: Asymptotic bound between the Green function and the restricted Green function on d equals 2}). Let us write $\strA_p = \{z \in \strF_p : n_3 = -4 \cdot 9^p\}$ for the left-hand base of $\strF_p$ and $\strB_p = \strF_p \setminus \strA_p$ for the remainder of the separating cylinder. For the time being, consider $\varepsilon > 0$ to be any positive number. By Lemmas (\ref{Lemma: Asymptotic bound for Green function restricted to Up}) and (\ref{Lemma: Probability STp lies in the left hand base}), there are indices $k_0$ and $p_0,$ only depending on $\varepsilon,$ and a universal constant $c > 0$ (not depending on any index) satisfying for all $p \geq p_0,$ $z_1 \in \strD_{p k_0},$ $z_2 \in \strU_{p - 1}$ $\sup\limits_{z_3 \in \strB_p} G(z_3, z_1 - z_2) \leq \varepsilon G(0, z_1 - z_2),$ $\Probability[^0]{S_{T_p} \in \strA_p} \leq \varepsilon$ 	and for all $z \in \strD_{p k_0} - \strU_{p - 1},$ $\sup\limits_{z_3 \in \strA_p} G(z_3, z) \leq cG(0, z).$ Write now, for $z \in \strD_{p k_0} - \strU_{p - 1},$
			\begin{align*}
				G(0, z) &= G_p(0, z) + \Expectation[^0]{G(S_{T_p}, z) \indic{\strA_p}(S_{T_p})} + \Expectation[^0]{G(S_{T_p}, z) \indic{\strB_p}(S_{T_p})} \\
				&\leq G_p(0, z) + c G(0, z) \Probability[^0]{S_{T_p} \in \strA_p} + \varepsilon G(0, z) \Probability[^0]{S_{T_p} \in \strB_p} \\
				&\leq G_p(0, z) + (c + 1) \varepsilon G(0, z).
			\end{align*}
		
			Hence, for any $\varepsilon > 0,$ there are indices $k_0$ and $p_0,$ depending only on $\varepsilon,$ so that the relations $p \geq p_0$ and $z \in \strD_{p k_0} - \strU_{p - 1}$ imply $G(0, z) \leq G_p(0, z) + \frac{\varepsilon}{1 + \varepsilon} G(0, z),$ which yields $G(0, z) \leq (1 + \varepsilon) G_p(0, z).$
		\end{proof}

	\subsection{Dimension $d = 2;$ ``logarithmic scale of crossings'' of $\Gamma_2(\lambda)$}
	We will establish a \emph{lower bound} on the probability of how frequently two random walk paths in $\Gamma_2(\lambda)$ will cross each other.

	Invoking Proposition (\ref{Proposition: Asymptotic bound between the Green function and the restricted Green function on d equals 2}), for $\varepsilon = 1$ there exists two indices $k_0$ and $p_0,$ depending only on $\varepsilon,$ such that the conclusions of this proposition hold. For the remainder of this chapter, we assume both $k_0$ and $p_0$ are fixed and given as stated.
	\begin{Theorem}\label{Theorem: Logarithmic scale of crossings}
		Consider two independent network random walks $S$ and $S'$ in $\Gamma_2(\lambda),$ with respective starting points $o$ and $o'.$ Define the random variables $T_p$ and $T_p'$ in (\ref{Equation: Hitting time of the frontier Fp}) using $S$ and $S',$ accordingly; now introduce $M_p$ as in (\ref{Equation: Definition of the number of crossing in the sections}). Consider now the $\sigma$-fields $\mathscr{E}_p  = \scrF^S_{T_p} \vee \scrF^{S'}_{T_p'}.$ Define the event $\strA_p = \{M_p > 0\}$ and the random variable $q_{p + 1} = \Probability{\strA_{p + 1} \mid \mathscr{E}_p}.$ Then, $\strA_p \in \mathscr{E}_p$ and there exists a positive integer $P = P(o, o')$ and a universal constant $c > 0$ (depending neither on $o$ nor on $o'$) such that, almost surely, $q_p > \frac{c}{p}$ for all $p \geq P.$
	\end{Theorem}
	\begin{proof}
		It is clear  that $\strA_p \in \mathscr{E}_p.$ Consider $p_1 = \max(p_o, p_{o'})$ (\ref{Equation: The first index such that S0 belong to Up}), that is, $p_1$ is the first index such that $o$ and $o'$ belong to $\strU_p.$ Suppose $p$ is any integer $\geq \max(p_0, p_1).$ By (\ref{Proposition: Markov property simultaneous for two Markov processes}), we have $q_{p + 1} = \Probability{\strA_{p + 1} \middle| S_{T_p}, S_{T_p'}'}.$ It is now straightforward to check that
		\begin{equation}\textstyle
			q_{p + 1} = \sum\limits_{(z_1, z_2) \in \strF_p^2} \Probability{\strA_{p + 1} \middle| S_{T_p} = z_1, S_{T_p'}' = z_2} \indic{\left\{S_{T_p} = z_1, S_{T_p'}' = z_2\right\}}. \tag{$\ast$}
		\end{equation}
		We will use (\ref{Equation: Second moment lower bound}). By the strong Markov property for a pair of independent Markov processes (\ref{Proposition: Markov property simultaneous for two Markov processes}) gives (inside the expectation on the right of the first equality we use abuse of notation and assume $(S_m)$ and $(S_{m'}')$ are independent network random walks started at zero)
		\begin{align*}\textstyle
			\Expectation{M_{p + 1} \middle| S_{T_p} = z_1, S_{T_p'}' = z_2} &\textstyle= \sum\limits_{z_3 \in \strD_{(p + 1) k_0}} \Expectation{ \sum\limits_{m = 0}^{T_{p + 1}} \sum\limits_{m' = 0}^{T_{p + 1}'} \indic{\{S_m + z_1 = z_3\}} \indic{ \left\{ S_{m'}' + z_2 = z_3 \right\} } } \\
			&\textstyle= \sum\limits_{z_3 \in \strD_{(p + 1) k_0}} G_{p + 1}(0, z_3 - z_1) G_{p + 1}(0, z_3 - z_2) \\
			&\textstyle\asymp \sum\limits_{z_3 \in \strD_{(p + 1) k_0}} G(0, z_3 - z_1) G(0, z_3 - z_2), \quad \text{by } (\ref{Proposition: Asymptotic bound between the Green function and the restricted Green function on d equals 2}) \\
			&\textstyle\asymp \sum\limits_{z_3 \in \strD_{(p + 1) k_0}} 9^{-2 (p + 1)k_0} \asymp 1, \quad \text{ by } (\ref{Remark: Principal calculations for the number of crossings in dimension d equals 2})(k).
		\end{align*}
		We shall now estimate the second moment of $M_p$ with respect to the conditional measure. We have
		\begin{displaymath}\textstyle
			M_{p + 1}^2 = \sum\limits_{(z_3, z_4) \in \strD_{(p + 1)k_0}^2} \sum\limits_{m_1, m_2 = T_p}^{T_{p + 1}} \sum\limits_{m_1', m_2' = T_p'}^{T_{p+1}'} \indic{\{S_{m_1} = z_3, S_{m_2} = z_4\}} \indic{ \left\{ S_{m_1'}' = z_3, S_{m_2'}' = z_4 \right\} }.
		\end{displaymath}
		Now apply the conditional expectation, by the strong Markov property (\ref{Proposition: Markov property simultaneous for two Markov processes}) and independence between $S$ and $S'$ (we again employ abuse of notation and assume $S$ and $S'$ are two independent network random walks started at zero after having used the strong Markov property for two independent Markov processes)
		\begin{align*}\textstyle
			\Expectation{M_{p + 1}^2 \middle| S_{T_p} = z_1, S_{T_p'}' = z_2} &\textstyle= \sum\limits_{(z_3, z_4) \in \strD_{(p + 1)k_0}^2} \Expectation{ \sum\limits_{m_1, m_2 = 0}^{T_{p + 1}} \indic{ \{S_{m_1} = z_3 - z_1, S_{m_2} = z_4 - z_1 \} } } \\
			&\textstyle\times \Expectation{ \sum\limits_{m_1', m_2' = 0}^{T_{p + 1}'} \indic{ \left\{S_{m_1'}' = z_3 - z_2, S_{m_2'}' = z_4 - z_2 \right\} } } \\
			&\textstyle\leq B_1 + B_2 + B_3 + B_4,
		\end{align*}
		in which
		\begin{displaymath}
			\begin{split}
				B_1 &\textstyle= \sum\limits_{ \substack{m_1 \leq m_2 \\ m_1' \leq m_2' \\ (z_3, z_4) \in \strD_{(p + 1)k_0}^2} } \Probability{S_{m_1} = z_3 - z_1, S_{m_2} = z_4 - z_1} \Probability{ S_{m_1'}' = z_3 - z_2, S_{m_2'}' = z_4 - z_2 }, \\
				B_2 &\textstyle= \sum\limits_{ \substack{m_1 \leq m_2 \\ m_1' \geq m_2' \\ (z_3, z_4) \in \strD_{(p + 1)k_0}^2} } \Probability{S_{m_1} = z_3 - z_1, S_{m_2} = z_4 - z_1 } \Probability{S_{m_1'}' = z_3 - z_2, S_{m_2'}' = z_4 - z_2 }, \\
				B_3 &\textstyle= \sum\limits_{ \substack{m_1 \geq m_2 \\ m_1' \leq m_2' \\ (z_3, z_4) \in \strD_{(p + 1)k_0}^2} } \Probability{ S_{m_1} = z_3 - z_1, S_{m_2} = z_4 - z_1 } \Probability{ S_{m_1'}' = z_3 - z_2, S_{m_2'}' = z_4 - z_2 }, 
			\end{split}
		\end{displaymath}
		and
		\begin{displaymath}\textstyle
			B_4 = \sum\limits_{ \substack{m_1 \geq m_2 \\ m_1' \geq m_2' \\ (z_3, z_4) \in \strD_{(p + 1)k_0}^2} } \Probability{ S_{m_1} = z_3 - z_1, S_{m_2} = z_4 - z_1 } \Probability{ S_{m_1'}' = z_3 - z_2, S_{m_2'}' = z_4 - z_2 }.
		\end{displaymath}

		We handle each sum separately.
		\begin{enumerate}
			\item By independent increments and writing $p^k(z, z') = \Probability{S_k = z' \mid S_0 = z},$
			\begin{align*}
					B_1 &\textstyle= \sum\limits_{ \substack{m_1 \leq m_2 \\ m_1' \leq m_2' \\ (z_3, z_4) \in \strD_{(p + 1)k_0}^2} } p^{m_1}(0, z_3 - z_1) p^{m_2 - m_1}(0, z_4 - z_3) p^{m_1'}(0, z_3 - z_2) p^{m_2'-m_1'}(0, z_4 - z_3) \\
				&\textstyle= \sum\limits_{(z_3, z_4) \in \strD_{(p + 1)k_0}^2} G(0, z_3 - z_1) G(0, z_4 - z_3)^2 G(0, z_3 - z_2),
			\end{align*}
			the last equality being an immediate consequence of Lebesgue-Tonelli's theorem and the definition of the Green's function of $\Gamma_2(\lambda).$ By virtue of (\ref{Equation: Size of the coordinates in Dpk0 - Up-1}) and (\ref{Proposition: Sum of the squares of the Green function}) we can proceed further,
			\begin{align*}
				B_1 &\textstyle\asymp \sum\limits_{(z_3, z_4) \in \strD_{(p + 1)k_0}^2} 9^{-2(p + 1)k_0} G(0, z_4 - z_3)^2 \leq c k_0 (p + 1) \sum\limits_{z_3 \in \strD_{(p + 1)k_0}} 9^{-2(p + 1)k_0} \asymp p.
			\end{align*}
	
			\item Here we get
			\begin{align*}
				B_2 &\textstyle= \sum\limits_{(z_3, z_4) \in \strD_{(p + 1)k_0}^2} G(0, z_3 - z_1) G(0, z_4 - z_3) G(0, z_4 - z_2) G(0, z_3 - z_4) \\
				&\textstyle\asymp \sum\limits_{z_3 \in \strD_{(p + 1)k_0}} \sum\limits_{z_4 \in \strD_{(p + 1)k_0}} 9^{-2(p + 1)k_0}  G(0, z_4 - z_3) G(0, z_3 - z_4) \\
				&\textstyle\asymp \sum\limits_{z_3 \in \strD_{(p + 1)k_0}} 9^{-2(p + 1)k_0}  \asymp 1,
			\end{align*}
			where the second $\asymp$ is a consequence of (\ref{Proposition: Sum of the squares of the Green function}).
	
			\item Similarly to $B_2,$ we get $B_3 \asymp 1.$
	
			\item As with $B_1,$ $B_4 \leq c p,$ for some universal constant $c > 0.$
		\end{enumerate}
	
		Therefore, we have shown that there exist two positive constants $c$ and $c',$ such that, whenever $p \geq \max(p_0, p_1)$ and $(z_1, z_2) \in \strF_p^2,$ $\Expectation{M_{p + 1} \middle| S_{T_p} = z_1, S_{T_p'}' = z_2} \geq c$ and $\Expectation{M_{p + 1}^2 \middle| S_{T_p} = z_1, S_{T_p'}' = z_2} \leq c' p.$ By (\ref{Equation: Second moment lower bound}), we finally reach that there exists a constant $c > 0$ and a $P = \max(p_0, p_1)$ such that, for all $p \geq P$ and all $(z_1, z_2) \in \strF_p^2,$ $\Probability{\strA_{p + 1} \middle| S_{T_p} = z_1, S_{T_p'}' = z_2} \geq \frac{c}{p},$ comparing with the expression for $q_{p + 1}$ in ($\ast$), we have substantiated the claims.
	\end{proof}

	\begin{Corollary}\label{Corollary: Lower bound in the probability of intersection of two sections of random walk in Gamma2}
		Let $S$ and $S'$ be two independent network random walks in $\Gamma_2(\lambda),$ with starting points $S_0 = o$ and $S_0' = o'.$ There is a universal constant $c > 0$ and an index $P = P(o, o'),$ such that for all $p \geq P,$ the probability that the two random walk paths $(S_m)_{m = T_{p - 1}, \ldots, T_p}$ and $(S_{m'}')_{m' = T_{p - 1}', \ldots, T_p'}$ cross each other inside $\strD_{p k_0}$ is at least $\frac{c}{p}.$
	\end{Corollary}

	The following corollary follows using (\ref{Proposition: Levys Borel Cantelli lemma}) and (\ref{Theorem: Logarithmic scale of crossings}).
	\begin{Corollary}\label{Corollary: In dimension d = 2 there are infinitely many crossings}
		Any two independent network random walks on $\Gamma_2(\lambda)$ will have infinitely many intersections.
	\end{Corollary}

	\begin{Remark}
		There is a much simpler way to obtain Corollary (\ref{Corollary: In dimension d = 2 there are infinitely many crossings}). Denote $S = (B, X),$ that is, $S_n = (B_n, X_n)$ for $n \in \Z_+;$ here $B$ is the biased random walk on $\Z$ and $X$ a random walk on $\Z^2.$ Define $\zeta_k = \tau_{\{k\} \times \Z^2}(S)$ for $k \in \Z_+,$ which is the time at which the random walk $B$ visits $k.$ With these hitting times, consider $Y_n = X_{\zeta_n}.$ Similarly, define $Y_n'$ using $S' = (B', X'),$ which is an independent copy of $S.$ In a similar manner as theorem (\ref{Theorem: In dimension d = 1 the difference of two random walks is recurrent}), it is seen that $Y_n = Y_n'$ for infinitely many indices $n.$ Having this, it is now clear that $S_n = S_{n'}'$ for infinitely many pairs $(n, n').$ Of course, this result is much coarser than (\ref{Corollary: Lower bound in the probability of intersection of two sections of random walk in Gamma2}).
	\end{Remark}

	\section{One endedness for $\Gamma_d(\lambda)$}\label{Section: One endedness}
Let $\strG$ be a directed infinite tree. A \textbf{ray} is an infinite path that does not repeat vertices; any two rays of $\strG$ will be disjoint or will eventually merge, this creates equivalence classes on the set of rays of $\strG.$ Any of these classes is called an \textbf{end}.

	\subsection{One endedness for $\Gamma_1(\lambda)$}
	We refer the reader to \cite{LyPe:PTN}, see Section 6.5 , Section 9.2 and Section 10.3 for the definitions and notations of planar graph and duality of graphs.

	\begin{Theorem}\label{Theorem: One end in dimension one}
		$\UST$ in $\Gamma_1(\lambda)$ (\ref{Theorem: Number of components in the forest}) has, almost surely, one end.
	\end{Theorem}
	\begin{proof}
		It is easy to see that the dual network $\Gamma_1(\lambda)^\dagger$ of $\Gamma_1(\lambda)$ is the network image of the latter under the reflection $(n, x) \mapsto (-n, x).$ In particular, the wired and free uniform spanning forest of $\Gamma_1(\lambda)^\dagger$ coincide and the number of $\USF$ trees in $\Gamma_1(\lambda)^\dagger$ is the same as that of $\Gamma_1(\lambda);$ that is to say, for almost every realisation, just one tree, by (\ref{Theorem: Number of components in the forest}). Bearing this in mind and the duality between $\UST$ in $\Gamma_1(\lambda)$ with $\UST$ in $\Gamma_1(\lambda)^\dagger,$ we see at once that the tree in $\Gamma_1(\lambda)$ is almost surely one ended; for if it were two ended, it would split the plane into two disconnected regions making it impossible for the dual tree to be a tree.
	\end{proof}

	\begin{Remark}
		The core idea of the proof of previous theorem is implicit in \cite[Theorem 10.36]{LyPe:PTN}.
	\end{Remark}

	\subsection{One endedness for $\Gamma_2(\lambda)$}
	We will prove the following theorem: $\UST$ in $\Gamma_2(\lambda)$ has, almost surely, one end.

	The proof is rather involved and to facilitate the reading of it we are going to develop it in several steps and draw the conclusion at the end.

		\subsection{Some preliminary results on Markov processes}
		Consider a discrete time Markov process $(S_n, \scrF_n)_{n \in \Z_+},$ with values on a metrisable separable space $\strE,$ endowed with Borels sets $\scrE$ as $\sigma$-algebra. Let $\eta$ be an a.s. finite stopping time for this Markov process. We may construct the $\sigma$-field $\scrF_\eta = \scrF_\eta^S$ of all events $\strG \in \scrE$ such that $\strG \cap \{\eta \leq n\} \in \scrF_n$ (``stopped filtration up to time $\eta$''). In a similar way, we may construct $S_\eta:\omega \mapsto S_{\eta(\omega)}(\omega).$ If no filtration is specified, it is assumed the canonical filtration is used $\scrF_n^S = \sigma(S_k; 0 \leq k \leq n).$ Consider now the product space $\strE^{\Z_+},$ which is also a metrisable separable space, and its Borel $\sigma$-field coincides with $\ds \bigotimes_{n \in \Z_+} \scrE.$ For any bounded measurable function $\varphi:\strE^{\Z_+} \to \R,$ there exists a version of the conditional expectation $\Expectation{ \varphi \left( S_j \right)_{j \in \Z_+} \middle| S_0 = x};$ denote this function as $\Expectation[^x]{\varphi \left( S_j \right)_{j \in \Z_+}}.$ The strong Markov property states that $\Expectation{\varphi \left( (S_j)_{j \geq \eta} \right) \middle| \scrF_\eta} = \Expectation[^{S_\eta}]{\varphi \left( (S_j)_{j \in \Z_+} \right)}.$

		The following results follow from well-known techniques (monotone-class, Dynkin's theorem, etc.). The reader may consult \cite{MaDib:USFwithDrift}, propositions (1.12.3) and (1.12.4).

		\begin{Proposition}\label{Proposition: Markov property simultaneous for two Markov processes}
			Suppose $S = (S_n)_{n \in \Z_+}$ and $S' = (S_n')_{n \in \Z_+}$ are two independent Markov processes, defined on some probability space $\left( \Omega, \scrF, \mathds{P} \right)$ with values on some metrisable, separable space $\strE,$ endowed with Borel $\sigma$-field $\scrE.$ Let $\psi:\strE^{\Z_+} \times \strE^{\Z_+} \to \R$ be a bounded measurable function, relative to the Borel sets of both $\R$ and $\strE^{\Z_+} \times \strE^{\Z_+}.$ Define $v_\psi(x, x') = \Expectation{\psi \left( S_j, S_{j'}' \right)_{(j, j') \in \Z_+^2} \middle| S_0 = x, S_0' = x'},$ for $(x, x') \in \strE^2,$ that is $v_\psi$ is a version of $\psi \left( S_j, S_{j'}' \right)_{(j, j') \in \Z_+^2}$ given $(S_0, S_0').$ Suppose $\eta$ and $\eta'$ are two stopping times, $\eta$ relative to $S$ and $\eta',$ to $S'.$ With $\psi$ as before, $\Expectation{\psi(S_j, S_{j'}')_{j \geq \eta, j' \geq \eta'} \middle| \scrF_\eta^S \vee \scrF_{\eta'}^{S'}} = v_\psi(S_\eta, S_\eta').$
		\end{Proposition}
	
		\begin{Proposition}\label{Proposition: Independence between past and future of two processes with stationary independent increments a stopping times}
			Let $\strE = \R^d$ (or any separable, metrisable, Fr\'{e}chet space). Suppose $(S_n)_{n \in \Z_+}$ and $(S'_n)_{n \in \Z_+}$ are two $\strE$-valued stochastic processes with stationary independent increments. Suppose that $\eta$ and $\eta'$ are stopping times relative to $S$ and $S',$ respectively. Then, $\scrF_\eta^S \vee \scrF_{\eta'}^{S'}$ is independent of $\scrG_\eta^S \vee \scrG_{\eta'}^{S'},$ with $\scrG_\eta^S = \sigma(S_{\eta + j} - S_\eta; j \in \Z_+),$ and a corresponding definition for $\scrG_{\eta'}^{S'}$ with $(S', \eta')$ replacing $(S, \eta).$ In particular, if $\EuScript{A}$ is an event depending on the paths $(S_j)_{0 \leq j \leq \eta}$ and $(S_{j'}')_{0 \leq j' \leq \eta'}$ and $\EuScript{B},$ depending on the paths $(S_{\eta + j} - S_\eta)_{j \in \Z_+}$ and $(S_{\eta' + j'}' - S_{\eta'}')_{j' \in \Z_+},$ then $\EuScript{A}$ and $\EuScript{B}$ are independent. 
		\end{Proposition}

		\subsubsection{Invariance of the chronology of Wilson's algorithm rooted at infinity.}
		We start with a theorem related to the study of the chronological construction of Wilson's algorithm rooted at infinity on general networks. Denote by $\scrS_\strG$ the (metrisable compact) topological space of spanning subgraphs of the graph $\strG$ and recall the notations from the introduction. With these notations, we construct the following random object $\left( \mathfrak{F}^\xi, \left( L_k^\xi \right)_{k \in \N} \right) = \left( \mathfrak{F}^\xi, L_1^\xi, \ldots, L_k^\xi, \ldots \right).$ This random object is $\scrS_\strG \times \scrS_\strG^\N = \scrS_\strG^{\Z_+}$-valued. Notice that the sequence $\left( L_k^\xi \right)_{k \in \N}$ may be referred to as a ``chronology of the $\WSF$-forest.'' Wilson's algorithm rooted at infinity gives $\mathfrak{F}^\xi \sim \mathfrak{F}^\eta.$ Observe that if $\xi$ and $\eta$ are two orderings satisfying $\xi(j) = \eta(j)$ for $1 \leq j \leq M$ then $\left( L_1^\xi, \ldots, L_M^\xi \right) = \left( L_1^\eta, \ldots, L_M^\eta \right)$ a.s.. Furthermore, in the finite case we may run Wilson's algorithm with stacks (\ref{Proposition: WA with stacks}) to obtain that $\mathfrak{F}^\xi = \mathfrak{F}^\eta$ a.s., thus in the finite case $\left( \mathfrak{F}^\xi, L_1^\xi, \ldots, L_M^\xi \right) = \left( \mathfrak{F}^\eta, L_1^\eta, \ldots, L_M^\eta \right).$

		\begin{Theorem}
			If $\xi$ and $\eta$ are two orderings of the vertex set of $\strG$ such that $\xi(j) = \eta(j)$ for $1 \leq j \leq M,$ then $	\left( \mathfrak{F}^\xi, L_1^\xi, \ldots, L_M^\xi \right) \sim \left( \mathfrak{F}^\eta, L_1^\eta, \ldots, L_M^\eta \right).$ (The two random objects have the same law.)
		\end{Theorem}
		The proof of this theorem follows from finite approximation, see \cite[Theorem (5.2.1)]{MaDib:USFwithDrift}.

		\begin{Corollary}\label{Corollary: Main extension of WA}
			The statistical properties of events depending on $\WSF$ and its first branches are not affected by reordering vertices not yet searched. In other words, if $\varphi:\scrS_\strG^{M + 1} \to \R$ is any bounded measurable function or any non negative measurable function, and $\xi, \eta:\N \to \strV$ are two orderings of the vertices such that $\xi(j) = \eta(j)$ for $1 \leq j \leq M,$ then $\Expectation{\varphi \left( \mathfrak{F}^\xi, L_1^\xi, \ldots, L_M^\xi \right) } = \Expectation{\varphi \left( \mathfrak{F}^\eta, L_1^\eta, \ldots, L_M^\eta \right)}.$
		\end{Corollary}

		The importance of the previous corollary is that one can see a partial forest, and to estimate \emph{probabilities} of successive branches, we can choose the next vertices to depend on the partial tree. (This is no surprise, since Wilson's algorithm with stacks allows for a stronger result in the finite case: the next vertices may depend on the partial tree and the final tree cannot change.) 

		\subsubsection{Preliminary estimates in the proof of one-endedness for $\Gamma_2(\lambda).$}
		We denote, for $r > 0,$
		\begin{displaymath}\textstyle
			\ball{x; r} = \{y \in \Z^d \mid \norm{x - y} \leq r\}, \quad	\sphere{x; r} = \{y \in \Z^d \mid r \leq \norm{x - y} \leq r + 1\}.
		\end{displaymath}
		We remark that there exists a constant $L = L(d) > 0$ such that for all $x \in \Z^d$ and all $r \geq 1,$
		\begin{displaymath}\textstyle
			\max \big( \card{\ball{x; r}}, \card{\sphere{x; r}} \big) \leq L \cdot r^d.
		\end{displaymath}
		Also, $\ball{x; r} = x + \ball{0; r}$ and $\sphere{x; r} = x + \sphere{0; r}.$ Further, we show that $\frontier \ball{x; r} \subset \sphere{x; r}.$ To see this, suffices to show this when $x = 0.$ In this case, the relation $y \in \frontier \ball{0; r}$ signifies there exists a $y' \in \ball{0; r}$ such that $\norm{y - y'} = 1,$ and then $\norm{y} \leq r + 1.$ The inequality $\norm{y} \geq r$ for $y \in \frontier \ball{0; r}$ follows from definition of the vertex boundary of a set.
	
		\begin{Proposition}\label{Proposition: Sum of the Green function on cylinder can be made small}
			For each $p \in \N,$ define the sets $\strA_p = [-p, p] \times \sphere{0; p},$ $\strB_p = \{p\} \times \ball{0; p}$ and $\strC_p = \{-p\} \times \ball{0; p}.$ Then, there is a pair of constants $c, c' > 0$ such that for all $p \in \N,$ $\sum\limits_{z \in \strA_p \cup \strB_p} G(z, 0) \leq c e^{-c' p}.$ Furthermore, there exists a another pair of constants $c, c' > 0$ and a positive number $B,$ obeying the following: for every $b \geq B$ we can split $\strC_p = \strC_{p, 0} \cup \strC_{p, 1},$ where $\strC_{p, 0} = \strC_{p, 0}(b) = \left( \{-p\} \times \ball{0; b p^{\frac{1}{2}}} \right) \cap \strC_p$ and $\strC_{p, 1} = \strC_{p, 1}(b) = \strC_p \setminus \strC_{p, 0},$ 	and with this division, $\sup\limits_{p \in \N} \sum\limits_{z \in \strC_{p, 1}} G(z, 0) \leq c e^{-c'b^2}.$ In particular, for any $\varepsilon > 0,$ there exists an index $p_1 = p_1(\varepsilon) \in \N$ and a positive number $b_1 = b_1(\varepsilon)$ such that, for any $p \geq p_1$ and $b \geq b_1,$ we can divide $\strC_p = \strC_{p, 0} \cup \strC_{p, 1}$ as before, and with this division $	\sup\limits_{\substack{p \geq p_1 \\ b \geq b_1}} \sum\limits_{z \in \strE_p} G(z, 0) \leq \varepsilon,$ where $\strE_p = \strA_p \cup \strB_p \cup \strC_{p, 1}.$
		\end{Proposition}
		\begin{proof}
			If $z \in \strA_p \cup \strB_p$ then the bounds on the Green's function (\ref{Theorem: Greens function bounds}) give at once $G(z, 0) = G(0, -z) \leq c_1 e^{-c_2 \|z\|} \|z\|^{-\frac{d}{2}} \leq c_1 e^{-c_2 p} p^{-\frac{d}{2}}.$ Now, $\card{\strA_p \cup \strB_p} \asymp p^d$ and so, for a constant $c_3 > 0,$ $\sum\limits_{z \in \strA_p \cup \strB_p} G(z, 0) \leq c_3 e^{-c_2 p} p^{\frac{d}{2}} \leq \left[ c_3 \sup\limits_{t \geq 1} e^{-\frac{c_2}{2} t} t^{\frac{d}{2}} \right] e^{-\frac{c_2}{2} p} = c_4 e^{-\frac{c_2}{2} p}.$ Now, for any $b > 0,$ we may construct $\strC_{p, 0}$ and $\strC_{p, 1}$ as in the statement; if $b \geq p^{\frac{1}{2}},$ then $\strC_{p, 1} = \varnothing,$ which is fine. Write $\strC_{p, 1} \subset \bigcup\limits_{b p^{\frac{1}{2}} \leq k \leq p} \{-p\} \times \sphere{0; k},$ so that $\sum\limits_{z \in \strC_{p, 1}} G(z, 0) \leq \sum\limits_{b p^{\frac{1}{2}} < k \leq p} \card{\sphere{0; k}} c_1 e^{-c_2 \frac{k^2}{p}} p^{-\frac{d}{2}} \leq \sum\limits_{b p^{\frac{1}{2}} < k \leq p} c_5 e^{-c_2 \frac{k^2}{p}} p^{-\frac{d}{2}} k^d,$ for a constant $c_5 > 0.$ The function $t \mapsto e^{-c \frac{t^2}{p}} t^d$ has a maximum at $t = \left( \frac{d}{2c} \right)^{\frac{1}{2}} p^{\frac{1}{2}}$ and it is decreasing for $t \geq \left( \frac{d}{2c} \right)^{\frac{1}{2}} p^{\frac{1}{2}};$ any value $b \geq \left( \frac{d}{2c} \right)^{\frac{1}{2}}$ allows bounding $\sum\limits_{z \in \strC_{p, 1}} G(z, 0) \leq c_5 \integral[b p^{\frac{1}{2}}]^p dt\ e^{-c_2 \frac{t^2}{p}} p^{-\frac{d}{2}} t^{d - 1} = c_5 \integral[b]^{\sqrt{p}} dt\ e^{-c_2 t^2} t^{d - 1} \leq c_6 e^{-\frac{c_2}{2} b^2}.$ Hence, putting all together, we have established that, for any $p \in \N$ and $b \geq \left( \frac{d}{2c} \right)^{\frac{1}{2}},$ $\sum\limits_{z \in \strA_p \cup \strB_p} G(z, 0) \leq c e^{-c' p}$ and $\sum\limits_{z \in \strC_{p, 1}} G(z, 0) \leq C e^{-C' b^2},$ where all these constants $c, c', C, C' > 0$ are independent of both $p$ and $b.$ The desired results are now clear.
		\end{proof}
	
		\begin{Proposition}\label{Proposition: Existence of sparse good sets}
			Let $\varepsilon > 0$ and we maintain the notations and conclusions of (\ref{Proposition: Sum of the Green function on cylinder can be made small}). Then, there exists a positive integer $p_2 = p_2(\varepsilon)$ and a family of sets $(\strC_{p, 0}') _{p \geq p_2},$ such that $\strC_{p, 0}' \subset \strC_{p, 0},$ $\card{\strC_{p, 0}'} \leq \varepsilon \cdot \card{\strC_{p, 0}}$ and $\sup\limits_{p \geq p_2} \max\limits_{z \in \strC_{p, 0}} \norm{z - \strC_{p, 0}'} < \infty.$
		\end{Proposition}
		The proof follows from elementary techniques and can be found in \cite[Proposition (5.2.3)]{MaDib:USFwithDrift}.
	
		For the next proposition, the reader should recall the definition of $\LE.$
	
		\begin{Proposition}\label{Proposition: Lower bound in the probability of intersection of LERW and RW}
			Let $S$ and $S'$ be two independent transient irreducible Markov chains on the same countable state space (having same transition density), with initial states $o$ and $o',$ respectively. Let $\strU$ be any subset of the states containing $o$ and $o'.$ Denote by $T$ and $T'$ the exit times of $\strU$ by $S$ and $S',$ respectively. If the event $\EuScript{H}_\strU$ defined by ``there exists $0 \leq m \leq T$ and $0 \leq m' \leq T'$ with $S_m = S_{m'}'$'' has positive probability, then $\Probability{\LE(S_{m'}')_{m' = 0, \ldots, T'} \cap (S_m)_{m = 0, \ldots, T} \neq \varnothing} \geq 2^{-8} \Probability{\EuScript{H}_\strU}.$ (This inequality obviously also holds if $\EuScript{H}_\strU$ has probability zero.)
		\end{Proposition}
		We only provide a reference for this proposition: see Remark 1.3 and Lemma 4.1 of \cite{LyPeSch:MCILERW}.

		In (\ref{Theorem: Logarithmic scale of crossings}) we proved a \emph{lower bound} for the probability that two independent sections of the network random walk intersect. That was an involved theorem. The following proposition provides an \emph{upper bound} for the probability that two network random walks intersect in a \emph{narrow} strip.
	
		\begin{Proposition}\label{Proposition: Probability of intersection of two random walks in the strips}
			Let $d = 2.$ Consider $k_0$ and $p_0$ as in (\ref{Proposition: Asymptotic bound between the Green function and the restricted Green function on d equals 2}) with $\varepsilon = 1$ and define the ``strips''
			\begin{displaymath}\textstyle
				\strJ_p^- = [4 \cdot 9^{p k_0} - 2 \cdot 3^{p k_0}, 4 \cdot 9^{p k_0}] \times \Z^2, \quad \strJ_p^+ = [4 \cdot 9^{p k_0}, 4 \cdot 9^{p k_0} + 2 \cdot 3^{p k_0}] \times \Z^2.
			\end{displaymath}
			With $S$ and $S'$ being two independent network random walks in $\Gamma_d(\lambda)$ started at $o$ and $0,$ respectively, denote by $I_p^\pm(o)$ the number of intersections of the paths $S$ and $S'$ inside the strip $\strJ_p^\pm.$ Then, there exists a constant $c > 0$ such that for all $\delta > 0,$ there exists an index $p_3 = p_3(\delta)$ such that, for all $p \geq p_3,$ $\max\limits_{\norm{o} \leq \delta} \Probability{I_p^\pm(o) > 0} \leq c 3^{-p k_0}.$
		\end{Proposition}
		\begin{proof}
			Take $p_{3, 1} = p_{3, 1}(\delta)$ the first integer $s$ such that the ball $\ball{0; \delta + 1}$ is a subset of $\strU_{s - 1}$ defined by (\ref{Equation: Region of restriction for the Green function}). Assume $p \geq p_{3, 1}.$ Write $I_p^\pm$ instead of $I_p^\pm(o)$ for simplicity and we will show the bounds are universal, provided $p$ is large enough (depending on $\delta$); since the argument is very similar, suffices to consider the case of $I_p^+.$ We have $I_p^+ = \sum\limits_{z \in \strJ_p^+} \sum\limits_{m, m' = 0}^\infty \indic{\{S_m = z\}} \indic{\{S_{m'}' = z\}},$ and by independence, $\Probability{I_p^+ > 0} \leq \Expectation{I_p^+} = \sum\limits_{z \in \strJ_p^+} G(o, z) G(0, z).$ Write now
			\begin{displaymath}\textstyle
				\begin{split}
					\strJ_{p, 0}^+ &= [4 \cdot 9^{p k_0}, 4 \cdot 9^{p k_0} + 2 \cdot 3^{p k_0}] \times \ball{0; 3^{p k_0}} \\
					\strJ_{p, t}^+ &= [4 \cdot 9^{p k_0}, 4 \cdot 9^{p k_0} + 2 \cdot 3^{p k_0}] \times \left( \ball{0; 3^{(t + 1) p k_0}} \setminus \ball{0; 3^{t p k_0}} \right) \quad (t \in \N),
				\end{split}
			\end{displaymath}
			so that $\strJ_p^+$ is the union of the pairwise disjoint sets $\strJ_{p, t}^+$ ($t \in \Z_+$).
		
			It easily follows that $G(o, z) \asymp G(0, z)$ for $z \in \strJ_p^+$ and all $p$ large enough with any implicit constant universal. Denote now $p_3 = \max(p_{3, 1}, p_{3, 2}, p_{3, 3}).$ There exists $c > 0$ such that $G(o, z) \leq c G(0, z).$ Substitute this into the equation above, $\Probability{I_p^+ > 0} \leq \sum\limits_{z \in \strJ_p^+} G(o, z) G(0, z) \leq c \sum\limits_{z \in \strJ_p^+} G(0, z)^2 = c \sum\limits_{t \in \Z_+} \sum\limits_{z \in \strJ_{p, t}^+} G(0, z)^2.$ For $t = 0,$ by (\ref{Theorem: Greens function bounds}), $\sum\limits_{z \in \strJ_{p, 0}^+} G(0, z) \leq c_1^2 9^{-2p k_0} \card{\strJ_{p, 0}^+} \asymp 3^{-4p k_0} \times 3^{p k_0} \cdot 3^{2p k_0} = 3^{-p k_0}.$ For $t = 1,$ and with $\card{\ball{x; r}} \leq L r^2$ for $x \in \Z^2$ and $r \geq 1,$ we have
			\begin{align*}\textstyle
				\sum\limits_{z \in \strJ_{p, 1}^+} G(0, z)^2 &\textstyle\leq c_1^2 \sum\limits_{n = 4 \cdot 9^{p k_0}}^{4 \cdot 9^{p k_0} + 2 \cdot 3^{p k_0}} \sum\limits_{r = 1}^{2 \cdot 3^{p k_0}} \sum\limits_{r 3^{p k_0} \leq \norm{x} \leq (r + 1) 3^{p k_0}} e^{-2 c_2 \frac{\norm{x}^2}{n}} n^{-2} \\
				&\textstyle\leq  2c_1^2 9^{-2 p k_0} \cdot 3^{p k_0} \sum\limits_{r = 1}^{2 \cdot 3^{p k_0}} L(r + 1)^2 9^{p k_0} e^{-c_2 r^2} \leq \left[ 2 c_1^2 L \sum\limits_{r = 1}^\infty (r + 1)^2 e^{-c_2 r^2} \right] 3^{-p k_0}.
			\end{align*}
			For $t \geq 2$ and $z = (n_z, x_z) \in \strJ_{p, t}^+,$ the condition $3^{t p k_0} \leq \norm{x_z} \leq 3^{(t + 1) p k_0}$ guarantees the bound $G(0, z) \leq c_1 e^{-c_2 \norm{z}},$ and we obtain bounds that are negligible compared with the two bounds already obtained. Indeed,
			\begin{displaymath}\textstyle
				\sum\limits_{t \geq 2} \sum\limits_{z \in \strJ_{p, 0}^+} G(0, z)^2 \textstyle\leq c_1^2 L \sum\limits_{t \geq 2} \sum\limits_{n = 4 \cdot 9^{p k_0}}^{4 \cdot 9^{p k_0} + 2 \cdot 3^{p k_0}} e^{-c_2 3^{t p k_0}} 3^{2(t + 1) p k_0} \leq \left[ 2 c_1^2 L \sum\limits_{t \geq 2} e^{-\frac{c_2}{2} 3^{t p k_0}} 3^{(2t + 3) p k_0} \right] e^{-\frac{c_2}{2} 3^{2 p k_0}}.
			\end{displaymath}
			This completes the proof.
		\end{proof}
	
		Before continuing with the next proposition, we introduce the following terminology. A set $\strH$ of the form $\strH = \{h\} \times \Z^d$ will be called a \textbf{splitting hyperplane} (\textbf{splitting plane} if $d = 2$) for the path  $(v_n)$ (either finite or infinite) if for one, and exactly one value of $m,$ we have $\pr{0}v_m = h.$ In this case, we will call $h$ to be a \textbf{splitting level} for this path $(v_m).$
	
		\begin{Proposition}\label{Proposition: Juxtaposition of random walks given splitting levels}
			Let $S$ be the network random walk on $\Gamma_d(\lambda)$ started at $S_0.$ Suppose that $\zeta \leq \eta$ are two stopping times for $S,$ with $\zeta < \infty$ almost surely ($\eta = \infty$ with positive probability not excluded). Denote by $\EuScript{S}_h(\zeta, \eta)$ the event that $h \in \Z$ is a splitting level for the path $S = (S_m)_{\zeta \leq m \leq \eta}$ (whenever $\eta = \infty,$ replace $m \leq \eta$ by $m < \eta$). On $\EuScript{S}_h(\zeta, \eta),$ denote by $m_0$ to be the unique positive integer $m$ for which $\pr{0} S_m = h.$ Then, on the event $\EuScript{S}_h(\zeta, \eta),$ $\LE(S) = \LE(S_m)_{\zeta \leq m \leq m_0} \vee \LE(S_m)_{m_0 \leq m \leq \eta}.$
		\end{Proposition}
		The proof is very easy and therefore, omitted.
	
		The next proposition shows that splitting levels do occur with a universal overwhelming probability if one let appropriate scales to be considered. In other words, if the random walk moves enough so that the initial and last state considered are $n$ levels apart,  \emph{at the very least}, one in roughly every $\sqrt{n}$ consecutive levels will be splitting. The gist of its proof lies in random walk in one dimension, see (\ref{Proposition: Probability for the existence of splitting levels}).
	
		\begin{Proposition}\label{Proposition: Lower bound in the probability of existence of splitting hyperplanes}
			Consider the network random walk $S$ on $\Gamma_d(\lambda)$ started at zero. For every positive integer $a,$ every pair $\zeta_p < \eta_p$ of stopping times for $S$ such that, for almost every realisation, $\zeta_p \leq \tau_{\{a\} \times \Z^2}(S),$ $\tau_{\{a + p^2\} \times \Z^2}(S) \leq \eta_p,$ and for $h \in \Z_+,$ consider the event $\EuScript{S}_h(\zeta_p, \eta_p)$ that $h$ is a splitting level for $(S_m)_{m = \zeta_p, \ldots, \eta_p}.$ There exists a universal constant $c > 0$ and an index $p_4 \in \N$ such that, for all $p \geq p_4,$ all systems $(a, \zeta_p, \eta_p)$ as above and all integers $b$ satisfying $a \leq b$ and $b + p \leq a + p^2,$ we have $\Probability{\bigcup\limits_{h = b}^{b + p} \EuScript{S}_h(\zeta_p, \eta_p)} \geq 1 - e^{-c p}.$
		\end{Proposition}
		\begin{proof}
			Let $\EuScript{K}_p$ be the event ``some $h = b, \ldots, b + p$ is a splitting level for the path $(S_m)_{m \in \Z_+}$.'' We have $\EuScript{K}_p \subset \bigcup\limits_{h = b}^{b + p} \EuScript{S}_h(\zeta_p, \eta_p)$ and therefore $\Probability{\bigcup\limits_{h = b}^{b + p} \EuScript{S}_h(\zeta_p, \eta_p)} \geq \Probability{\EuScript{K}_p}.$ The result follows at once from (\ref{Proposition: Probability for the existence of splitting levels}).
		\end{proof}
	
		\begin{Proposition}\label{Proposition: Probability of intersection of two sections of random walk}
			Let $d = 2.$ Consider $S$ and $S'$ two independent network random walks in $\Gamma_2(\lambda)$ started at $o$ and $0,$ respectively. Define the following stopping time $\sigma_p = \tau_{\strH_p}(S),$ where $\strH_p = \{4 \cdot 9^{p k_0}\} \times \Z^2,$ in other words, $\sigma_p$ is the hitting time of the ``plane'' $\strH_p$ by $S.$ Define $\sigma_p'$ similarly, using $S'$ in lieu of $S.$ Denote $L' = \LE(S_{m'}')_{m' = \sigma_{2(p - 1)}', \ldots, \sigma_{2p}'},$ which is the loop-erasure of the section of $S'$ between the planes $\strH_{2(p - 2)}$ and $\strH_{2p};$ with this, define $\beta_p(o) = \sum\limits_{m = \sigma_{2(p - 1)}}^{\sigma_{2p}} \sum\limits_{z' \in L'} \indic{\{S_m = z\}},$ which is the number of times the section $(S_m)_{m = \sigma_{2(p - 1)}, \ldots, \sigma_{2p}}$ crosses the path $L'.$ Then, there exists a constant $c > 0$ such that for every $\delta > 0,$ there exists a $p_5 = p_5(\delta) \in \N$ satisfying $\min\limits_{\norm{o} \leq \delta} \Probability{\beta_p(o) > 0} \geq \frac{c}{p},$ for all $p \geq p_5.$
		\end{Proposition}
		\begin{proof}
			Define $M_p(o)$ as in (\ref{Equation: Definition of the number of crossing in the sections}) using $S$ and $S'.$ Define $p_{5, 1} = p_{5, 1}(\delta) \in \N$ to be the first index $p$ such that $\ball{0; \delta + 1} \subset \strU_{p - 1}.$ Notice that for $p \geq p_{5, 1},$ $T_p \leq \sigma_p$ and $T_p' \leq \sigma_p'.$ Since $\sigma_{p - 1}$ and $\sigma_{p - 1}'$ are the hitting times of the plane $\strH_{p - 1},$ we find that
			\begin{displaymath}\textstyle
				M_p(o) = \sum\limits_{z \in \strD_{p k_0}} \sum\limits_{m = \sigma_{p - 1}}^{T_p} \sum\limits_{m' = \sigma_{p - 1}'}^{T_p'} \indic{\{S_m = z\}} \indic{\{S_{m'}' = z\}} \leq \sum\limits_{z \in \strD_{p k_0}} \sum\limits_{m = \sigma_{p - 1}}^{\sigma_p} \sum\limits_{m' = \sigma_{p - 1}'}^{\sigma_p'} \indic{\{S_m = z\}} \indic{\{S_{m'}' = z\}} \mathop{=}\limits^{\mathrm{def.}} N_p(o).
			\end{displaymath}
			 Thus, $\Probability{M_p(o) > 0} \leq \Probability{N_p(o) > 0}$ and Theorem (\ref{Theorem: Logarithmic scale of crossings}) gives $\Probability{M_p(o) > 0} \geq \frac{c}{p},$ provided $p \geq P(o, 0),$ the constant $c$ being universal. Define $\ds p_{5, 2} = p_{5, 2}(\delta) = \max_{\norm{o} \leq \delta} P(o, 0).$ Finally, (\ref{Proposition: Lower bound in the probability of intersection of LERW and RW}) gives $\Probability{\beta_p(o) > 0} \geq 2^{-8} \Probability{\tilde{N}_p(o) > 0},$ where $\tilde{N}_p(o)$ is the number of intersection between the two sections $(S_m)_{m = \sigma_{2(p - 1)}, \ldots, \sigma_{2p}}$ and $(S_{m'}')_{m' = \sigma_{2(p - 1)}', \ldots, \sigma_{2p}'}.$ Clearly $\tilde{N}_p(o) \geq N_p(o),$ and the result follows with $p_5 = \max(p_{5, 1}, p_{5, 2}).$
		\end{proof}
	
		\begin{Proposition}\label{Proposition: Limit probability that the LERW and a section of SRW will not intersect}
			Let $d = 2.$ With the hypotheses and notations of (\ref{Proposition: Probability of intersection of two sections of random walk}). Define $\alpha_p(o)$ to be the number of intersections of $\LE(S_{m'}')_{m' = 0, \ldots, \sigma_{2p}'}$ and $(S_m)_{m = 0, \ldots, \sigma_{2p}}.$ For every $\delta > 0,$ $\lim\limits_{p \to \infty} \max\limits_{\norm{o} \leq \delta} \Probability{\alpha_p(o) = 0} = 0.$
		\end{Proposition}
		\setcounter{auxCounter}{\value{equation}}
		\setcounter{equation}{0}
		\renewcommand{\theequation}{\Alph{equation}}
		\renewcommand{\theSubCounter}{\thesection.\theauxCounter.\arabic{equation}}
		\begin{proof}
			Let $\varepsilon > 0.$ We can then find the indices $p_1(\varepsilon)$ and $p_2(\varepsilon),$ and a constant $b_1(\varepsilon)$ satisfying the conclusions of (\ref{Proposition: Sum of the Green function on cylinder can be made small}) and (\ref{Proposition: Existence of sparse good sets}), we take $b = b_1(\varepsilon).$ The $\delta$ in the hypothesis, together with (\ref{Proposition: Probability of intersection of two random walks in the strips}) and (\ref{Proposition: Probability of intersection of two sections of random walk}) give two indices $p_3(\delta)$ and $p_5(\delta),$ satisfying the corresponding conclusions of these results. And (\ref{Proposition: Lower bound in the probability of existence of splitting hyperplanes}) gives yet another index $p_4$ satisfying the outcome of that proposition. Observe that the choice of $p_3(\delta)$ guarantees that the closed ball $\ball{0; \delta + 1}$ is contained in $\strU_{p_3(\delta) - 1}.$ Let $P = \max(p_1(\varepsilon), p_2(\varepsilon), p_3(\delta), p_4, p_5(\delta)).$
		
			Write $\alpha_p$ in place of $\alpha_p(o)$ for simplicity, all the bounds will be shown to be independent of $\delta.$ Assume $p \geq P.$ Define the events
			\begin{equation}\textstyle\label{EqInside: Definition of the left and right splitting levels}
				\begin{split}
					\EuScript{L}_p' &\textstyle= \bigcup\limits_{h = 4 \cdot 9^{2(p - 1) k_0}}^{4 \cdot 9^{2(p - 1) k_0} + 2 \cdot 9^{(p - 1) k_0}} \left\{ h\text{ is a splitting level for } \left( S'_{m'} \right)_{m' = \sigma_{2(p - 1)}', \ldots, \sigma_{2p}'} \right\} \\
					\EuScript{R}_p' &\textstyle= \bigcup\limits_{h = 4 \cdot 9^{2p k_0} - 2 \cdot 9^{p k_0}}^{4 \cdot 9^{2p k_0}} \left\{ h\text{ is a splitting level for } \left( S'_{m'} \right)_{m' = \sigma_{2(p - 1)}', \ldots, \sigma_{2p}'} \right\},
				\end{split}
			\end{equation}
			that is, the events where one of the $2 \cdot 9^{(p - 1) k_0}$ left-most levels between the planes $\strH_{2(p - 1)}$ and $\strH_{2p}$ and one of the $2 \cdot 9^{p k_0}$ right-most levels in the same region is a splitting level for the stated section of $S'.$ Denote by $\EuScript{R}_p'^\complement$ and $\EuScript{L}_p'^\complement$ their corresponding complements; we also define $\beta_p = \beta_p(o)$ from (\ref{Proposition: Probability of intersection of two sections of random walk}). Then,
			\begin{align*}
				\Probability{\alpha_p = 0} &= \Probability{\{\alpha_p = 0\} \cap \EuScript{L}_p'} + \Probability{\{\alpha_p = 0\} \cap \EuScript{L}_p'^\complement} \\
				&\leq \Probability{\{\alpha_p = 0\} \cap \EuScript{L}_p' \cap \{\beta_p = 0\}} + \Probability{\{\alpha_p = 0\} \cap \EuScript{L}_p' \cap \{\beta_p > 0\}} + \Probability{\EuScript{L}_p'^\complement}.
			\end{align*}
			We analyse now what happens on the intersection $\{\alpha_p = 0\} \cap \EuScript{L}_p' \cap \{\beta_p > 0\}.$ Assume that this event occurs. Then, there exists a splitting level $h \in \left\{ 4 \cdot 9^{2(p - 1)k_0}, \ldots, 4 \cdot 9^{2(p - 1)k_0} + 2 \cdot 9^{(p - 1)k_0} \right\}$ for the section $(S'_{m'})_{m' = \sigma_{2(p - 1)}', \ldots, \sigma_{2p}'}.$ Denote by $h_p'$ the first of such splitting levels and by $m_p'$ the unique index for which $\pr{0} S_{m_p'}' = h_p'.$ According to (\ref{Proposition: Juxtaposition of random walks given splitting levels}), we have
			\begin{equation}\textstyle\label{EqInside: Crutial splitting partial section}
				\LE \left( S'_{m'} \right)_{m' = \sigma_{2(p - 1)}', \ldots, \sigma_{2p}'} = \LE(S'_{m'})_{m' = \sigma_{2(p - 1)}', \ldots, m_p'} \vee \LE(S'_{m'})_{m' = m_p', \ldots, \sigma_{2p}'}.
			\end{equation}
			\noindent Since $\sigma_{2(p - 1)}'$ is the hitting time of the set $\strH_{2(p - 1)},$ it is then obvious that $h_p'$ is also a splitting level for the section $(S_{m'}')_{m' = 0, \ldots, \sigma_{2p}'}$ and so, we also have
			\begin{equation}\textstyle\label{EqInside: Crutial splitting whole section}
				\LE(S'_{m'})_{m' = 0, \ldots, \sigma_{2p}'} = \LE(S'_{m'})_{m' =0, \ldots, m_p'} \vee \LE(S'_{m'})_{m' = m_p', \ldots, \sigma_{2p}'}.
			\end{equation}
			\noindent That $\alpha_p = 0$ entails that the section $(S_m)_{m = 0, \ldots, \sigma_{2p}}$ does not intersect the set in (\ref{EqInside: Crutial splitting whole section}) and the relation $\beta_p > 0$ means that it does intersect the set in (\ref{EqInside: Crutial splitting partial section}). Therefore,  on the event under consideration, the intersection of $(S_m)_{m = 0, \ldots, \sigma_{2p}}$ with $\LE(S'_{m'})_{m' = \sigma_{2(p - 1)}', \ldots, \sigma_{2p}'}$ occurs at some $S_{m'}'$ with $m' \in \left\{ \sigma_{2(p - 1)}',\ldots, m_p' \right\};$ in particular, $S$ and $S'$ intersect inside the strip $\strJ_{2(p - 1)}^+$ (notation of (\ref{Proposition: Probability of intersection of two random walks in the strips})). In other words, with the notation $I_p^+$ of (\ref{Proposition: Probability of intersection of two random walks in the strips}), $\Probability{\{\alpha_p = 0\} \cap \EuScript{L}_p' \cap \{\beta_p > 0\}} \leq \Probability{I_{2(p - 1)}^+ > 0}.$ Whence, by the arguments of the previous two paragraphs,
			\begin{equation}\textstyle\label{Equation: First split of the probability alpha = 0}
				\Probability{\alpha_p = 0} \leq \Probability{ \{\alpha_p = 0\} \cap \EuScript{L}_p' \cap \{\beta_p = 0\} } + \Probability{\EuScript{L}_p'^\complement} + \Probability{I_{2(p - 1)}^+ > 0}.
			\end{equation}
			Now observe that
			\begin{align*}
				\{\alpha_p = 0\} \cap \EuScript{L}_p' \cap \{\beta_p = 0\} &= \big( \{\alpha_p = 0, \alpha_{p - 1} = 0\} \cap \EuScript{L}_p' \cap \{\beta_p = 0\} \big) \\
				&\cup \big( \{\alpha_p = 0, \alpha_{p - 1} > 0\} \cap \EuScript{L}_p' \cap \{\beta_p = 0\} \big) \\
				&\subset \big( \{\alpha_{p - 1} = 0\} \cap \EuScript{L}_p' \cap \{\beta_p = 0\} \big) \cup \{\alpha_p = 0, \alpha_{p - 1} > 0\},
			\end{align*}
			and from this it follows that
			\begin{displaymath}\textstyle
				\Probability{\{\alpha_p = 0\} \cap \EuScript{L}_p' \cap \{\beta_p = 0\}} \leq \Probability{ \{\alpha_{p - 1} = 0\} \cap \EuScript{L}_p' \cap \{\beta_p = 0\} } + \Probability{\alpha_p = 0, \alpha_{p - 1} > 0}.
			\end{displaymath}
			By (\ref{Proposition: Independence between past and future of two processes with stationary independent increments a stopping times}), $\{\alpha_{p - 1}= 0\}$ and $\EuScript{L}_p' \cap \{\beta_p = 0\}$ are independent events. Write now
			\begin{displaymath}\textstyle
				\{\alpha_p = 0, \alpha_{p - 1} > 0\} = \big( \{\alpha_p = 0, \alpha_{p - 1} > 0\} \cap \EuScript{R}_{p - 1}' \big) \cup \left( \left\{ \alpha_p = 0, \alpha_{p - 1} > 0 \right\} \cap \EuScript{R}_{p - 1}'^\complement \right);
			\end{displaymath}
			as done before, on the set $\{\alpha_p = 0, \alpha_{p - 1} > 0\} \cap \EuScript{R}_{p - 1}',$ the two random walks $S$ and $S'$ intersect in the strip $\strJ_{2(p - 1)}^-,$ yielding $\Probability{\{\alpha_p = 0, \alpha_{p - 1} > 0\} \cap \EuScript{R}_{p - 1}'} \leq \Probability{I_{2(p - 1)}^- > 0}.$ With the aforementioned considerations and substituting into (\ref{Equation: First split of the probability alpha = 0}), we reach
			\begin{align*}
				\Probability{\alpha_p = 0} &\leq \Probability{\alpha_{p - 1} = 0} \Probability{\beta_p = 0} \\
				&+ \Probability{I_{2(p - 1)}^- > 0}+\Probability{\EuScript{R}_{p - 1}'^\complement}+ \Probability{\EuScript{L}_p'^\complement} + \Probability{I_{2(p - 1)}^+ > 0}.
			\end{align*}
			By virtue of propositions (\ref{Proposition: Probability of intersection of two random walks in the strips}) and (\ref{Proposition: Lower bound in the probability of existence of splitting hyperplanes}), we have $\Probability{I_{2(p - 1)}^- > 0}+\Probability{\EuScript{R}_{p - 1}'^\complement}+ \Probability{\EuScript{L}_p'^\complement} + \Probability{I_{2(p - 1)}^+ > 0} \leq c e^{-c' p}$ and by (\ref{Proposition: Probability of intersection of two sections of random walk}) $\Probability{\beta_p = 0} \leq 1 - \frac{c''}{p},$ for $p \geq P,$ where the constants $c, c'$ and $c''$ are all universal. We can conclude then that $\Probability{\alpha_{P + p} = 0} \leq \prod\limits_{q = 0}^p \left(1 - \frac{c''}{P + q} \right) + c e^{-c' P}\sum\limits_{q = 0}^p e^{-c' q}.$ The conclusion of the proposition is now clear given the last inequality and bearing in mind that all estimates are uniform in $\norm{o} \leq \delta.$
		\end{proof}
	\numberwithin{SubCounter}{equation}
	\numberwithin{equation}{section}
	\setcounter{equation}{\value{auxCounter}}
	
		\begin{Proposition}\label{Proposition: Limit probability that the LERW and SRW will not intersect}
			Let $d = 2.$ With the hypotheses of (\ref{Proposition: Limit probability that the LERW and a section of SRW will not intersect}), let $\tilde\alpha_p(o)$ be the number of intersections of $\LE(S_{m'}')_{m' \in \Z_+}$ and $(S_m)_{m = 0, \ldots, \sigma_{2p}}.$ Then, for every $\delta > 0,$ $\lim_{p \to 0} \max\limits_{\norm{o} \leq \delta} \Probability{\tilde\alpha_p(o) = 0} = 0.$
		\end{Proposition}
		\begin{proof}
			We use the notation of (\ref{Proposition: Lower bound in the probability of existence of splitting hyperplanes}). Simply write
			\begin{align*}
				\Probability{\tilde\alpha_p(o) = 0} &= \Probability{\tilde\alpha_p(o) = 0, \alpha_p(o) = 0} + \Probability{\tilde\alpha_p(o) = 0, \alpha_p(o) > 0} \\
				&\leq \Probability{\alpha_p(o) = 0} + \Probability{ \left\{ \tilde\alpha_p(o) = 0, \alpha_p(o) > 0 \right\} \cap \EuScript{R}_p} \\
				&+ \Probability{ \left\{ \tilde\alpha_p(o) = 0, \alpha_p(o) > 0 \right\} \cap \EuScript{R}_p^\complement }.
			\end{align*}
			Exactly as in the proof of (\ref{Proposition: Limit probability that the LERW and a section of SRW will not intersect}), the two relations $\tilde\alpha_p(o) = 0$ and $\alpha_p(o) > 0$ on the event $\EuScript{R}_p$ imply that two independent random walks will intersect in the strip $\strJ_{2p}^-,$ hence the result follows immediately upon invoking (\ref{Proposition: Probability of intersection of two random walks in the strips}), (\ref{Proposition: Lower bound in the probability of existence of splitting hyperplanes}) and (\ref{Proposition: Limit probability that the LERW and a section of SRW will not intersect}).
		\end{proof}

		\subsubsection{The proof of the main theorem.}
		\begin{Theorem}\label{Theorem: One end in dimension two}
			For almost every realisation of $\UST$ in $\Gamma_2(\lambda),$ the tree has one end.
		\end{Theorem}
		\setcounter{auxCounter}{\value{equation}}
		\setcounter{auxSubCounter}{\value{SubCounter}}
		\setcounter{equation}{0}
		\setcounter{SubCounter}{0}
		\renewcommand{\theequation}{\Alph{equation}}
		\renewcommand{\theSubCounter}{\thesection.\theauxCounter.\arabic{SubCounter}}
		\begin{proof}
			Denote by $\UST$ the measure of uniform spanning tree on the network $\Gamma_2(\lambda).$ When we define an event on spanning trees, we use $\UST$ to refer to the tree itself. We know that $\UST$ has at least one end, up to a $\UST$-null event. Consider $\scrN = \{\UST\text{ has at least two ends}\}.$ We aim at showing $\UST(\scrN) = 0.$ For a vertex $z \in \Z^3,$ consider $\scrN_z = \{\text{there are two disjoint rays starting at }z\text{ in }\UST\}.$ Then, $\scrN = \bigcup\limits_{z \in \Z^3} \scrN_z.$ Further, the translation invariance of the probability kernel of the network random walk of $\Gamma_2(\lambda)$ shows that $\UST(\scrN_z)$ is a value independent of $z \in \Z^3.$ Hence, \emph{to show $\scrN$ is $\UST$-null, it suffices to show $\scrN_0$ is $\UST$-null.}
	
			Consider now a set of vertices $\strK$ of $\Z^3$ with the following property (a ``cutset'' between $0$ and infinity): every ray from $0$ must cross $\strK$ at some vertex. Recall an ``edge of $\strK$'' is an edge of the subgraph induced on $\strK.$ Then, $\scrN_0$ is contained in the event
			\begin{displaymath}\textstyle
				\scrC_\strK = \left\{\text{there are two disjoint paths in }\UST\text{ from }\strK\text{ to 0 that use no edge of }\strK \right\}.
			\end{displaymath}
			In particular, we have the following: $\UST(\scrN_0) \leq \inf\limits_{\substack{\strK\text{ is a cutset}\\\text{between }0\text{ and }\infty}} \UST(\scrC_\strK).$
	
			As in proposition (\ref{Proposition: Sum of the Green function on cylinder can be made small}), define the sets $\strA_p = [-p, p] \times \sphere{0; p},$ $\strB_p =\{p\} \times \ball{0; p}$ and $\strC_p =\{-p\} \times \ball{0;p}.$ Their union $\strK_p$ is a cutset between $0$ and $\infty.$ In particular,
			\begin{equation}\textstyle\label{SubEquation: Upper bound in the probability of N0}
				\UST(\scrN_0) \leq \inf\limits_{p \in \N} \UST(\scrC_{\strK_p}).
			\end{equation}
			Construct $\mathfrak{T}$ the $\UST$ of $\Z^3$ using Wilson's algorithm rooted at infinity and following some starting at zero predefined ordering of the vertices of $\Z^3,$ and assume it is defined on some probability space $(\Omega, \scrF, \mathds{P}).$ (For instance, we can follow the order in which we start with $0,$ and having searched all vertices such that $\norm{z}_1 = k,$ we search the vertices satisfying $\norm{z}_1 = k + 1$ arranging them lexicographically.) Thus, we have a family of independent network random walks $(S^z)_{z \in \Z^3}$ defined on this probability space with $S^z$ started at $z.$ Since $0$ is the first vertex searched, $L_0 = \LE(S^0_m)_{m \in \Z_+} \subset \mathfrak{T}$ (both trees canonically identified with their sets of edges). Up to a $\mathds{P}$-null event, $L_0$ is an infinite branch. Thus, $\{\mathfrak{T} \in \scrC_{\strK_p}\} = \{\text{some vertex of }\strK_p\text{ is connected to $0$ in }\mathfrak{T} \setminus L_0\}.$ Denote $\EuScript{C}(z) = \{z\text{ is connected to $0$ in }\mathfrak{T} \setminus L_0\},$ so that $\{\mathfrak{T} \in \scrC_{\strK_p}\} = \bigcup\limits_{z \in \strK_p} \EuScript{C}(z).$ By virtue of (\ref{Proposition: Sum of the Green function on cylinder can be made small}) we reach the existence of an index $p_1(\varepsilon)$ and a number $b_1(\varepsilon)$ such that we may divide $\strC_p = \strC_{p, 0} \cup \strC_{p, 1}$ as in this proposition (taking $b = b_1(\varepsilon)$) and for all $p \geq p_1(\varepsilon),$ $\sum\limits_{z \in \strE_p} G(z, 0) \leq \varepsilon,$ 	with the notation $\strE_p = \strA_p \cup \strB_p \cup \strC_{p, 1}$ of the proposition (in particular, $G$ is the Green's function of $\Gamma_2(\lambda)$). Next,
			\begin{displaymath}\textstyle
				\UST(\scrC_{\strK_p}) = \Probability{\mathfrak{T} \in \scrC_{\strK_p}} \leq \sum\limits_{z \in \strE_p} \Probability{\EuScript{C}(z)} + \Probability{\bigcup\limits_{z \in \strC_{p, 0}} \EuScript{C}(z)}.
			\end{displaymath}
			By corollary (\ref{Corollary: Main extension of WA}), when calculating $\Probability{\EuScript{C}(z)}$ we may assume that the order in which $\mathfrak{T}$ was constructed is $(0, z, \ldots).$ If the order in the construction of $\mathfrak{T}$ were $(0, z, \ldots),$ then $\EuScript{C}(z)$ would be the event where $S^z$ hits $L_0$ for the first time at zero, and then $\EuScript{C}(z)$ would be contained in the event where $\tau_0(S^z) < \infty.$ Then, $\Probability{\EuScript{C}(z)} \leq G(z, 0).$ Substituting this inequality above we obtain that, for all $p \geq p_1(\varepsilon),$
			\begin{equation}\textstyle\label{SubEquation: First bound in the probability of CK}
				\UST(\scrC_{\strK_p}) \leq \varepsilon + \Probability{\bigcup\limits_{z \in \strC_{p, 0}} \EuScript{C}(z)}.
			\end{equation}
	
			Using (\ref{Proposition: Existence of sparse good sets}), there is an index $p_2(\varepsilon)$ and a family of sets $(\strC_{p, 0}')_{p \geq p_2(\varepsilon)},$ for which the conclusions of this proposition hold. In particular, we may define $\delta_0 < \infty$ as follows
			\begin{displaymath}\textstyle
				\delta_0 = \sup\limits_{p \geq p_2(\varepsilon)} \max\limits_{z \in \strC_{p, 0}} \norm{z - \strC_{p, 0}'}.
			\end{displaymath}
			For each $z' \in \strC_{p, 0}',$ consider the set $\strC_{p, 0}(z') = \{-p\} \times \ball{z'; \delta_0}.$ The definition of $\delta_0$ shows at once $\strC_{p, 0} \subset \bigcup\limits_{z' \in \strC_{p, 0}'} \strC_{p, 0}(z').$ Then, $\Probability{\bigcup\limits_{z \in \strC_{p, 0}} \EuScript{C}(z)} \leq \sum\limits_{z' \in \strC_{p, 0}'} \Probability{\bigcup\limits_{z \in \strC_{p, 0}(z')} \EuScript{C}(z)}.$ By definition, $\bigcup\limits_{z \in \strC_{p, 0}(z')} \EuScript{C}(z)$ is an event depending solely on $\mathfrak{T}$ and $L_0.$ Thus, (\ref{Corollary: Main extension of WA}) shows that its probability is independent of the ordering of the vertices of $\Z^3 \setminus \{0\},$ in particular, we may assume that $\mathfrak{T}$ was constructed using the ordering $(0, z', \ldots).$ Consider now the following event $\EuScript{H}(z') = \left\{ S^{z'}_{ \tau_{L_0} \left( S^{z'} \right) } \neq 0 \right\},$ in other words, $\EuScript{H}(z')$ is the event where $S^{z'}$ hits $L_0$ for the first time anywhere except zero (recall $S^{z'}$ will hit $L_0$ a.s.), which is the same as the event where $z'$ is not connected to $0$ in $\mathfrak{T} \setminus L_0$ (with the assumed order $(0, z', \ldots)$ of $\Z^3$). Then, $1 - \Probability{\EuScript{H}(z')} = \Probability{ S^{z'}_{ \tau_{L_0} \left( S^{z'} \right) } = 0 } \leq \Probability{\tau_0 \left( S^{z'} \right) < \infty} \leq G(z', 0).$ Since $z' \in \strC_{p, 0}' \subset \strC_{p, 0} = \{-p\} \times \ball{0; b p^{\frac{1}{2}}},$ we have by theorem (\ref{Theorem: Greens function bounds}), $	G(z', 0) \leq c_1 \norm{z'}^{-1} \leq c_1 p^{-1}.$ Then, $\Probability{\bigcup\limits_{z \in \strC_{p, 0}(z')} \EuScript{C}(z)} \leq \Probability{\bigcup\limits_{z \in \strC_{p, 0}(z')} \EuScript{C}(z) \cap \EuScript{H}(z')} +  c_1 p^{-1}.$ Summing over $z' \in \strC_{p, 0}',$
			\begin{displaymath}\textstyle
				\sum\limits_{z' \in \strC_{p, 0}'} \Probability{\bigcup\limits_{z \in \strC_{p, 0}(z')} \EuScript{C}(z)} \leq \sum\limits_{z' \in \strC_{p, 0}'} \Probability{\bigcup\limits_{z \in \strC_{p, 0}(z')} \EuScript{C}(z) \cap \EuScript{H}(z')} + c_1 p^{-1} \cdot \card{\strC_{p, 0}'}.
			\end{displaymath}
			By (\ref{Proposition: Existence of sparse good sets}), we have that for $p \geq p_2(\varepsilon),$ $\card{\strC_{p, 0}'} \leq cp \cdot \varepsilon,$ where $c$ is a universal constant. By (\ref{SubEquation: First bound in the probability of CK}), we have shown hitherto the existence of a universal constant $c$ such that for all $p \geq \max(p_1(\varepsilon), p_2(\varepsilon)),$
			\begin{equation}\textstyle\label{SubEquation: Refinement in the probability of CK}
				\UST(\scrC_{\strK_p}) \leq c \cdot \varepsilon + \sum\limits_{z' \in \strC_{p, 0}'} \sum\limits_{z \in \strC_{p, 0}(z')} \Probability{\EuScript{C}(z) \cap \EuScript{H}(z')}.
			\end{equation}
	
			By virtue of (\ref{SubEquation: Upper bound in the probability of N0}) and (\ref{SubEquation: Refinement in the probability of CK}), we reduce the proof of the theorem to showing following: \emph{there exists a function $\varphi(p)$ such that
			\begin{equation}\textstyle
				\sum\limits_{z' \in \strC_{p, 0}'} \sum\limits_{z \in \strC_{p, 0}(z')} \Probability{\EuScript{C}(z) \cap \EuScript{H}(z')} \leq \varphi(p)
			\end{equation}
			and $\ds \varphi(p) \to 0$ via some subsequence.}
	
			\vspace{0.25cm}
	
			We begin by estimating $\Probability{\EuScript{C}(z) \cap \EuScript{H}(z')}.$ Again by (\ref{Corollary: Main extension of WA}), to calculate the probability of the event $\EuScript{C}(z) \cap \EuScript{H}(z')$ we may assume that $\mathfrak{T}$ was constructed using the ordering $(0, z', z, \ldots)$ of $\Z^3.$ Assuming this ordering of $\Z^3,$ and on the event $\EuScript{C}(z) \cap \EuScript{H}(z'),$ the first three steps of the construction of $\mathfrak{T}$ proceed as follows:
			\begin{enumerate}\renewcommand{\labelenumi}{\arabic{enumi})}
				\item The first branch at zero is created, we call it $L_0.$
				\item The random walk $S^{z'}$ runs, it hits the first branch at zero at some non zero vertex. Denote by $L_{z'}$ this second branch; by definition, $L_{z'} = \LE \left( S^{z'}_m \right)_{m = 0, \ldots, \tau_{L_0} \left( S^{z'} \right)}$
				\item The random walk $S^z$ runs, and it hit $L_0$ at $0$ before hitting $L_{z'}.$
			\end{enumerate}
	
			Introduce the event $\EuScript{N}(z, z') = \left\{\tau_0(S^z) < \tau_{L_{z'}}(S^z) \right\},$ which is the event where $S^z$ hits $L_0$ at zero before it touches $L_{z'}$ anywhere. The aforementioned first three steps show that
			\begin{equation}\textstyle\label{SubEquation: The event C and H is contained in N}
				\EuScript{C}(z) \cap \EuScript{H}(z') \subset \EuScript{N}(z, z').
			\end{equation}
	
			In what follows, we will consider the planes $\strH_p = \left\{ 4 \cdot 9^{2p k_0} \right\} \times \Z^2,$ and for $z, z' \in \Z^3,$ the hitting times $\sigma_p(z, z') = \tau_{\strH_p + z'}(S^z).$ As such, we will make a ``change of scales'' and consider from now onwards indices of the form $4 \cdot 9^{2p k_0}.$ Assume then that $z' \in \strC_{4 \cdot 9^{2p k_0}, 0}'$ and $z \in \strC_{4 \cdot 9^{2p k_0}, 0}(z').$ Denote $a_p = -4 \cdot 9^{2p k_0} + 4 \cdot 9^{2(p - 1) k_0},$ so that $\strH_{p - 1} + z' = \{a_p\} \times \Z^2.$ Consider the events
			\begin{displaymath}
				\begin{split}
					\EuScript{E}_p &\textstyle= \left\{ S^0\text{ ever reaches }-9^{p k_0} \right\}, \\
					\EuScript{S}_p(z') &\textstyle= \bigcup\limits_{h = a_p - 2 \cdot 9^{(p - 1) k_0}}^{a_p} \left\{ h\text{ is a splitting level for } \left( S^{z'}_{m'} \right)_{m' \in \Z_+} \right\}.
				\end{split}
			\end{displaymath}
			By (\ref{Proposition: Probability for the existence of splitting levels}) and (\ref{Proposition: Lower bound in the probability of existence of splitting hyperplanes}) we know there exists a constant $c > 0$ and an index $p_4 \in \N$ such that for all $p \geq p_4,$ $\Probability{\EuScript{E}_p} = e^{-\lambda 9^{pk_0}}$ and $\Probability{\EuScript{S}_p(z')} \geq 1 - e^{-c 9^{p k_0}}.$ Hence,
			\begin{equation}\textstyle\label{SubEquation: First bound in the probability of N}
				\Probability{\EuScript{N}(z, z')} \leq \Probability{\EuScript{N}(z, z') \cap \EuScript{S}_p(z') \cap \EuScript{E}_p^\complement} + e^{-c 9^{p k_0}} + e^{-\lambda 9^{2p k_0}}.
			\end{equation}
			We introduce the following ``good event'' $\EuScript{G}_p(z') = \EuScript{S}_p(z') \cap \EuScript{E}_p^\complement;$ on this event $S^0$ never goes left of the plane $\{-9^{p k_0}\} \times \Z^2$ and $S^{z'}$ has a splitting level in the stated range. On the good event $\EuScript{G}_p,$ the random walk $S^{z'}$ must have crossed the plane $\strH_{p - 1} + z'$ before hitting the first branch. In particular, $\sigma_{p - 1}(z', z') < \tau_{L_0} \left( S^{z'} \right) < \infty \quad \mathrm{on} \quad \EuScript{G}_p(z').$ Also, on the good event $\EuScript{G}_p,$ there is a splitting plane for $\left( S^{z'}_{m'} \right)_{m' \in \Z_+}$; introduce $h_p(z')$ and $m_p(z') \leq \sigma_{p - 1}(z', z')$ to be the last splitting level (amongst the levels $h = a_p - 2 \cdot 9^{(p - 1) k_0}, \ldots, a_p$)  and the only index $m'$ such that $\pr{0} \left( S^{z'}_{m'} \right) = h_p(z').$ By (\ref{Proposition: Juxtaposition of random walks given splitting levels}), on the good event,
			\begin{displaymath}\textstyle
				L_{z'} = \LE \left( S^{z'}_{m'} \right)_{m' = 0, \ldots, \tau_{L_0 \left( S^{z'} \right)}} = \LE \left( S^{z'}_{m'} \right)_{m' = 0, \ldots, m_p(z')} \vee \LE \left( S^{z'}_{m'} \right)_{m' = m_p(z'), \ldots, \tau_{L_0 \left( S^{z'} \right)}}.
			\end{displaymath}
			Consider now the ``partial branch at $z'$ (until its hitting time of $\strH_{p - 1} + z'$)'' which is
			\begin{displaymath}\textstyle
				L^{z'}_{p - 1} \mathop{=}\limits^{\mathrm{def.}} \LE \left( S^{z'}_{m'} \right)_{m' = 0, \ldots,\sigma_{p - 1}(z', z')}.
			\end{displaymath}
			On $\EuScript{G}_p(z'),$
			\begin{displaymath}\textstyle
				L^{z'}_{p - 1} = \LE \left( S^{z'}_{m'} \right)_{m' = 0, \ldots, \sigma_{p - 1}(z', z')} = \LE \left( S^{z'}_{m'} \right)_{m' = 0, \ldots, m_p(z')} \vee \LE \left( S^{z'}_{m'} \right)_{m' = m_p(z'), \ldots, \sigma_{p - 1}(z', z')}.
			\end{displaymath}
			Define $\alpha_{p - 1}(z, z') = \sum\limits_{m = 0, \ldots, \sigma_{p - 1}(z, z')} \sum\limits_{z'' \in L^{z'}_{p  - 1}} \indic{\{S_m^z = z''\}},$ which is the number of intersections between the partial branch at $z',$ and the section $(S_m^z)_{m = 0, \ldots, \sigma_{p - 1}(z, z')}.$ Similarly, define $\beta_{p - 1}(z, z')$ to be the number of intersections between $L_{z'}$ and the same section of $S^z,$ that is, $\beta_{p - 1}(z, z') = \sum\limits_{m = 0, \ldots, \sigma_{p - 1}(z, z')} \sum\limits_{z'' \in L_{z'}} \indic{\{S_m^z = z''\}}.$ We have the following limit $\lim\limits_{p \to \infty} \max\limits_{\norm{z - z'} \leq \delta_0} \Probability{\alpha_{p - 1}(z, z') = 0} = 0,$ where the maximum runs over all pairs $(z, z') \in \strC_{4 \cdot 9^{2p k_0}, 0} \times \strC_{4 \cdot 9^{2p k_0}, 0}'.$ This last limit is a corollary of proposition (\ref{Proposition: Limit probability that the LERW and a section of SRW will not intersect}): the translation invariance of the probability kernel shows at once that the maximum inside the limit does not depend on the pair $(z, z')$ above but solely on the distance separating them.
			\begin{Lemma}
				$\ds \lim_{p \to \infty} \max_{\norm{z - z'} \leq \delta_0} \Probability{\EuScript{G}_p(z') \cap \{\beta_{p - 1}(z, z') = 0\}} = 0,$ where the maximum runs over all pairs $(z, z') \in \strC_{4 \cdot 9^{2p k_0}, 0} \times \strC_{4 \cdot 9^{2p k_0}, 0}'.$
			\end{Lemma}
			\begin{proof}[Proof of lemma.]
				Observe that
				\begin{displaymath}\textstyle
					\Probability{\EuScript{G}_p(z') \cap \{\beta_{p - 1}(z, z') = 0\}} \leq \Probability{\EuScript{G}_p(z') \cap \{\alpha_{p - 1}(z, z') > 0, \beta_{p - 1}(z, z') = 0\}} + \Probability{\alpha_{p - 1}(z, z') = 0}.
				\end{displaymath}
				The second of these two terms tends uniformly to zero for the set of stated pairs $(z, z'),$ as was shown before the lemma. On the good event $\EuScript{G}_p$ we can apply the two decompositions above for $L_{z'}$ and $L^{z'}_{p - 1},$ giving
				\begin{displaymath}\textstyle
					 \EuScript{G}_p \cap \{\alpha_{p - 1}(z, z') > 0, \beta_{p - 1}(z, z') = 0\} \subset \left\{ S^z\text{ cross }S^{z'}\text{ inside }\strJ_{p - 1} \right\},
				\end{displaymath}
				where $\strJ_{p - 1} = \left\{ a_p - 2 \cdot 9^{(p - 1) k_0}, \ldots, a_p \right\} \times \Z^2.$ By (\ref{Proposition: Probability of intersection of two random walks in the strips}), we have that there exists a universal constant $c > 0$ and an index $p_3$ such that for $p \geq p_3,$ $\Probability{S^z\text{ cross }S^{z'}\text{ inside }\strJ_{p - 1}} \leq c 3^{-p k_0},$ this bound being uniform for $z' \in \strC_{4 \cdot 9^{2p k_0, 0}}'$ and $z \in \strC_{4 \cdot 9^{2p k_0}, 0}(z').$
			\end{proof}
	
			We continue now with the proof of the theorem. Recall we consider $z' \in \strC_{4 \cdot 9^{2p k_0}, 0}'$ and $z \in \strC_{4 \cdot 9^{2p k_0}, 0}(z').$ Introduce the following random time $\theta_{p - 1}(z, z') = \min \left( \sigma_{p - 1}(z, z'), \tau_{L_{z'}}(S^z) \right)$ which is the hitting time of $(\strH_{p - 1} + z') \cup L_{z'}$ by $S^z.$ Then, $\theta_{p - 1}(z, z')$ is a stopping time for the filtration
			\begin{displaymath}\textstyle
				\left( \scrX_k^{z'}(z) = \sigma \left( S^z_m\ \middle|\ 0 \leq m \leq k \right) \vee \scrF_\infty \left( \{0, z'\} \right) \right)_{k \in \Z_+},
			\end{displaymath}
			where $\scrF_\infty(\strT) = \sigma(S^t_m \mid t \in \strT, m \in \Z_+).$ Observe that $\EuScript{G}_p(z') \in \scrF_\infty(\{0, z'\}).$ We apply the strong Markov property to the process $\left( S^z_k, \scrX_k^{z'}(z) \right)_{k \in \Z_+},$
			\begin{equation}\textstyle\label{SubEquation: Conditional expectation of N and the good event}
				\Expectation{ \indic{ \EuScript{N}(z, z') \cap \EuScript{G}_p(z') } } = \Expectation{ \Expectation{ \indic{ \EuScript{N}(z, z') } \indic{ \EuScript{G}_p(z') }\ \middle|\ \mathscr{X}_{\theta_{p - 1}(z, z')}^{z'}(z)}} = \Expectation{ \indic{ \EuScript{G}_p(z') } \Probability{\EuScript{N}(z, z') \mid S_{\theta_{p - 1}(z, z')}^z} }.
			\end{equation}
			The definition of the stopping time $\theta_{p - 1}(z, z')$ shows that $\Probability{ \EuScript{N}(z, z') \middle| S_{\theta_{p - 1}(z, z')}^z  = z''} = 0$ for all $z'' \in L_{z'}$ lying to the left or on the plane $\strH_{p - 1} + z'.$ Therefore, we have
			\begin{displaymath}\textstyle
				\Probability{ \EuScript{N}(z, z') \middle| S_{\theta_{p - 1}(z, z')}^z} = \Probability{ \EuScript{N}(z, z') \middle| S_{\theta_{p - 1}(z, z')}^z} \indic{ \left\{ S_{\theta_{p - 1}(z, z')}^z \notin L_{z'} \right\} }.
			\end{displaymath}
			Clearly $\EuScript{N}(z, z') \subset \{\tau_0(S^z) < \infty\},$ also, $\theta_{p - 1}(z, z') = \sigma_{p - 1}(z, z'),$ on the event $\left\{ S_{\theta_{p - 1}(z, z')}^z \notin L_{z'}\right\};$ this implies $S_{\theta_{p - 1}(z, z')}^z \in \strH_{p - 1} + z'$ on this event. Whence,
			\begin{equation}\textstyle\label{SubEquation: Conditional probability of N given S sub theta}
				\Probability{ \EuScript{N}(z, z') \middle| S_{\theta_{p - 1}(z, z')}^z} \leq \left( \sup\limits_{z'' \in \strH_{p - 1} + z'} \Probability{ \tau_0(S^z) < \infty \mid S_0^z = z''} \right) \indic{ \left\{ S_{\theta_{p - 1}(z, z')}^z \notin L_{z'} \right\} }.
			\end{equation}
			We know from the bounds of the Green's function that
			\begin{equation}\textstyle\label{SubEquation: Upper bound in the probabiltiy that Sz will visit zero starting from Hp}
				\sup\limits_{z'' \in \strH_{p - 1} + z'} \Probability{\tau_0(S^z) < \infty \mid S^z_0 = z''} \leq c_1 |a_p|^{-1} \leq c 9^{-2p k_0},
			\end{equation}
			where $c > 0$ is a universal constant. 
	
			Use the bound of (\ref{SubEquation: Upper bound in the probabiltiy that Sz will visit zero starting from Hp}) in (\ref{SubEquation: Conditional probability of N given S sub theta}), and then apply this to (\ref{SubEquation: Conditional expectation of N and the good event}). Thus, we have shown up to this point in the argument the existence of a universal constant $c > 0$ such that 
			\begin{displaymath}\textstyle
				\Probability{\EuScript{N}(z, z') \cap \EuScript{G}_p(z')} \leq c 9^{-2p k_0} \Probability{ \EuScript{G}_p(z') \cap \left\{ S_{\theta_{p - 1}(z, z')}^z \notin L_{z'} \right\} }.
			\end{displaymath}
			It is also obvious that $\EuScript{G}_p(z') \cap \left\{ S_{\theta_{p - 1}(z, z')}^z \notin L_{z'} \right\} \subset \EuScript{G}_p(z') \cap \{\beta_{p - 1}(z, z') = 0\}.$ Hence,
			\begin{equation}\textstyle\label{SubEquation: Final bound in the probability of N and the good event}
				\Probability{\EuScript{N}(z, z') \cap \EuScript{G}_p(z')} \leq c 9^{-2p k_0} \Delta_p,
			\end{equation}
			with $\Delta_p = \max\limits_{\norm{z - z'} \leq \delta_0} \Probability{\EuScript{G}_p(z') \cap \{\beta_{p - 1}(z, z') = 0\}} \to 0$ as $p \to \infty,$ the pairs $(z, z') \in \strC_{4 \cdot 9^{2p k_0}, 0} \times \strC_{4 \cdot 9^{2p k_0}, 0}'.$ 
	
			We are in conditions to conclude the proof of the theorem. By the foregoing,
			\begin{align*}
				\Probability{\EuScript{C}(z) \cap \EuScript{H}(z')} &\leq \Probability{\EuScript{N}(z, z')}, &\text{use (\ref{SubEquation: The event C and H is contained in N})} \\
				&\leq \Probability{ \EuScript{N}(z, z') \cap \EuScript{G}_p(z') } + e^{-c' 9^{p k_0}} + e^{-\lambda 9^{-2p k_0}}, &\text{use (\ref{SubEquation: First bound in the probability of N})} \\
				&\leq c 9^{-2p k_0} \Delta_p + e^{-c' 9^{p k_0}} + e^{-\lambda 9^{-2p k_0}}, &\text{use (\ref{SubEquation: Final bound in the probability of N and the good event})}
			\end{align*}
			with the two constants $c, c' > 0$ being universal. We consider now $z' \in \strC_{4 \cdot 9^{2p k_0}, 0}'$ and sum over $z \in \strC_{4 \cdot 9^{2p k_0}, 0}(z'),$ this latter set has cardinality bounded above by $c'' \delta_0^2,$ with $c'' > 0$ another universal constant.  This shows, for every $z' \in \strC_{4 \cdot 9^{2p k_0}, 0}',$
			\begin{displaymath}\textstyle
				\sum\limits_{z \in \strC_{4 \cdot 9^{2p k_0}, 0}(z')} \Probability{\EuScript{C}(z) \cap \EuScript{H}(z')} \leq c'' \delta_0^2 \left( c 9^{-2p k_0} \Delta_p + e^{-c' 9^{p k_0}} + e^{-\lambda 9^{-2p k_0}} \right),
			\end{displaymath}
			this bound being uniform for $z'.$ We know also that $\card{ \strC_{4 \cdot 9^{2p k_0}, 0}' } \leq \card{ \strC_{4 \cdot 9^{2p k_0}, 0} } \leq c''' \cdot 9^{2p k_0},$ for yet another universal constant $c''' > 0.$ Whence, if we sum over $z' \in \strC_{4 \cdot 9^{2p k_0}, 0}',$ we reach a bound of the form
			\begin{displaymath}\textstyle
				\sum\limits_{z' \in \strC_{4 \cdot 9^{2p k_0}, 0}'} \sum\limits_{z \in \strC_{4 \cdot 9^{2p k_0}, 0}(z')} \Probability{\EuScript{C}(z) \cap \EuScript{H}(z')} \leq c \left( \Delta_p + 9^{2p k_0} \left( e^{-c' 9^{p k_0}} + e^{-\lambda 9^{2p k_0}} \right) \right),
			\end{displaymath}
			for a constant $c > 0,$ not depending on $p.$ The right hand side of this inequality tends to zero, and this completes the proof of Theorem (\ref{Theorem: One end in dimension two}) by (D).
		\end{proof}
	\numberwithin{SubCounter}{equation}
	\numberwithin{equation}{section}
	\setcounter{equation}{\value{auxCounter}}

	\subsection{One endedness for $\Gamma_d(\lambda),$ with $d \geq 3$}
	We finally deal with the case of ``high dimensions.'' Here, the proof is considerably simpler than $d = 2.$ We will make use of several general results of electrical networks, see \cite{LyPe:PTN} or \cite[Appendix C]{MaDib:USFwithDrift}.
	\begin{Theorem}\label{Theorem: One end in high dimensions}
		Let $d \geq 3.$ Then, for almost every realisation of the uniform spanning forest in $\Gamma_d(\lambda),$ all the components in this realisation have one end.
	\end{Theorem}
	\begin{proof}
		For every $p \in \N,$ consider the sets $\strA_p,$ $\strB_p,$ $\strC_p,$ defined in proposition (\ref{Proposition: Sum of the Green function on cylinder can be made small}).  Using the notation of this proposition, $p^{\frac{1}{4}}$ will be larger than $B$ starting at some index $p^*.$ As such, we may divide $\strC_p = \strC_{p, 0} \cup \strC_{p, 1},$ with $\strC_{p, 0} = \{-p\} \times \ball{0; p^{\frac{3}{4}}}$ and $\strC_{p, 1} = \strC_p \setminus \strC_{p, 0}.$
		\begin{EqInside}\label{EqInside: Sum of the Green function on Ep for big d}
			For every $p \geq p^*,$ $\sum\limits_{z \in \strA_p \cup \strB_p} G(z, 0) \leq c e^{-c' p}$ and $\sum\limits_{z \in \strC_{p, 1}} G(z, 0) \leq c e^{-c' p^{\frac{1}{2}}},$ where $c, c'$ are positive constants.			
		\end{EqInside}
	
		Assume $\mathfrak{F}$ is the random spanning forest of $\Gamma_d(\lambda)$ as constructed by Wilson's algorithm rooted at infinity. To be more precise, assume that on some probability space $(\Omega, \mathscr{F}, \mathds{P})$ there is a family of independent network random walks of $\Gamma_d(\lambda),$ say $(S^z)_{z \in \Z^{d + 1}}$ with $S^z$ started at $z,$ and then construct $\mathfrak{F}$ using some predefined order starting at zero of $\Z^{d + 1}$ following Wilson's algorithm. (For instance, start at zero and then search the $\norm{\cdot}_1$-spheres by first increasing the radii one unit at a time and then lexicographically ordering each sphere.)  We know then that $\mathfrak{F}$ is a random spanning forest following law $\USF$ of $\Gamma_d(\lambda);$ this signifies that $\Probability{\mathfrak{F} \in \scrE} = \USF(\scrE),$ for all events $\scrE$ in the probability space of $\USF$ spanning forests of $\Gamma_d(\lambda).$
	
		Because of the order chosen to construct $\mathfrak{F},$ $L_0 = \LE(S^0_m)_{m \in \Z_+}$ is (a.s.) an infinite path, we will call it ``first branch of the random forest $\mathfrak{F}.$'' Define the following set in $\Z^{d + 1}$ (consider $z \in \Z^{d + 1}$) and events in $\scrF:$
		\begin{displaymath}\textstyle
			\begin{split}
				\strE_p &= \strA_p \cup \strB_p \cup \strC_{p, 1}, \\
				\EuScript{E}(z) &= \{z\text{ is connected to }0\text{ in }\mathfrak{F}\setminus L_0\} \\
				\EuScript{E}_p &\textstyle= \bigcup\limits_{z \in \strE_p} \EuScript{E}(z) = \{\text{some vertex in } \strE_p \text{ is connected to }0\text{ in }\mathfrak{F} \setminus L_0\}.
			\end{split}
		\end{displaymath}
		Then, $\Probability{\EuScript{E}_p} \leq \sum\limits_{z \in \strE_p} \Probability{\EuScript{E}(z)}.$ To calculate $\Probability{\EuScript{E}(z)}$ we may apply corollary (\ref{Corollary: Main extension of WA}) and assume that $\mathfrak{F}$ was constructed using the ordering $(0, z, \ldots)$ of $\Z^{d + 1}.$ With this ordering, the event $\EuScript{E}(z)$ is the event where $S^z$ hits $L_0$ for the first time at vertex $0,$ whence $\Probability{\EuScript{E}_p} \leq \sum\limits_{z \in \strE_p} \Probability{\tau_0(S^z) < \infty} \leq \sum\limits_{z \in \strE_p} G(z, 0) < c e^{-c' p^{\frac{1}{2}}},$ where the last inequality follows from (\ref{EqInside: Sum of the Green function on Ep for big d}). By virtue of Borel-Cantelli lemma, we reach the existence of a random index $P:\Omega \to \Z_+$ (measurable relative to $\mathscr{F}$) such that for $\mathds{P}$-a.e. $\omega \in \Omega,$ $\omega \notin \bigcup\limits_{p \geq P(\omega)} \EuScript{E}_p.$
	
		Let $\mathfrak{U}$ denote the (random connected) component of $0$ in the forest $\mathfrak{F} \setminus L_0.$ Denote by $\mathfrak{U}_p$ the graph induced by $\mathfrak{U}$ in the cylinder $\strV_p = [-p, p] \times \ball{0; p}$ of $\Z \times \Z^d.$
	
			\begin{Lemma}\label{Lemma: In a transitive network, if the component of zero minus its the first branch is finite, then the forest has one end}
			If $\Probability{\mathfrak{U} \text{ is finite}} = 1,$ then $\Probability{\text{all the components of } \mathfrak{F}\text{ have one end}} = 1.$
		\end{Lemma}
		\begin{proof}[Proof of lemma.]
			For every $z \in \Z^{d + 1},$ denote by $\EuScript{T}(z)$ the event where the component of $z$ in $\mathfrak{F}$ is not one ended. Then,
			\begin{displaymath}\textstyle
				\EuScript{T} \mathop{=}\limits^{\text{def.}} \{\text{there exists a component of }\mathfrak{F}\text{ that is not one ended}\} = \bigcup\limits_{z \in \Z^{d + 1}} \EuScript{T}(z).
			\end{displaymath}
			We want to prove $\Probability{\EuScript{T}} = 0.$ By translation invariance, $\Probability{\EuScript{T}(z)}$ is independent of $z,$ and therefore, it suffices to prove that $\EuScript{T}(0)$ is a null event.
		
				Now, since $0$ is a vertex of $L_0 \subset \mathfrak{F}$ and this is an infinite branch (a.s.), the relation that the component of $0$ in $\mathfrak{F}$ is not one ended signifies that $\mathfrak{U}$ does not possess a ray, and this is equivalent to $\mathfrak{U}$ being finite since it is a locally finite tree. We established that $\EuScript{T}(0)^\complement = \{\mathfrak{U} \text{ is finite}\}.$ The conclusion of the lemma is now clear.
			\end{proof}

			We know that for $p \geq P,$ no edge in $\mathfrak{U}$ is adjacent to $\strE_p,$ and the only way $\mathfrak{U}$ can have infinitely many edges is by using edges adjacent to $\strC_{p, 0}.$ Let $\mathfrak{K}_{p - 1}$ be the (random) set of edges of $\strV_p = [-p, p] \times \ball{0; p}$ that are adjacent to $\mathfrak{U}_{p - 1}$ and $\strC_{p, 0}.$ We need to study the size of the sets $\mathfrak{K}_p,$ the following lemma will be needed.
	
		\begin{Lemma}
			Suppose that $\mathfrak{A}$ is a random spanning subgraph of $\strG$ defined on some probability space $(\Omega, \mathscr{F}, \mathds{P}).$ Suppose that $\mathfrak{A}$ satisfies the following: if $v$ and $v'$ are two vertices that are end-points of edges of $\mathfrak{A},$ then $v$ and $v'$ are in the same connected-component relative to $\mathfrak{A}.$ (Think of $\mathfrak{A}$ as a connected subgraph of $\strG$ plus all the vertices of $\strG.$) Denote by $(\strV_n)$ an exhaustion of $\strG$ such that $\closure{\strV}_n \subset \strV_{n + 1}$ for every $n \in \N.$ Let $\mathscr{F}_n$ denote the $\sigma$-algebra generated by the events $\{e \in \mathfrak{A}\}$ as $e$ runs through the edges of the graph induced on $\strV_n.$ Define $\mathfrak{A}_n$ to be the spanning subgraph of $\strG$ whose edge set consists of the edges of $\mathfrak{A}$ with both end-points belonging to $\strV_n$ (with $\mathfrak{A}_0 = \varnothing$). Consider the events ($n \in \N$) $\EuScript{G}_n = \{\mathfrak{A}_n \neq \mathfrak{A}_{n - 1}\},$ so that, $\EuScript{G}_n$ is the event where there was ``growth by	edges'' of $\mathfrak{A}$ in $\strV_n.$ Set $Y_n = \Probability{\EuScript{G}_n \mid \mathscr{F}_{n - 1}}$ ($n \in \N$). Then $(\EuScript{G}_n)_{n \in \N}$ is a decreasing sequence, if $\EuScript{G}_\infty = \bigcap\limits_{n \in \N} \EuScript{G}_n,$ then $\EuScript{G}_\infty$ is the event where $\mathfrak{A}$ has infinite edges; furthermore, on $\EuScript{G}_\infty,$ $Y_n \to 1.$
		\end{Lemma}
		\begin{proof}[Proof of lemma.]
			It is clear that $\EuScript{G}_n^\complement \subset \EuScript{G}_{n + 1}^\complement,$ hence $(\EuScript{G}_n)$ is decreasing. Next, $1 \geq \Probability{\bigcup\limits_{n \in \N} \EuScript{G}_n^\complement \cap \EuScript{G}_{n - 1}} = \sum\limits_{n \in \N} \Expectation{ \Expectation{ \indic{\EuScript{G}_n^\complement}\indic{\EuScript{G}_{n - 1}} \middle| \mathscr{F}_{n - 1} } }.$ It is clear that $\EuScript{G}_{n - 1} \in \mathscr{F}_{n - 1}$ and so, $1 \geq \Expectation{\sum\limits_{n \in \N} \indic{\EuScript{G}_{n - 1}} (1 - Y_n)}.$ Therefore, the positive random variable $\sum\limits_{n \in \N} \indic{\EuScript{G}_{n - 1}} (1 - Y_n)$ is finite for $\mathds{P}$-a.e. realisation $\omega \in \Omega.$ Thus, $\sum\limits_{n \in \N} (1 - Y_n)$ is finite $\mathds{P}$-a.s. on the event $\EuScript{G}_\infty,$ which immediately implies the conclusion to be reached, namely, that $Y_n \to 1$ on $\EuScript{G}_\infty.$
		\end{proof}

		Let us return to the proof of Theorem (\ref{Theorem: One end in high dimensions}). We apply the previous lemma to the random graph $\mathfrak{U},$ and the exhaustion $\strV_p = [-p, p] \times \ball{0; p}.$ Then, the events constructed in the lemma are $\EuScript{G}_p = \{\mathfrak{U}_p \neq \mathfrak{U}_{p - 1}\}$ and the $\sigma$-algebras are $\mathscr{F}_p = \sigma(\{e \in \mathfrak{U}_p\}; e \text{ is an edge of the graph induced on } \strV_p  \text{ by } \Z^{d + 1}).$ By construction, $\mathfrak{K}_{p - 1}$ is $\mathscr{F}_{p - 1}$-measurable. Set $Y_p = \Probability{\EuScript{G}_p \mid \mathscr{F}_{p - 1}},$ so that $\mathds{P}$-a.s. $Y_p \to 1$ on the event where $\mathfrak{U}$ is infinite. Now,
		\begin{align*}
			1 - Y_p &= \Probability{\mathfrak{U} = \mathfrak{U}_{p - 1} \mid \mathscr{F}_{p - 1}} = \Probability{\mathfrak{U} = \mathfrak{U}_{p - 1}, p \geq P \mid \mathscr{F}_{p - 1}} + \Probability{\mathfrak{U} = \mathfrak{U}_{p - 1}, p < P \mid \mathscr{F}_{p - 1}} \\
			&= \Probability{\mathfrak{F} \cap \mathfrak{K}_{p - 1} = \varnothing, p \geq P \mid \mathscr{F}_{p - 1}} + \Probability{\mathfrak{U} = \mathfrak{U}_{p - 1}, p < P \mid \mathscr{F}_{p - 1}} \\
			&= \Probability{\mathfrak{F} \cap \mathfrak{K}_{p - 1} = \varnothing \mid \mathscr{F}_{p - 1}} - \Probability{\mathfrak{F} \cap \mathfrak{K}_{p - 1} = \varnothing, p < P \mid \mathscr{F}_{p - 1}} \\
			&+ \Probability{\mathfrak{U} = \mathfrak{U}_{p - 1}, p < P \mid \mathscr{F}_{p - 1}}.
		\end{align*}
		Observe that $\{\mathfrak{U} = \mathfrak{U}_{p - 1}\} \subset \{\mathfrak{F} \cap \mathfrak{K}_{p - 1} = \varnothing\},$ since, by definition, $\mathfrak{K}_{p - 1}$ are the edges adjacent to both $\strC_{p, 0}$ and $\mathfrak{U}_{p - 1}$ in $\strV_p.$ Thus, $1 - Y_p = \Probability{\mathfrak{F} \cap \mathfrak{K}_{p - 1} = \varnothing \mid \mathscr{F}_{p - 1}} - \Probability{\mathfrak{F} \cap \mathfrak{K}_{p - 1} = \varnothing, \mathfrak{U} \neq \mathfrak{U}_{p - 1}, p < P \mid \mathscr{F}_{p - 1}}.$ If $\mathfrak{U}$ grew at stage $p,$ but it did not so through $\mathfrak{K}_{p - 1},$ then $\mathfrak{U}$ grew through some edge adjacent to $\strE_p,$ hence
		\begin{displaymath}\textstyle
			\Probability{\mathfrak{F} \cap \mathfrak{K}_{p - 1} = \varnothing, \mathfrak{U} \neq \mathfrak{U}_{p - 1}, p < P \mid \mathscr{F}_{p - 1}} \leq \Probability{\EuScript{E}_p \mid \mathscr{F}_{p - 1}}.
		\end{displaymath}
		We know that $\Probability{\EuScript{E}_p} \to 0,$ so the sequence $(\Probability{\EuScript{E}_p \mid \mathscr{F}_{p - 1}})$ converges to zero in $\mathscr{L}^1.$ We may assume, passing through a subsequence should the need arise, that $(\Probability{\EuScript{E}_p \mid \mathscr{F}_{p - 1}})$ converges to zero $\mathds{P}$-a.s. Therefore, $1 - Y_p \geq \Probability{\mathfrak{F} \cap \mathfrak{K}_{p - 1} = \varnothing \mid \mathscr{F}_{p - 1}} - \Probability{\EuScript{E}_p \mid \mathscr{F}_{p - 1}}$ thus, on the event where $\mathfrak{U}$ is an infinite component, $\Probability{\mathfrak{F} \cap \mathfrak{K}_{p - 1} = \varnothing \mid \mathscr{F}_{p - 1}} \to 0.$
	
		Denote by $\mathfrak{G}_p$ the random graph obtained from $\Z^{d + 1}$ where all edges $e \in \mathfrak{U}_p$ are shorted and the edges $e$ in $\strV_p$ that are not in $\mathfrak{U}_p$ are cut ($p \in \N$). Elementary properties of wired spanning forest measures show that $\Probability{\mathfrak{F} \cap \mathfrak{K}_{p - 1} = \varnothing \mid \mathscr{F}_{p - 1}} = \WSF_{\mathfrak{G}_{p - 1}}\left(\mathfrak{F} \cap \mathfrak{K}_{p - 1} = \varnothing \right),$ the right hand side meaning the $\WSF$ measure on the random graph $\mathfrak{G}_{p - 1};$ in this graph however, $\mathfrak{K}_{p - 1}$ is not random, while $\mathfrak{F}$ is its $\WSF$-distributed random object. We can apply \cite[Lemma 10.40]{LyPe:PTN} and conclude
		\begin{displaymath}\textstyle
			\WSF_{\mathfrak{G}_{p - 1}}\left(\mathfrak{F} \cap \mathfrak{K}_{p - 1} = \varnothing \right) \geq \prod\limits_{e \in \mathfrak{K}_{p - 1}} \frac{1}{1 + \mu(e) \effRes[_{\mathfrak{G}_{p - 1} \setminus \mathfrak{K}_{p - 1}}^W]{e^-, e^+}},
		\end{displaymath}
		where $\effRes[_\strG^W]{e^-, e^+}$ stands for the wired effective resistance between $e^-$ and $e^+$ in the network $\strG.$ We know that $\mu(e) \asymp e^{-\lambda p}$ for $e \in \mathfrak{K}_{p - 1}$ and $\effRes[_{\mathfrak{G}_{p - 1} \setminus \mathfrak{K}_{p - 1}}^W]{e^-, e^+} \leq \effRes[_{\mathfrak{G}_{p - 1} \setminus \mathfrak{K}_{p - 1}}]{e^-, \infty} + \effRes[_{\mathfrak{G}_{p - 1} \setminus \mathfrak{K}_{p - 1}}]{e^+, \infty}.$ The following lemma is crux in the proof and shows at once why $d \geq 3$ is needed (which has not been needed so far).
		\begin{Lemma}
			Let $d \geq 3.$ For $e \in \mathfrak{K}_{p - 1},$ we have $\effRes[_{\mathfrak{G}_{p - 1} \setminus \mathfrak{K}_{p - 1}}]{e^\pm, \infty} \asymp e^{\lambda p}$ (any implicit constant being universal).
		\end{Lemma}
		\begin{proof}[Proof of lemma.]
			By Thompson's principle $\effRes[_\Gamma]{v, \infty} = \inf\limits_{\theta} \energyFlow{\theta},$ where $\theta$ runs on the sets of unit flows from $v$ to $\infty.$ Since $d \geq 3,$ we have that the standard random walk of $\Z^d$ is transient and T. Lyon's criterion \cite{LyonsTe:Criterion} shows the existence of a unit flow $\theta_0$ from $0$ to infinity having finite energy in the graph $\Z^d.$ We know that the network $\mathfrak{G}_{p - 1} \setminus \mathfrak{K}_{p - 1}$ contains the hyperplane $\strH_p = \{-p\} \times \Z^d.$ Observe that $\strH_p$ as a subnetwork of $\mathfrak{G}_{p - 1} \setminus \mathfrak{K}_{p - 1}$ has constant resistances equal to $e^{\lambda p}.$ Translate $\theta_0$ to $\strH_p$ in the obvious way and denote the new flow by $\theta_p.$ Clearly, $\energyFlow[_{\strH_p}]{\theta_p} = e^{\lambda p} \energyFlow[_{\Z^d}]{\theta_0},$ where $\energyFlow[_\strG]{\theta}$ is the energy of the flow $\theta$ on the network $\strG.$ Extend $\theta_p$ to $\mathfrak{G}_{p - 1} \setminus \mathfrak{K}_{p - 1}$ by sending zero flow through the edges outside $\strH_p.$ The energy of the extension of $\theta_p,$ denoted by $\closure{\theta}_p,$ does not change, that is to say, $\energyFlow[_{\strH_p}]{\theta_p} = \energyFlow[_{\mathfrak{G}_{p - 1} \setminus \mathfrak{K}_{p - 1}}]{\closure{\theta}_p}.$ Thus, we proved that there exists a $c = \energyFlow[_{\Z^d}]{\theta_0} > 0$ such that for all $p \in \N$ and all edges $e \in \mathfrak{K}_{p - 1},$ $\effRes[_{\mathfrak{G}_{p - 1} \setminus \mathfrak{K}_{p - 1}}]{e^\pm, \infty} \leq c e^{\lambda p}.$
		
				The reverse inequality is not needed but we prove it for completeness. Here we use that $\effCon{v, \infty} = \effRes{v, \infty}^{-1}$ and we bound from above the effective conductance. We apply Dirichlet's principle. Consider a vertex $z = (-p, x)$ and define $\varphi = \indic{z}.$ Then, $\effCon[_{\mathfrak{G}_{p - 1} \setminus \mathfrak{K}_{p - 1}}]{z, \infty} \leq \energyDir[_{\mathfrak{G}_{p - 1} \setminus \mathfrak{K}_{p - 1}}]{\varphi} \leq \sum\limits_{\substack{e \text{ and edge of } \Gamma_d(\lambda) \\ e^- = z}} \mu(e) \leq 2(d + 1) e^\lambda e^{-\lambda p}.$ Thus, for $c^{-1} = 2(d + 1)e^\lambda > 0,$ we have $\effRes[_{\mathfrak{G}_{p - 1} \setminus \mathfrak{K}_{p - 1}}]{z, \infty} \geq c e^{\lambda p}.$
		\end{proof}

		We may continue with the proof of one endedness in $\Gamma_d(\lambda).$ By virtue of the previous lemma, there exists a universal $\delta > 0$ such that $\mathds{P}$-a.s. $\WSF_{\mathfrak{G}_{p - 1}}\left(\mathfrak{F} \cap \mathfrak{K}_{p - 1} = \varnothing \right) \geq \delta^{\card{\mathfrak{K}_{p - 1}}}.$ On the event $\mathfrak{U}$ is infinite, we proved that $\WSF_{\mathfrak{G}_{p - 1}}\left(\mathfrak{F} \cap \mathfrak{K}_{p - 1} = \varnothing \right) \to 0$ and so, on this event, $\card{\mathfrak{K}_{p - 1}} \to \infty.$
	
		Consider then the random variables $Z_{p - 1} = \card{\mathfrak{K}_{p - 1}}.$ Notice that $(-p, x) \in \strC_{p, 0}$ is adjacent to an edge of $\mathfrak{K}_{p - 1}$ if and only if $(-p + 1, x)$ is adjacent to an edge in $\mathfrak{U}_{p - 1}.$ Thus,
		\begin{displaymath}\textstyle
			\Expectation{Z_{p - 1}} = \sum\limits_{z \in \strC_{p - 1, 0}} \Probability{z \text{ is connected to } 0 \text{ in } \mathfrak{U}_{p - 1}} \leq \sum\limits_{z \in \strC_{p - 1, 0}} \Probability{z \text{ is connected to } 0 \text{ in } \mathfrak{F} \setminus L_0}.
		\end{displaymath}
		Using (\ref{Corollary: Main extension of WA}), we may calculate $\Probability{z \text{ is connected to } 0 \text{ in } \mathfrak{F} \setminus L_0}$ assuming that the order in which the vertices of $\Z^{d + 1}$ were searched was $(0, z, \ldots).$ Then, the event where $z$ is connected to $0$ in $\mathfrak{F} \setminus L_0$ is the event where $S^z$ visits $0$ before hitting $L_0$ anywhere else. In particular, $\Expectation{Z_{p - 1}} \leq \sum\limits_{z \in \strC_{p - 1, 0}} \Probability{\tau_0(S^z) < \infty} \leq \sum\limits_{z \in \strC_{p - 1, 0}} G(z, 0) \leq c$ by (\ref{EqInside: Sum of the Green function on Ep for big d}). This implies $\Probability{Z_p \to \infty} = 0$ (by Fatou's lemma). Finally, we showed that $\Probability{\card{\mathfrak{U}} = \infty} \leq \Probability{Y_p \to 1} \leq \Probability{Z_p \to \infty} = 0.$ The proof of (\ref{Theorem: One end in high dimensions}) is complete by Lemma (\ref{Lemma: In a transitive network, if the component of zero minus its the first branch is finite, then the forest has one end}).
	\end{proof}

	\subsection{On the rays of $\Gamma_d(\lambda)$}
	We will say that a path $\gamma = (v_j)_{j \in \Z_+}$ is \textbf{ultimately inside a region $\strR$} if there exists an almost surely finite random index $J$ such that $v_j \in \strR$ for all $j \geq J.$

	Denote $\strP = \left\{ z = (n, x) \in \Z \times \Z^d;\ \norm{x}^2 \leq n \right\}$ and $\strP^\varepsilon = \left\{ z;\ \norm{x}^2 \leq n^{1 + \varepsilon} \right\}.$
	\begin{Theorem}
		Let $z$ be any vertex of $\Gamma_d(\lambda).$ Then, there exists one and only one ray in $\USF$ ($\UST$ in $d = 1, 2$) starting at $z.$ This ray is a random object (is $\USF$-measurable). Furthermore, for every $\varepsilon > 0,$ this ray is ultimately inside $z + \strP^\varepsilon,$ for every $\varepsilon > 0.$
	\end{Theorem}
	\begin{proof}
		By Theorems (\ref{Theorem: One end in dimension one}), (\ref{Theorem: One end in dimension two}) and (\ref{Theorem: One end in high dimensions}) there exists one and only one ray, this proves at once such ray is $\USF$-measurable. There remains to prove that said ray is ultimately inside $z + \strP^\varepsilon$ for every $\varepsilon > 0.$ We may construct $\UST$ in $d = 1, 2$ or $\USF$ for $d \geq 3$ using Wilson's algorithm rooted at infinity with the ordering of $\Z^{d + 1}$ to be $(z, \ldots).$ Thus, suffices to show that if $S^z$ is the network random walk of $\Gamma_d(\lambda)$ started at $z,$ then $S^z$ is ultimately inside $z + \strP^\varepsilon;$ furthermore, by translation invariance, we may assume $z = 0.$
	
		Let $S$ denote the network random walk of $\Gamma_d(\lambda)$ started at zero. To show that $S$ is ultimately inside $\strP$ is the same as showing that, for every $\varepsilon > 0,$ $S$ will only be finitely many times outside $\strP^\varepsilon.$ We may apply Borel-Cantelli lemma to show this and as such, it all reduces to show that $\sum\limits_{n \in \Z_+} \Probability{S_n \notin \strP^\varepsilon} < \infty.$ The Green's function estimates (\ref{Theorem: Greens function bounds}) give at once $\sum\limits_{n \in \Z_+} \Probability{S_n \notin \strP^\varepsilon} \leq \sum\limits_{z \notin \strP^\varepsilon} G(0, z) = \sum\limits_{z \in \strS_1} G(0, z) + \sum\limits_{z \in \strS_2 \cap \strP^\varepsilon} G(0, z),$ where $\strS_1 = \{z; n \leq \norm{x}\}$ and $\strS_2 = \strS_1^\complement.$ Since $\sum\limits_{z \in \strS_1} G(0, z) \leq \sum\limits_{z \in \Z^{d + 1}} e^{-c \norm{z}} < \infty,$ we focus on the second term. We have (write $z = (n, x)$) $\sum\limits_{z \in \strS_2 \cap \strP^\varepsilon} G(0, z) \leq c_1 \sum\limits_{n \in \Z_+} \sum\limits_{n^{1 + \varepsilon} \leq \norm{x}^2 \leq n^2} e^{-c_2 \frac{\norm{x}^2}{n}} \leq c \sum\limits_{n \in \Z_+} e^{-c_2 n^\varepsilon} n^d < \infty.$
	\end{proof}

	\section[Probability of same component]{Probability of belonging in the same component}\label{Section: Probability of same component}
In this section we will provide a function $u:\Z^{d + 1} \to \R_+$ with ``polynomial decay'' such that for all $z \in \Z^{d + 1},$ the probability that $0$ and $z$ are in the same tree of $\Gamma_d(\lambda)$ is $\asymp u(z);$ the important part is that the decay depends solely on dimension. Recall that for $d = 1$ or $d = 2,$ on $\Gamma_d(\lambda)$ we have a tree, thus the probability that $0$ and any vertex $z$ are connected is $1,$ thus no decay.

	\subsection[Probability that two vertices are $\USF$-connected]{The probability that two vertices are connected in the $\USF$}
	As usual, we will write $z = (n, x), z' = (n', x'), z_1 = (n_1, x_1),$ etc. to denote the vertices of $\Z^{d + 1} = \Z \times \Z^d.$ Define
	\begin{equation}\textstyle\label{Equation: Metric in the probability of same component}
		\eta(z) = \max \left( |n|^\frac{1}{2}, \norm{x} \right), \quad z = (n, x) \in \Z^{d + 1}.
	\end{equation}
	This $\eta$ has the following properties, all of which are obvious.
	\begin{enumerate}
		\item $\eta(z) \geq 0,$ and $\eta(z) = 0$ is equivalent to $z = 0.$
		\item $\eta(z) = \eta(-z);$
		\item $\eta(z + z') \leq \eta(z) + \eta(z').$
	\end{enumerate}
	Thus, the function $(z, z') \mapsto \eta(z - z')$ is a metric on $\Z^{d + 1}.$
	\begin{Theorem}\label{Theorem: Probability of two points belonging in the same component}
		Let $d \geq 3.$ Then, $\USF(0\text{ is connected to }z) \asymp \eta(z)^{-(d - 2)},$ with any constant depending solely on dimension.
	\end{Theorem}
		\setcounter{auxCounter}{\value{equation}}
		\setcounter{auxSubCounter}{\value{SubCounter}}
		\setcounter{equation}{0}
		\setcounter{SubCounter}{0}
		\renewcommand{\theequation}{\Alph{equation}}
		\renewcommand{\theSubCounter}{\thesection.\theauxCounter.\arabic{SubCounter}}
	\begin{proof}
		We remark that $\eta(z)^{-(d - 2)} \asymp \max \left( n^{-\frac{d}{2} + 1}, \norm{x}^{-(d - 2)} \right),$ and thus it suffices to show the bounds for any of these expressions. We also remark that for every $c > 0,$ there exists $c' > 0$ such that $e^{-c \norm{z}} \leq c' \eta(z)^{-(d - 2)};$ thus, if at some point, we show an upper exponential bound in an estimate, we are done with that particular estimate.
	
		Consider the event $\scrC(0, z)$ that $0$ is connected to $z$ in $\USF.$ Since the measure $\USF$ is translation invariant in $\Gamma_d(\lambda),$ $\USF(\scrC(0, z)) = \USF(\scrC(z', z + z'))$ for any $z' \in \Z^{d + 1}.$ It is also clear that $\scrC(0, z) = \scrC(z, 0).$ Thus, we may assume $z \in \Z_+ \times \Z^d.$ Construct now $\USF$ using Wilson's algorithm rooted at infinity using the order $(0, z, \ldots)$ of $\Z^{d + 1},$ denote by $\mathfrak{F}$ the random forest constructed in this way. Thus, we have on some probability space $(\Omega, \scrF, \mathds{P})$ a family of independent network random walks $(S^v)_{v \in \Z^{d + 1}},$ with $S^v$ started at $v.$ We know $\mathfrak{F} \sim \USF,$ and by construction $\left\{ \mathfrak{F} \in \scrC(0, z) \right\} = \left\{ S^z \cap \LE(S^0_m)_{m \in \Z_+} \neq \varnothing \right\}.$ Thus, and with the aid of (\ref{Proposition: Lower bound in the probability of intersection of LERW and RW}),
		\begin{equation}\textstyle\label{SubEquation: Bounds above and below of the probability of same component}
			\USF(\scrC(0, z)) = \Probability{S^z \cap \LE(S^0_m)_{m \in \Z_+} \neq \varnothing} \asymp \Probability{S^z \cap S^0 \neq \varnothing}.
		\end{equation}
		Denote by $K_z$ the number of intersections between $S^z$ and $S^0.$ In other words,
		\begin{equation}\textstyle\label{SubEquation: Definition of Kz}
			K_z = \sum\limits_{z' \in \Z^{d + 1}} \sum\limits_{(p, q) \in \Z_+^2} \indic{\left\{S^z_p = z', S^0_q = z' \right\}}.
		\end{equation}
		Clearly, $\Probability{K_z > 0} = \Probability{S^z \cap S^0 \neq \varnothing}.$ Independence and the Lebesgue-Tonelli theorem show at once
		\begin{equation}\textstyle\label{SubEquation: Expectation of Kz as sums of Greens function}
			\Expectation{K_z} = \sum\limits_{z' \in \Z^{d + 1}} G(z, z') G(0, z').
		\end{equation}
		\begin{Lemma}\label{Lemma: The probability of Kz is positive is asymptotic to its expectation}
			We have $\Probability{K_z > 0} \asymp \Expectation{K_z},$ where any implicit constant is universal (does not depend on $z$).
		\end{Lemma}
		\begin{proof}[Proof of Lemma]
			Since $K_z$ is a $\Z_+$-valued random variable, we have $\Probability{K_z > 0} \leq \Expectation{K_z}.$ There remains to prove a lower bound of the form $\Probability{K_z > 0} \geq c \Expectation{K_z},$ for a universal constant $c > 0.$ Using the second moment inequality (\ref{Equation: Second moment lower bound}), suffices to establish that $\Expectation{K_z^2} \leq c \Expectation{K_z},$ for some universal constant $c > 0.$ By definition,
			\begin{displaymath}\textstyle
				K_z^2 = \sum\limits_{(z', z'') \in \Z^{d + 1} \times \Z^{d + 1}} \sum\limits_{(p, q, a, b) \in \Z_+^4} \indic{\left\{ S^z_p = z', S^z_a = z'' \right\}} \indic{\left\{S^0_q = z', S^0_b = z'' \right\}}.
			\end{displaymath}
			By adding terms corresponding to $p = a$ or $q = b,$ then expanding the sum as to whether $p \leq a$ or $p \geq a$ and $q \leq b$ or $q \geq b,$ we reach a bound with four sums
			\begin{displaymath}\textstyle
				K_z^2 \leq L_1 + L_2 + L_3 + L_4,
			\end{displaymath}
			and
			\begin{align*}
				\Expectation{L_1} &\textstyle= \sum\limits_{(z', z'') \in \Z^{d + 1} \times \Z^{d + 1}} G(z, z') G(z', z'')^2 G(0, z') \\
				\Expectation{L_2} &\textstyle= \sum\limits_{(z', z'') \in \Z^{d + 1} \times \Z^{d + 1}} G(z, z') G(z', z'') G(z'', z') G(0, z'') \\
				\Expectation{L_3} &\textstyle= \sum\limits_{(z', z'') \in \Z^{d + 1} \times \Z^{d + 1}} G(z, z'') G(z'', z') G(z', z'') G(0, z') \\
				\Expectation{L_4} &\textstyle= \sum\limits_{(z', z'') \in \Z^{d + 1} \times \Z^{d + 1}} G(z, z'') G(z'', z')^2 G(0, z'')
			\end{align*}
			It is clear that $\Expectation{L_1} = \Expectation{L_4}$ and $\Expectation{L_2} = \Expectation{L_3},$ thus
			\begin{displaymath}\textstyle
				\Expectation{K_z^2} \leq 2 \Big[ \Expectation{L_1} + \Expectation{L_2} \Big]
			\end{displaymath}
			We first handle the expectation of $L_1,$ we obtain
			\begin{displaymath}\textstyle
				\Expectation{L_1} = \sum\limits_{z' \in \Z^{d + 1}} G(z, z') G(0, z') \sum\limits_{z'' \in \Z^{d + 1}} G(z', z'')^2 = \left[ \sum\limits_{z'' \in \Z^{d + 1}} G(0, z'')^2 \right] \sum\limits_{z' \in \Z^{d + 1}} G(z, z') G(0, z'),
			\end{displaymath}
			and $\sum\limits_{z'' \in \Z^{d + 1}} G(0, z'')^2$ is a universal constant by the bubble condition (\ref{Theorem: The bubble condition}).

			We now handle $L_2.$ First, the translation invariance of the Green's function shows at once that 
			\begin{displaymath}\textstyle
				\Expectation{L_2} = \sum\limits_{(z', z'')} G(z, z') G(z', z'') G(z'', z') G(0, z'') = \sum\limits_{z'} G(z, z') \sum\limits_{z''}  G(0, z'') G(z'', 0) G(0, z' + z'').
			\end{displaymath}
			In the following calculations we are going to establish that the inner sum in the line above, when viewed as a function of $z',$ possesses the same type of upper bounds as the Green's function's upper bounds (\ref{Theorem: Greens function bounds}). There are two cases to consider. Write $z' = (n', x')$ and similarly $z'' = (n'', x'').$ 
			\begin{description}
				\item[Case 1.] Here we assume $\norm{x'} \leq n'.$ We then divide the inner sum into two, the first consists of all $z''$ satisfying $|n''| \leq \frac{n'}{2}$ and $\norm{x''} \leq \frac{\norm{x'}}{2}$ and the second sum is over all other $z''.$ In the first sum, we reach $\norm{x' + x''} \asymp \norm{x'}$ and $|n' + n''| \asymp n' \asymp \norm{z'},$ then $G(0, z' + z'') \leq c e^{-c'\frac{\norm{x'}^2}{n'}} \norm{z'}^{-\frac{d}{2}}$ for some constants $c, c' > 0.$ This entails
				\begin{displaymath}\textstyle
					\sum\limits_{|n''| \leq \frac{n'}{2}, \norm{x''} \leq \frac{\norm{x'}}{2}} G(z'', 0)G(0, z'') G(0, z' + z'') \leq c e^{-c'\frac{\norm{x'}^2}{n'}} \norm{z'}^{-\frac{d}{2}} \sum\limits_{z'' \in \Z^{d + 1}} G(z'', 0)G(0, z''),
				\end{displaymath}
				for a pair of constants $c, c' > 0.$ We may apply the inequality $ab \leq a^2 + b^2,$ and obtain that $\sum\limits_{z'' \in \Z^{d + 1}} G(z'', 0)G(0, z'') < \infty$ by the bubble condition (\ref{Theorem: The bubble condition}). This shows that
				\begin{displaymath}\textstyle
					\sum\limits_{|n''| \leq \frac{n'}{2}, \norm{x''} \leq \frac{\norm{x'}}{2}} G(z'', 0)G(0, z'') G(0, z' + z'') \leq c e^{-c'\frac{\norm{x'}^2}{n'}} \norm{z'}^{-\frac{d}{2}},
				\end{displaymath}
				for a pair of positive constants $c, c' > 0.$ Next, we are going to handle the second sum, which is over all $z''$ for which either $|n''| > \frac{n'}{2}$ or $\norm{x''} > \frac{\norm{x'}}{2}.$ Here, we will use that one of the two factors $G(z'', 0),$ $G(0, z'')$ has an upper bound of the form $e^{-c\norm{z''}}.$ We will write $e^{-c \norm{z''}} = e^{-\frac{c}{2} \norm{z''}} e^{-\frac{c}{2} \norm{z''}}$ and use the fact that the Green's function is bounded by some constant $L > 0.$ Since $\norm{z'} \asymp n',$ we may write
				\begin{displaymath}\textstyle
					\sum\limits_{|n''| > \frac{n'}{2}, x'' \in \Z^d} G(z'', 0) G(0, z'') G(0, z' + z'') \leq L^2 e^{-\frac{c}{4} n'} \sum\limits_{z'' \in \Z^{d + 1}} e^{-\frac{c}{2} \norm{z''}} \leq c' e^{-c'' \norm{z'}},
				\end{displaymath}
				for positive constants $c, c', c'' > 0.$ It is clear that $c' e^{-c'' \norm{z'}} \leq C e^{-C' \frac{\norm{x'}^2}{n'}} \norm{z'}^{-\frac{d}{2}}$ for another pair of constants $C, C' > 0.$ We finally deal with the sum corresponding to $|n''| \leq \frac{n'}{2}$ and $\norm{x''} > \frac{\norm{x'}}{2}.$ Here we decompose again $e^{-c \norm{z''}} = \left( e^{-\frac{c}{2} \norm{z''}} \right)^2.$ The assumption $\norm{x'} \leq n'$ shows that $\norm{x'} \geq \frac{\norm{x'}^2}{n'},$ and thus
				\begin{displaymath}\textstyle
					e^{-\frac{c}{2} \norm{z''}} \leq e^{-\frac{c}{2} \norm{x''}} \leq e^{-\frac{c}{4} \norm{x'}} \leq e^{-\frac{c}{6} \frac{\norm{x'}^2}{n'}}.
				\end{displaymath}
				Also, it is clear that $G(0, z' + z'') \leq c' \norm{z'}^{-\frac{d}{2}}$ for some positive constant $c' > 0$ since $\norm{z' + z''} \asymp \norm{z'}.$ Thus,
				\begin{displaymath}\textstyle
					\sum\limits_{|n''| \leq \frac{n'}{2}, \norm{x''} > \frac{\norm{x'}}{2}} G(z'', 0) G(0, z'') G(0, z' + z'') \leq c' e^{-\frac{c}{6} \frac{\norm{x'}^2}{n'}} \norm{z'}^{-\frac{d}{2}},
				\end{displaymath}
				for a pair of constants.

				\item[Case 2.] Here we assume $|n'| < \norm{x'}.$ Split $\sum\limits_{z''} G(z'', 0) G(0, z'') G(0, z' + z'') = A + B,$ where $A$ is the sum corresponding to all $z'' = (n'', x'')$ satisfying $|n''| \leq \frac{\norm{x'}}{2}$ and $\norm{x''} \leq \frac{\norm{x'}}{2},$ $B$ is the sum over all other $z''.$ For $z''$ in the range of $A,$ we obtain $\norm{x' + x''} \asymp \norm{x'} \asymp \norm{z'}$ and $|n' + n''| \leq c \norm{x'} \leq c \norm{z'},$ and this shows that, on this range, $G(0, z' + z'') \leq c e^{-c' \norm{z'}}$ for a pair of constants $c, c' >0.$ Then, there exists two universal constants $c, c' > 0$ such that (recall from \textbf{Case 1} that $G(z'', 0) G(0, z'')$ is a summable family for $z'' \in \Z^{d + 1}$) $A \leq c e^{-c' \norm{z'}}.$ In the range of $B,$ we have that either $|n''| > \frac{\norm{x'}}{2} \asymp \norm{z'}$ or $\norm{x''} > \frac{\norm{x'}}{2} \asymp \norm{z'},$ thus $\norm{z''} \geq c \norm{z'}$ for a universal constant $c > 0.$ Next, one of the two factor $G(z'', 0)$ or $G(0, z'')$ is bounded from above by $e^{-c \norm{z''}}$ by (\ref{Theorem: Greens function bounds}). Let $L > 0$ be a bound for the Green's function of $\Gamma_d(\lambda).$ Then $B \leq L^2 e^{-c' \norm{z'}} \sum\limits_{z'' \in \Z^{d + 1}} e^{-\frac{c}{2} \norm{z''}} \leq C e^{-c' \norm{z'}}.$ 
			\end{description}
			By virtue of the two cases above, we have established that there exists two constants $c, c' > 0$ such that (write $z' = (n', x')$)
			\begin{displaymath}\textstyle
				\sum\limits_{z'' \in \Z^{d + 1}} G(z'', 0) G(0, z'') G(z' + z'') \leq
				\begin{cases}
					c e^{-c' \norm{z'}} &\text{ if } \norm{x'} > n', \\
					c e^{-c' \frac{\norm{x'}^2}{n'}} \norm{z'}^{-\frac{d}{2}} &\text{ if } \norm{x'} \leq n'.
				\end{cases}
			\end{displaymath}
			The proof of (\ref{Lemma: The probability of Kz is positive is asymptotic to its expectation}) is now an immediate consequence of the next lemma. Indeed, the next lemma shows at once that $\Expectation{K_z} \asymp \eta(z)^{-(d - 2)}$ and that $\Expectation{L_2} \leq c \eta(z)^{-(d - 2)}.$ Therefore, $\Expectation{K_z^2} \leq c \eta(z)^{-(d - 2)} \leq c' \Expectation{K_z},$ as wanted.
		\end{proof}

		\begin{Lemma}\label{Lemma: Expectation of Kz is proportional to eta of z to the power of negative d minus two}
			Suppose that $\varphi$ and $\psi$ are two bounded functions defined on $\Z^{d + 1}$ and with positive values. Furthermore, suppose that there exists two constants $c_1, c_2 > 0$ such that
			\begin{displaymath}\textstyle
				\max(\varphi(z), \psi(z)) \leq
				\begin{cases}
					c_1 e^{-c_2 \norm{z}} &\text{ if } \norm{x} > n, \\
					c_1 e^{-c_2 \frac{\norm{x}^2}{n}} \norm{z}^{-\frac{d}{2}} &\text{ if } \norm{x} \leq n.
				\end{cases}
			\end{displaymath}
			Then, there exists a universal constant $c > 0$ such that for any $z \in \Z^{d + 1},$ $\sum\limits_{z' \in \Z^{d + 1}} \varphi(z' - z) \psi(z') \leq c \eta(z)^{-(d - 2)}.$ Likewise, if there exists two constants $c_3, c_4 > 0$ such that
			\begin{displaymath}\textstyle
				\min(\varphi(z), \psi(z)) \geq
				\begin{cases}
					c_3 e^{-c_4 \norm{z}} &\text{ if } \norm{x} > n, \\
					c_3 e^{-c_4 \frac{\norm{x}^2}{n}} \norm{z}^{-\frac{d}{2}} &\text{ if } \norm{x} \leq n,
				\end{cases}
			\end{displaymath}
			then there exists a universal constant $c > 0$ such that for any $z \in \Z^{d + 1},$ $\sum\limits_{z' \in \Z^{d + 1}} \varphi(z' - z) \psi(z') \geq c \eta(z)^{-(d - 2)}.$
		\end{Lemma}
		\begin{proof}[Proof of Lemma]
		We divide the proof of the upper and lower bounds. We write $z = (n, x),$ $z' = (n', x')$ and so on with other affixes (such as superscripts).

		\noindent \textsc{I. Upper bound.}
			There is no loss of generality to assume $\varphi = \psi.$ Decompose $\sum\limits_{z' \in \Z^{d + 1}} \varphi(z' - z) \varphi(z') = H + I + J,$ where
			\begin{equation}\textstyle\label{SubEquation: Definition of the sums I, H and J}
				H = \sum\limits_{n' \leq 0} \sum\limits_{x' \in \Z^d} \varphi(z' - z) \varphi(z'), \quad I = \sum\limits_{1 \leq n' \leq n - 1} \sum\limits_{x' \in \Z^d} \varphi(z' - z) \varphi(z'), \quad J = \sum\limits_{n' \geq n} \sum\limits_{x' \in \Z^d} \varphi(z' - z) \varphi(z')
			\end{equation}
			We shall bound now each of the terms $H, I, J$ and we will see, at the end, that the largest of these bounds is of the claimed form. We will use the bounds in the hypothesis and the estimates of \S \ref{Section: Estimates of sums} to handle the resulting sums.

			We start by bounding $H.$ Since $z \in \Z_+ \times \Z^d$ and $z' \in \Z_- \times \Z^d,$ we have the bounds $\varphi(z' - z) \varphi(z') \leq c e^{-c'(\norm{z - z'} + \norm{z'})},$ for a pair of positive constants. Now,
			\begin{displaymath}\textstyle
				c e^{-c'(\norm{z - z'} + \norm{z'})} \leq c e^{-\frac{c'}{2}(\norm{z} - \norm{z'}) - c'\norm{z'}} \leq c e^{-\frac{c'}{2} \norm{z}} e^{-\frac{c'}{2} \norm{z'}},
			\end{displaymath}
			summing over $z' \in \Z_- \times \Z^d,$ we reach
			\begin{equation}\textstyle\label{SubEquation: Upper bound of H}
				H \leq c e^{-c' \norm{z}},
			\end{equation}
			for a pair of constants $c, c' > 0.$

			We now handle $I$ from (\ref{SubEquation: Definition of the sums I, H and J}). Divide into two parts as follows
			\begin{displaymath}\textstyle
				I = I_1 + I_2 \mathop{=}\limits^\mathrm{def.} \sum\limits_{\substack{1 \leq n' \leq n - 1 \\ \norm{x'} \leq n'}} \varphi(z' - z) \varphi(z') + \sum\limits_{\substack{1 \leq n' \leq n - 1 \\ \norm{x'} > n'}} \varphi(z' - z) \varphi(z').
			\end{displaymath}
			Then, $I_1 \leq c \sum\limits_{\substack{1 \leq n' \leq n - 1 \\ \norm{x'} \leq n'}} \exp\left(-c'\left( \norm{z - z'} + \frac{\norm{x'}^2}{n'} \right) \right) (n')^{-\frac{d}{2}}.$ Notice $\norm{z - z'}^2 = (n - n')^2 + \norm{x - x'}^2 \geq \frac{\left( |n - n'| + \norm{x - x'} \right)^2}{2},$ thus $\norm{z - z'} \geq \frac{|n - n'| + \norm{x - x'}}{\sqrt{2}}.$ Hence,
			\begin{equation}\textstyle\label{SubEquation: First bound of the integral I1}
				I_1 \leq c \sum\limits_{n' = 1}^{n - 1} (n')^{-\frac{d}{2}} e^{-\frac{c'}{\sqrt{2}}(n - n')} \sum\limits_{\norm{x'} \leq n'} \exp\left(-\frac{c'}{\sqrt{2}}\left( \norm{x - x'} + \frac{\norm{x'}^2}{n'} \right) \right).
			\end{equation}
			For convenience, we now divide into two cases. But first, by (\ref{Proposition: Sums of Greens function, ver 1}), we can bound
			\begin{displaymath}\textstyle
				\sum\limits_{\norm{x'} \leq n'} \exp\left(-\frac{c'}{\sqrt{2}}\left( \norm{x - x'} + \frac{\norm{x'}^2}{n'} \right) \right) \leq
				\begin{cases}
					c e^{-c' \norm{x}} &\text{ if } \norm{x} > n' \\
					c e^{-c'\frac{\norm{x}^2}{n'}} &\text{ if } \norm{x} \leq n'.
				\end{cases}
			\end{displaymath}
			We may now proceed with the two cases.

					\begin{description}
				\item[Case 1.] Here we assume $n \leq \norm{x}.$ Then, there exists two constants $c, c' > 0,$
				\begin{displaymath}\textstyle
					I_1 \leq c e^{-c' \norm{x}} \sum\limits_{n' = 1}^{n - 1} (n')^{-\frac{d}{2}} e^{-c'(n - n')} \leq c \left[ \sum\limits_{k = 1}^\infty k^{-\frac{d}{2}} \right] e^{-c' \norm{x}} \leq C e^{-C' \norm{z}}.
				\end{displaymath}
		
				\item[Case 2.] Here we assume $\norm{x} < n.$ Then, we divide the sum in the right of (\ref{SubEquation: First bound of the integral I1}) into two parts. The first part consists on all terms for which	 $1 \leq n' \leq \norm{x}$ and the second, $\norm{x} < n' \leq n - 1.$ Then, for a pair of constants $C, C' > 0,$
				\begin{displaymath}\textstyle
					\sum\limits_{1 \leq n' \leq \norm{x}} (n')^{-\frac{d}{2}} e^{-\frac{c'}{\sqrt{2}}(n - n')} \sum\limits_{\norm{x'} \leq n'} e^{-\frac{c'}{\sqrt{2}}\left( \norm{x - x'} + \frac{\norm{x'}^2}{n'} \right)} \leq C \sum\limits_{1 \leq n' \leq \norm{x}} (n')^{-\frac{d}{2}} e^{-\frac{c'}{\sqrt{2}}(n - n')} e^{-C' \norm{x}}.
				\end{displaymath}
				By setting $c'' = \min\left( \frac{c'}{\sqrt{2}}, C' \right),$ it is clear that $e^{-\frac{c'}{\sqrt{2}}(n - n')} e^{-C' \norm{x}} \leq e^{-c''(n - n' + \norm{x})} \leq e^{-c'' n},$ 	for $1 \leq n' \leq \norm{x}.$ Then, the first part in the division of $I_1$ is bounded above by $c e^{-c' \norm{z}},$ for a pair of constants $c, c' > 0.$ As for the second part,
				\begin{displaymath}
					\begin{aligned}\textstyle
						\sum\limits_{\norm{x} < n' \leq n - 1} (n')^{-\frac{d}{2}} e^{-\frac{c'}{\sqrt{2}}(n - n')} &\textstyle\sum\limits_{\norm{x'} \leq n'} e^{-\frac{c'}{\sqrt{2}}\left( \norm{x - x'} + \frac{\norm{x'}^2}{n'} \right)} \\
						&\textstyle\leq C \sum\limits_{\norm{x} < n' \leq n - 1} (n')^{-\frac{d}{2}} e^{-\frac{c'}{\sqrt{2}}(n - n')} e^{-C' \frac{\norm{x}^2}{n'}},
					\end{aligned}
				\end{displaymath}
				where $C, C' > 0$ are two constants. We have shown up to this point that under the assumption $\norm{x} < n,$ there exists two constants $c, c' > 0,$ such that
				\begin{equation}\textstyle\label{SubEquation: Second bound of the integral I1}
					I_1 \leq ce^{-c' \norm{z}} + c \sum\limits_{\norm{x} < n' \leq n - 1} (n')^{-\frac{d}{2}} e^{-c' \left( n - n' + \frac{\norm{x}^2}{n'} \right)}.
				\end{equation}
				We are going to consider two further cases for the sum appearing in the previous line.
				\begin{enumerate}
					\item Assume first $\norm{x} \leq n^{\frac{1}{2}}.$ Here, we will bound $e^{-c \frac{\norm{x}^2}{n'}} \leq 1.$ Then,
					\begin{displaymath}\textstyle
						\sum\limits_{\norm{x} < n' \leq n - 1} (n')^{-\frac{d}{2}} e^{-c' \left( n - n' + \frac{\norm{x}^2}{n'} \right)} = \sum\limits_{\norm{x} < n' \leq n^{\frac{1}{2}}} (n')^{-\frac{d}{2}} e^{-c' \left( n - n' \right)} + \sum\limits_{n^{\frac{1}{2}} < n' \leq n} (n')^{-\frac{d}{2}} e^{-c' \left( n - n' \right)}.
					\end{displaymath}
					In the first of these sums, $e^{-c'(n - n')} \leq e^{-c'\left(n - n^{\frac{1}{2}} \right)} \leq C e^{-C' n^{\frac{1}{2}}},$ for two constants $C, C' > 0.$ We reach $\sum\limits_{\norm{x} < n' \leq n^{\frac{1}{2}}} (n')^{-\frac{d}{2}} e^{-c' \left( n - n' \right)} \leq C \left[ \sum\limits_{k = 1}^\infty k^{-\frac{d}{2}} \right] e^{-C' n^{\frac{1}{2}}}.$ Observe now that if $n^{\frac{1}{2}} < n' \leq n,$ there exists a unique integer $j$ between $1$ and $n^{\frac{1}{2}}$ such that $j n^{\frac{1}{2}} < n' \leq (j + 1) n^{\frac{1}{2}}$ (there are $n^{\frac{1}{2}}$ indices $n'$ in this range). Thus,
					\begin{align*}\textstyle
						\sum\limits_{n^{\frac{1}{2}} < n' \leq n} (n')^{-\frac{d}{2}} e^{-c' \left( n - n' \right)} 
						&\textstyle\leq n^{-\frac{d}{4}} n^{\frac{1}{2}} \sum\limits_{1 \leq j \leq n^{\frac{1}{2}}} j^{-\frac{d}{2}} e^{-c' \left( n - \min \left( (j + 1)n^{\frac{1}{2}}, n \right) \right)} \\
						&\textstyle= n^{-\frac{d}{4}+\frac{1}{2}} \sum\limits_{1 \leq j \leq n^{\frac{1}{2}}} j^{-\frac{d}{2}} e^{-c' n^{\frac{1}{2}} \left( n^{\frac{1}{2}} - \min \left( j + 1, n^{\frac{1}{2}} \right) \right)}.
					\end{align*}
					In the last sum, consider first the range $1 \leq j \leq \frac{n^{\frac{1}{2}}}{4}$ to obtain $j + 1 \leq 2j \leq \frac{n^{\frac{1}{2}}}{2},$ giving $e^{-c' n^{\frac{1}{2}} \left( n^{\frac{1}{2}} - \min \left( j + 1, n^{\frac{1}{2}} \right) \right)} \leq e^{-\frac{c'}{2} n}.$ Next, consider the range $\frac{n^{\frac{1}{2}}}{4} < j \leq n^{\frac{1}{2}},$ so that $j \asymp n^{\frac{1}{2}}.$ Then,
					\begin{displaymath}\textstyle
						\sum\limits_{n^{\frac{1}{2}} < n' \leq n} (n')^{-\frac{d}{2}} e^{-c' \left( n - n' \right)}  \leq \left[ \sum\limits_{k = 1}^\infty k^{-\frac{d}{2}} \right] n^{-\frac{d}{4} + \frac{1}{2}} e^{-\frac{c'}{2} n} + c n^{-\frac{d}{4} + \frac{1}{2}} n^{-\frac{d}{4}} n^\frac{1}{2} \leq C n^{-\frac{d}{2} + 1}.
					\end{displaymath}
					The foregoing shows $\sum\limits_{\norm{x} \leq n' < n - 1} (n')^{-\frac{d}{2}} e^{-c' \left( n - n' + \frac{\norm{x}^2}{n'} \right)} \leq c e^{-c' n} + C n^{-\frac{d}{2} + 1} \leq C'' n^{-\frac{d}{2} + 1}.$
		
					\item Now assume $n^{\frac{1}{2}} < \norm{x}.$ We bound the sum appearing in (\ref{SubEquation: Second bound of the integral I1}). We use (\ref{Proposition: Estimates in the p series with an exponential factor}),
					\begin{displaymath}\textstyle
						\sum\limits_{\norm{x} < n' \leq n - 1} (n')^{-\frac{d}{2}} e^{-c' \left( n - n' + \frac{\norm{x}^2}{n'} \right)} \leq \sum\limits_{n' \leq \norm{x}^2} (n')^{-\frac{d}{2}} e^{-c' \frac{\norm{x}^2}{n'}} \leq c \norm{x}^{-(d - 2)}.
					\end{displaymath}
				\end{enumerate}
				Bearing in mind the previous two items, and (\ref{SubEquation: Second bound of the integral I1}), we have established that $I_1 \leq c \eta(z)^{-(d - 2)},$ for a universal constant $c > 0.$ This completes \textbf{Case 2.}
			\end{description}

			We now need to bound $I_2.$ Again, by the bounds in the hypothesis,
			\begin{displaymath}\textstyle
				I_2 = \sum\limits_{1 \leq n' \leq n - 1} \sum\limits_{\norm{x'} > n'} \varphi(z' - z) \varphi(z') \leq c_1^2 \sum\limits_{1 \leq n' \leq n - 1} \sum\limits_{\norm{x'} > n'} e^{-c_2(\norm{z - z'} + \norm{z'})} \leq c e^{-c' \norm{z}},
			\end{displaymath}
			where we proceeded as in (\ref{SubEquation: Upper bound of H}). Thus, summing the upper bounds of $I_1$ and $I_2$ we have shown that
			\begin{equation}\textstyle\label{SubEquation: Upper bound for I}
				I \leq c \eta(z)^{-(d - 2)},
			\end{equation}
			for a universal constant $c > 0,$ as desired.
		
			We now proceed to upper bound the sum $J$ of (\ref{SubEquation: Definition of the sums I, H and J}). First we rewrite, $J = \sum\limits_{n' \geq 0} \sum\limits_{x' \in \Z^d} \varphi(0, z') \varphi(0, z' + z) = J_1 + J_2,$ where $J_1 = \sum\limits_{n' \geq 0} \sum\limits_{\norm{x'} \leq n'} \varphi(0, z') \varphi(0, z' + z)$ and $J_2 = \sum\limits_{n' \geq 0} \sum\limits_{\norm{x'} > n'} \varphi(0, z') \varphi(0, z' + z).$ We further divide $J_1$ as follows: $J_1 = J_{1, 1} + J_{1, 2},$ with $J_{1, 1} = \sum\limits_{n' \geq 0} \sum\limits_{ \substack{ \norm{x'} \leq n' \\ \norm{x + x'} \leq n + n' } } \varphi(0, z') \varphi(0, z' + z)$ and $J_{1, 2} = \sum\limits_{n' \geq 0} \sum\limits_{ \substack{ \norm{x'} \leq n' \\ \norm{x + x'} > n + n' } } \varphi(0, z') \varphi(0, z' + z).$		
			Let us handle $J_{1, 1}$ first. In this case, we have the bound
			\begin{displaymath}\textstyle
				J_{1, 1} \leq c_1^2 \sum\limits_{n' \geq 0} (n')^{-\frac{d}{2}} (n + n')^{-\frac{d}{2}} \sum\limits_{ \substack{ \norm{x'} \leq n' \\ \norm{x + x'} \leq n + n' } } e^{-c_2 \left( \frac{\norm{x'}^2}{n'} + \frac{\norm{x + x'}^2}{n + n'} \right)}.
			\end{displaymath}
			We finish the case $\norm{x} \leq n^{\frac{1}{2}}$ first. Here we may bound
			\begin{displaymath}\textstyle
				J_{1, 1} \leq c_1^2 \sum\limits_{n' \geq 0} (n')^{-\frac{d}{2}} (n + n')^{-\frac{d}{2}} \sum\limits_{\norm{x'} \leq n'} e^{-c_2 \frac{\norm{x'}^2}{n'}} \leq c \sum\limits_{n' \geq 0} (n + n')^{-\frac{d}{2}} \leq c' n^{-\frac{d}{2} + 1},
			\end{displaymath}
			by (\ref{Proposition: Asymptotic estimates on Eulerian sums}) and (\ref{Proposition: Estimates in the p series}) with all $c, c', c''$ positive constants. We may assume for the rest of the bounding of $J_{1, 1}$ that $\norm{x} > n^{\frac{1}{2}}.$ By virtue of (\ref{Proposition: Sums of Greens function, ver 2}), in the range $0 \leq n' \leq \frac{\norm{x}}{2},$ we may bound
			\begin{displaymath}\textstyle
				\sum\limits_{ \substack{ \norm{x'} \leq n' \\ \norm{x + x'} \leq n + n' } } e^{-c_2 \left( \frac{\norm{x'}^2}{n'} + \frac{\norm{x + x'}^2}{n + n'} \right)} \leq c (n')^{\frac{d}{2}} e^{-c \frac{\norm{x}^2}{n + n'}},
			\end{displaymath}
			with $c$ a universal constant. Next,
			\begin{displaymath}\textstyle
				\sum\limits_{0 \leq n' \leq \frac{\norm{x}}{2}} (n')^{-\frac{d}{2}} (n + n')^{-\frac{d}{2}} (n')^{\frac{d}{2}} e^{-c' \frac{\norm{x}^2}{n + n'}} = c \sum\limits_{0 \leq n' \leq \frac{\norm{x}}{2}} (n + n')^{-\frac{d}{2}} e^{-c' \frac{\norm{x}^2}{n + n'}} = c \sum\limits_{n \leq m \leq n + \frac{\norm{x}}{2}} m^{-\frac{d}{2}} e^{-c' \frac{\norm{x}^2}{m}}.
			\end{displaymath}
			By the assumption $\norm{x} > n^{\frac{1}{2}},$ we have $\sum\limits_{n \leq m \leq n + \frac{\norm{x}}{2}} m^{-\frac{d}{2}} e^{-c' \frac{\norm{x}^2}{m}} \leq \sum\limits_{m \leq 2 \norm{x}^2} m^{-\frac{d}{2}} e^{-c' \frac{\norm{x}^2}{m}} \leq c \norm{x}^{-(d - 2)},$ by (\ref{Proposition: Estimates in the p series with an exponential factor}). To sum up, we have established so far in the bounding of $J_{1, 1}$ (with the assumption that $\norm{x} > n^{\frac{1}{2}}$),
			\begin{displaymath}\textstyle
			J_{1, 1} \leq c \eta(z)^{-(d - 2)} + c \sum\limits_{n' \geq \frac{\norm{x}}{2}} (n')^{-\frac{d}{2}} (n + n')^{-\frac{d}{2}} \sum\limits_{ \substack{ \norm{x'} \leq n' \\ \norm{x + x'} \leq n + n' } } e^{-c_2 \left( \frac{\norm{x'}^2}{n'} + \frac{\norm{x + x'}^2}{n + n'} \right)}
			\end{displaymath}

			If $\norm{x} \leq 2n',$ we obtain from (\ref{Proposition: Sums of Greens function, ver 2})
			\begin{displaymath}\textstyle
				\sum\limits_{ \substack{ \norm{x'} \leq n' \\ \norm{x + x'} \leq n + n' } } e^{-c_2 \left( \frac{\norm{x'}^2}{n'} + \frac{\norm{x + x'}^2}{n + n'} \right)} \leq c \left[ (n')^{\frac{d}{2}} e^{-c' \frac{\norm{x}^2}{n + n'}} + (n + n')^{\frac{d}{2}} e^{-c' \frac{\norm{x}^2}{n'}} \right].
			\end{displaymath}
			Substituting this result above,
			\begin{align*}
				J_{1, 1} &\textstyle\leq c \eta(z)^{-(d - 2)} + c \sum\limits_{n' \geq \frac{\norm{x}}{2}} \left[ (n + n')^{-\frac{d}{2}} e^{-c' \frac{\norm{x}^2}{n + n'}} + (n')^{-\frac{d}{2}} e^{-c' \frac{\norm{x}^2}{n'}} \right] \\
				&\textstyle\leq c \eta(z)^{-(d - 2)} + 2c \sum\limits_{n' > \frac{\norm{x}}{2}} (n')^{-\frac{d}{2}} e^{-c' \frac{\norm{x}^2}{n'}}.
			\end{align*}
			Split the appearing sum in the line above into $\frac{\norm{x}}{2} \leq n' \leq \norm{x}^2$ and $n' > \norm{x}^2;$ bound the first resulting sum by (\ref{Proposition: Estimates in the p series with an exponential factor}) to obtain an upper bound of the form $c \norm{x}^{-(d - 2)},$ and the second, by $\sum\limits_{n' > \norm{x}^2} (n')^{-\frac{d}{2}} e^{-c' \frac{\norm{x}^2}{n'}} \leq  \sum\limits_{n' > \norm{x}^2} (n')^{-\frac{d}{2}}  \asymp \norm{x}^{-(d - 2)},$ the last $\asymp$ was obtained from (\ref{Proposition: Estimates in the p series}). We have established,
			\begin{equation}\textstyle\label{SubEquation: Upper bound for J11}
				J_{1, 1} \leq c \eta(z)^{-(d - 2)},
			\end{equation}
			$c$ being a universal constant depending solely on dimension.

			We now will handle $J_{1, 2}.$ The hypotheses give
			\begin{align*}
				J_{1, 2} &\textstyle\leq c_1^2 \sum\limits_{n' \geq 0} \sum\limits_{ \substack{ \norm{x'} \leq n' \\ \norm{x + x'} > n + n' } } (n')^{-\frac{d}{2}} e^{-c_2 \left( \norm{z + z'} + \frac{\norm{x'}^2}{n'} \right)} \\
				&\textstyle\leq c_1^2 e^{-\frac{c_2}{\sqrt{2}} n} \sum\limits_{n' \geq 0} (n')^{-\frac{d}{2}} e^{-\frac{c_2}{\sqrt{2}} n'} \sum\limits_{\norm{x'} \leq n'} e^{-\frac{c_2}{\sqrt{2}} \left( \norm{x + x'} + \frac{\norm{x'}^2}{n'} \right)}.
			\end{align*}
			It is now clear that we may bound $J_{1, 2} \leq c e^{-c' n},$ with a pair of constants $c, c' > 0$ depending solely on dimension. This establishes the case $n^{\frac{1}{2}} \geq \norm{x}$ and thus, we may assume now $\norm{x} > n^{\frac{1}{2}}$ in what follows and bound $e^{-\frac{c_2}{\sqrt{2}} n} \leq 1.$ Consider an integer $0 \leq n' < \norm{x},$ invoking (\ref{Proposition: Sums of Greens function, ver 1}) we reach $\sum\limits_{\norm{x'} \leq n'} e^{-\frac{c_2}{\sqrt{2}} \left( \norm{x + x'} + \frac{\norm{x'}^2}{n'} \right)} \leq c e^{-c' \norm{x}},$ \emph{a fortiori} $c_1^2 \sum\limits_{0 \leq n' < \norm{x}} (n')^{-\frac{d}{2}} \sum\limits_{\norm{x'} \leq n'} e^{-\frac{c_2}{\sqrt{2}} \left( \norm{x + x'} + \frac{\norm{x'}^2}{n'} \right)} \leq c e^{-c' \norm{x}},$ for universal constants $c, c' > 0.$ If $n' \geq \norm{x},$ bound $\norm{x + x'} \geq 0$ and use (\ref{Proposition: Asymptotic estimates on Eulerian sums}): $\sum\limits_{\norm{x'} \leq n'} e^{-\frac{c_2}{\sqrt{2}} \left( \norm{x + x'} + \frac{\norm{x'}^2}{n'} \right)} \leq c (n')^\frac{d}{2},$ and then $\sum\limits_{n' \geq \norm{x}} (n')^{-\frac{d}{2}} e^{-\frac{c_2}{\sqrt{2}} n'} \sum\limits_{\norm{x'} \leq n'} e^{-\frac{c_2}{\sqrt{2}} \left( \norm{x + x'} + \frac{\norm{x'}^2}{n'} \right)} \leq c \sum\limits_{n' \geq \norm{x}} e^{-\frac{c_2}{\sqrt{2}} n'} \asymp e^{-\frac{c_2}{\sqrt{2}} \norm{x}}.$ Thus, we may assert, there exists a constant $c > 0$ such that for all $z =(n, x)$ satisfying $\norm{x} > n^{\frac{1}{2}},$ $J_{1, 2} \leq c \eta(z)^{-(d - 2)};$ the same type of bound for $\norm{x} \leq n^{\frac{1}{2}}$ was established already. Bearing in mind these bounds, (\ref{SubEquation: Upper bound for J11}), and since $J_1 = J_{1, 1} + J_{1, 2},$ we have finally established
			\begin{equation}\textstyle\label{SubEquation: Upper bound for J1}
				J_1 \leq c \eta(z)^{-(d - 2)},
			\end{equation}
			with $c > 0$ is a constant that does not depend on $z.$
		
			We will now prove the upper bound for $J_2.$ Similarly as was done with $J_1,$ we split $J_2 = J_{2, 1} + J_{2, 2},$ being $J_{2, 1} = \sum\limits_{n' \geq 0} \sum\limits_{ \substack{ \norm{x'} > n' \\ \norm{x + x'} \leq n + n' } } \varphi(z') \varphi(z' + z)$ and $J_{2, 2} = \sum\limits_{n' \geq 0} \sum\limits_{ \substack{ \norm{x'} > n' \\ \norm{x + x'} > n + n' } } \varphi(z') \varphi(z' + z).$ Here, $J_{2, 2}$ is easy for we have for we may repeat the calculations of (\ref{SubEquation: Upper bound of H}) and obtain
			\begin{equation}\textstyle\label{SubEquation: Upper bound for J22}
				J_{2, 2} \leq c_1^2 \sum\limits_{z' \in \Z^{d + 1}} e^{-c_2(\norm{z'} + \norm{z + z'})} \leq c e^{-c' \norm{z}},
			\end{equation}
			with $c_1, c_2$ the constants in the hypothesis and $c, c'$ two positive universal constants. There remains to handle $J_{2, 1}.$ Here we may bound
			\begin{align*}
				J_{2, 1} &\textstyle\leq c_1^2 \sum\limits_{n' \geq 0} \sum\limits_{\norm{x + x'} \leq n + n'} e^{-c_2 \left( \norm{z'} + \frac{\norm{x + x'}^2}{n + n'} \right)} (n + n')^{-\frac{d}{2}} \\
				&\textstyle\leq c_1^2 \sum\limits_{n' \geq 0} (n + n')^{-\frac{d}{2}} e^{-\frac{c_2}{\sqrt{2}} n'} \sum\limits_{\norm{x + x'} \leq n + n'} e^{-\frac{c_2}{\sqrt{2}} \left( \norm{x'} + \frac{\norm{x + x'}^2}{n + n'} \right)} \\
				&\textstyle= c_1^2 \sum\limits_{n' \geq 0} (n + n')^{-\frac{d}{2}} e^{-\frac{c_2}{\sqrt{2}} n'} \sum\limits_{\norm{x'} \leq n + n'} e^{-\frac{c_2}{\sqrt{2}} \left( \norm{x' - x} + \frac{\norm{x'}^2}{n + n'} \right)}.
			\end{align*}
			Suppose first that $\norm{x} \leq n^{\frac{1}{2}}.$ The bounds of (\ref{Proposition: Sums of Greens function, ver 1}) give $\sum\limits_{\norm{x'} \leq n + n'} e^{-\frac{c_2}{\sqrt{2}} \left( \norm{x' - x} + \frac{\norm{x'}^2}{n + n'} \right)} \leq c e^{-c' \frac{\norm{x}^2}{n + n'}} \leq c,$ for a pair of constants $c, c' > 0.$ Then, $J_{2, 1} \leq c \sum\limits_{n' \geq 0} (n + n')^{-\frac{d}{2}} \leq c' n^{-\frac{d}{2} + 1} = c' \eta(z)^{-(d - 2)},$ where the constants $c, c'$ are positive and depending on nothing but, perhaps, dimension. This would finish the bound for $J_{2, 1}$ under the assumption $\norm{x} \leq n^{\frac{1}{2}}.$ Next, assume $n^{\frac{1}{2}} < \norm{x}.$ Here, we further split (we remind the reader the sums are over $n' \in \Z_+$ and $x' \in \Z^d$ as stated)
			\begin{align*}
				J_{2, 1} &\textstyle\leq c_1^2 \sum\limits_{0 \leq n' + n \leq \norm{x}} (n + n')^{-\frac{d}{2}} e^{-\frac{c_2}{\sqrt{2}} n'} \sum\limits_{\norm{x'} \leq n + n'} e^{-\frac{c_2}{\sqrt{2}} \left( \norm{x' - x} + \frac{\norm{x'}^2}{n + n'} \right)} \\
				&\textstyle+ c_1^2 \sum\limits_{n' + n > \norm{x}} (n + n')^{-\frac{d}{2}} e^{-\frac{c_2}{\sqrt{2}} n'} \sum\limits_{\norm{x'} \leq n + n'} e^{-\frac{c_2}{\sqrt{2}} \left( \norm{x' - x} + \frac{\norm{x'}^2}{n + n'} \right)}.
			\end{align*}
			If $\norm{x} \geq n' + n$ then the inner sum is bounded by an expression of the form $c e^{-c' \norm{x}},$ again by (\ref{Proposition: Sums of Greens function, ver 1}). This yields $c_1^2 \sum\limits_{0 \leq n' + n \leq \norm{x}} (n + n')^{-\frac{d}{2}} e^{-\frac{c_2}{\sqrt{2}} n'} \sum\limits_{\norm{x'} \leq n + n'} e^{-\frac{c_2}{\sqrt{2}} \left( \norm{x' - x} + \frac{\norm{x'}^2}{n + n'} \right)} \leq c e^{-c' \norm{x}} \leq C \eta(z)^{-(d - 2)},$ with $c, c', C > 0$ universal constants. If $\norm{x} < n + n',$ then the inner sum is bounded by an expression of the form $c e^{-c' \frac{\norm{x}^2}{n + n'}},$ by the same point. Thus, there exist two constants $c, c'$ such that for all $z,$
			\begin{displaymath}\textstyle
					c_1^2 \sum\limits_{n' + n > \norm{x}} (n + n')^{-\frac{d}{2}} e^{-\frac{c_2}{\sqrt{2}} n'} \sum\limits_{\norm{x'} \leq n + n'} e^{-\frac{c_2}{\sqrt{2}} \left( \norm{x' - x} + \frac{\norm{x'}^2}{n + n'} \right)} \leq c \sum\limits_{n' + n > \norm{x}} (n + n')^{-\frac{d}{2}} e^{-c' \left(n' + \frac{\norm{x}^2}{n + n'} \right)}.
			\end{displaymath}
			We may now use (\ref{Proposition: Estimates in the p series}) and (\ref{Proposition: Estimates in the p series with an exponential factor}) as follows:
			\begin{align*}
				\textstyle\sum\limits_{n' + n >\norm{x}} (n + n')^{-\frac{d}{2}} e^{-c' \left(n' + \frac{\norm{x}^2}{n + n'} \right)}
				&\textstyle\leq \sum\limits_{n' + n \leq \norm{x}^2} (n + n')^{-\frac{d}{2}} e^{-c' \frac{\norm{x}^2}{n + n'}} + \sum\limits_{n' + n > \norm{x}^2} (n + n')^{-\frac{d}{2}} \\
				&\leq c \norm{x}^{-(d - 2)} + c' \norm{x}^{-(d - 2)} = (c + c') \norm{x}^{-(d - 2)}.
			\end{align*}
			Thus, we showed that there exists a constant $c > 0$ such that for all $z \in \Z^{d + 1},$
			\begin{equation}\textstyle\label{SubEquation: Upper bound for J21}
				J_{2, 1} \leq c \eta(z)^{-(d - 2)}.
			\end{equation}
			Then, by (\ref{SubEquation: Upper bound for J22}) and (\ref{SubEquation: Upper bound for J21}), we reach $J_2 \leq c \eta(z)^{-(d - 2)},$ with $c$ not depending on $z.$ Furthermore, the previous bound on $J_2$ together with (\ref{SubEquation: Upper bound for J1}) show that there exists a $c > 0$ such that  for all vertices $z$ in $\Gamma_d(\lambda),$
			\begin{equation}\textstyle\label{SubEquation: Upper bound of J}
				J \leq c \eta(z)^{-(d - 2)}.
			\end{equation}
			Combining the results of (\ref{SubEquation: Upper bound of H}), (\ref{SubEquation: Upper bound for I}) and (\ref{SubEquation: Upper bound of J}), we finally reach that $\sum\limits_{z' \in \Z^{d + 1}} \varphi(z' - z) \varphi(z') \leq c \eta(z)^{-\frac{d}{2} + 1}.$

			\noindent \textsc{II. Lower bound.} Assume the hypothesis for the lower bound. Then
			\begin{displaymath}\textstyle
				\sum\limits_{z' \in \Z^{d + 1}} \varphi(z' - z) \psi(z') \geq \sum\limits_{n' > \max \left( n, \norm{x}^2 \right) } \sum\limits_{\norm{x' - x}^2 \leq n' - n} \varphi(z' - z) \psi(z').
			\end{displaymath}
			The two relations $\norm{x}^2 < n'$ and $\norm{x' - x}^2 \leq n' - n$ imply $\norm{x'}^2 \leq 2 \norm{x}^2 + 2\norm{x' - x}^2 \leq 4n'.$ Then, the bounds of the hypothesis give (the last $\asymp$ follows from (\ref{Proposition: Estimates in the p series}))
			\begin{displaymath}\textstyle
				\begin{aligned}
					\textstyle\sum\limits_{n' > \max \left( n, \norm{x}^2 \right)}
					&\textstyle\sum\limits_{\norm{x' - x}^2 \leq n' - n} \varphi(z' - z) \psi(z') \\
					&\textstyle\geq c \sum\limits_{n' > \max \left( n, \norm{x}^2 \right)} (n')^{-\frac{d}{2}} (n' - n)^{-\frac{d}{2}} \sum\limits_{\norm{x' - x}^2 \leq n' - n} e^{-c' \frac{\norm{x'}^2}{n'}} \\
					&\textstyle\geq c \sum\limits_{n' > \max \left( n, \norm{x}^2 \right)} (n')^{-\frac{d}{2}} (n' - n)^{-\frac{d}{2}} \sum\limits_{\norm{x' - x}^2 \leq n' - n} e^{-4c'} \\
					&\textstyle\geq c'' \sum\limits_{n' > \max \left( n, \norm{x}^2 \right)} (n')^{-\frac{d}{2}} \asymp \max \left( n, \norm{x}^2 \right)^{-\frac{d}{2} + 1}.
				\end{aligned}
			\end{displaymath}
			The lower bounds have been substantiated keeping in mind (\ref{SubEquation: Bounds above and below of the probability of same component}).
		\end{proof}

		With the two lemmas proved, we are in conditions to finish the theorem. By virtue of (\ref{SubEquation: Bounds above and below of the probability of same component}), (\ref{SubEquation: Definition of Kz}), (\ref{SubEquation: Expectation of Kz as sums of Greens function}) and (\ref{Lemma: The probability of Kz is positive is asymptotic to its expectation}) we have that $\USF(\scrC(0, z)) \asymp \Expectation{K_z} = \sum\limits_{z' \in \Z^{d + 1}} G(z, z') G(0, z').$ Applying (\ref{Lemma: Expectation of Kz is proportional to eta of z to the power of negative d minus two}) (bear in mind the Green's function bounds (\ref{Theorem: Greens function bounds})), we reach that $\USF(\scrC(0,z)) \asymp \eta(z)^{-(d - 2)}.$
	\end{proof}
	\numberwithin{SubCounter}{equation}
	\numberwithin{equation}{section}
	\setcounter{equation}{\value{auxCounter}}

	\subsection[Upper bounds on the probability of $\USF$-connectedness]{An upper bound on the probability that finitely many points are $\USF$-connected}
	We are going to generalise the previous result to any subset of vertices. We will follow closely \cite{BKPS:GeomUSF}. First, we introduce the definition of \emph{spread} of a set of vertices. For any vertex $z \in \Z^{d + 1}$ and any finite subset $\strW \subset \Z^{d + 1}$ (with two or more vertices) we define their spreads (relative to the metric $\eta$) as
	\begin{equation}\textstyle
		\begin{split}
			\spread{z} &= \max(1, \eta(z)) \\
			\spread{\strW} &\textstyle= \min\limits_{\strE} \prod\limits_{\{z, z'\} \in \strE} \spread{z - z'},
		\end{split}
	\end{equation}
	where the minimum runs over \emph{all $\strE \subset \scrP(\strW)$ making $(\strW, \strE)$ an undirected spanning tree.} If $\strW$ consists of two (different) vertices $\alpha$ and $\beta,$ then $\spread{\strW} = \eta(\alpha - \beta);$ in this case we will write $\spread{\alpha \beta}$ in lieu of $\spread{\strW}.$ Observe that (\ref{Equation: Japanese bracket}) defines the spread of a point relative to the Euclidean metric. For simplicity, if $\strW = \{\alpha, \beta, \ldots, \xi\},$ we will write $\spread{\alpha\beta \cdots \xi},$ in particular, $\spread{0z} = \spread{z}.$ Observe that \cite{BKPS:GeomUSF} also defined the concept of ``spread of subset of vertices of $\Z^d$'' but they use the Euclidean metric (see their Definition 2.1). When $\strW = \{\alpha, \beta, \gamma\},$ then
	\begin{displaymath}\textstyle
		\spread{\strW} = \min\{ \spread{\alpha \beta} \spread{\alpha \gamma}, \spread{\alpha \beta} \spread{\beta \gamma}, \spread{\alpha \gamma} \spread{\beta \gamma}\}.
	\end{displaymath}

	\begin{Proposition}{\normalfont \cite[Lemma 2.6]{BKPS:GeomUSF}}\label{Proposition: Lemma 2 point 6 of BKPS04}
		For every $L > 0$ there exists a constant $c > 0$ such that for every $z \in \Z^{d + 1},$ and every $\strW \subset \Z^{d + 1}$ with $\card{\strW} \leq L,$ we have $\spread{\strW \cup \{z\}} \leq \spread{\strW} \min\limits_{w \in \strW} \spread{wz} \leq c \spread{\strW \cup \{z\}}.$
	\end{Proposition}
	\begin{proof} The proof of this proposition proceeds word by word like that of Lemma 2.6 in \cite{BKPS:GeomUSF}. \end{proof}

	We may now generalise theorem (\ref{Theorem: Probability of two points belonging in the same component}) to any finite subset of vertices. We will only prove an upper bound.

	\begin{Theorem}\label{Theorem: The probability of finitely many vertices in the same component}
		Let $d \geq 3.$ For any $\strW \subset \Z^{d + 1},$ denote by $\scrC_\strW$ the event that all the points in $\strW$ are in the same $\USF$ component. For every integer $L > 0,$ there exists a constant $c > 0,$ depending only on $L$ and dimension, such that for whatever subset $\strW$ with at most $L$ vertices, $\USF(\scrC_\strW) \leq c \spread{\strW}^{-(d - 2)}.$
	\end{Theorem}
	\begin{proof}
		We apply the same reasoning as that used in the proof of Theorem 4.3 of \cite{BKPS:GeomUSF}. We will apply induction on the cardinality of $\strW.$ The case of two points is precisely theorem (\ref{Theorem: Probability of two points belonging in the same component}). Assume then the result holds for $\card{\strW} - 1 \geq 2.$ Let $z_1, z_2$ be two different vertices of $\strW$ and denote $\scrC_\strW(z_1, z_2)$ the event which is the intersection of $\scrC_\strW$ with the event where the $\USF$-path connecting $z_1$ and $z_2$ is edge-disjoint from the $\USF$-rays starting at the vertices of $\strV = \strW \setminus \{z_1, z_2\}.$ Let $\scrC_\strW^z(z_1, z_2)$ denote the event which is the intersection of $\scrC_\strW(z_1, z_2)$ and the event where the ray at $z_1$ meets the ray at $z_2$ at the vertex $z.$ Then
		\begin{displaymath}\textstyle
			\USF(\scrC_\strW(z_1, z_2)) \leq \sum\limits_{z \in \Z^{d+ 1}} \USF(\scrC_\strW^z(z_1, z_2)).
		\end{displaymath}
		To calculate the probability of $\USF(\scrC_\strW^z(z_1, z_2))$ we may assume that we construct the $\USF$-forest $\mathfrak{F}$ following Wilson's algorithm rooted at infinity by first ordering all the vertices in $\ds \strV z \mathop{=}^{\mathrm{def.}} \strV \cup \{z\},$ followed by the vertices $z_1$ and $z_2,$ and then the rest of the vertices of $\Z^{d + 1}.$ We see that
		\begin{displaymath}\textstyle
			\left\{ \mathfrak{F} \in \scrC_\strW^z(z_1, z_2) \right\} \subset \left\{ \mathfrak{F} \in \scrC_{\strV z} \right\} \cap \left\{ \tau_z \left( S^{z_1} \right) < \infty \right\} \cap \left\{ \tau_z \left( S^{z_2} \right) < \infty \right\},
		\end{displaymath}
		and $\left\{ \mathfrak{F} \in \scrC_{\strV z} \right\}$ is an event depending solely in the random walks $S^v$ for $v \in \strV z.$ By independence and induction,
		\begin{align*}
			\USF(\scrC_\strW^z(z_1, z_2)) &= \Probability{\mathfrak{F} \in \scrC_\strW^z(z_1, z_2)} \\
			&\leq \Probability{\mathfrak{F} \in \scrC_{\strV z}} \Probability{ \tau_z \left( S^{z_1} \right) < \infty } \Probability{ \tau_z \left( S^{z_2} \right) < \infty } \\
			&=\USF(\scrC_{\strV z}) G(z_1, z) G(z_2, z) \\
			&\leq c \spread{\strV z}^{-(d - 2)} G(z_1, z) G(z_2, z) \\
			&\textstyle\leq c' \spread{\strV}^{-(d - 2)} \sum\limits_{v \in \strV} \spread{vz}^{-(d - 2)} G(z_1, z) G(z_2, z),
		\end{align*}
		the last inequality by virtue of (\ref{Proposition: Lemma 2 point 6 of BKPS04}). Thus,
		\begin{equation}\textstyle
			\USF(\scrC_\strW(z_1, z_2)) \leq c \spread{\strV}^{-(d - 2)} \sum\limits_{v \in \strV} \sum\limits_{z \in \Z^{d + 1}} \spread{vz}^{-(d - 2)} G(z_1, z) G(z_2, z). \tag{$\ast$}
		\end{equation}
		We isolate some calculations first.
	
		\begin{EqInside}\label{EqInside: Sum of the spread inside a ball}
			For any $a \in \R$ and $r > 2,$ we have $\sum\limits_{\spread{z} \leq r} \spread{z}^{a} \asymp 1 + r^{d + 2 + a},$ the implicit constants depending only on $a$ and dimension.
		\end{EqInside}
		\begin{proof}[Proof of (\ref{EqInside: Sum of the spread inside a ball}).]
			Divide the cases $a < 0$ and $a \geq 0,$ use the fact that there exists an integer $k \geq 1$ such that $2^{k - 1} \leq r < 2^k$ and use that 
			\begin{displaymath}\textstyle
				\bigcup\limits_{j = 1}^{k - 1} \left\{ 2^{j - 1} \leq \spread{z} < 2^j \right\} \subset \left\{ \spread{z} \leq r \right\} \subset \bigcup\limits_{j = 1}^k \left\{ 2^{j - 1} \leq \spread{z} < 2^j \right\}.
			\end{displaymath}
			Keep in mind that $\spread{0} = 1.$
		\end{proof}
	
		\begin{EqInside}\label{EqInside: Weighted sum of Greens function by the spread}
			There exists a universal constant $c > 0$ such that for any two vertices $z_1, z_2 \in \Z^{d + 1},$ we have
			\begin{displaymath}\textstyle
				\sum\limits_{z \in \Z^{d + 1}} \spread{z}^{-(d - 2)} G(z_1, z) G(z_1, z) \leq c \spread{0z_1z_2}^{-(d - 2)}.
			\end{displaymath}
		\end{EqInside}
		\begin{proof}[Proof of (\ref{EqInside: Weighted sum of Greens function by the spread})]
			By symmetry, we may assume $\spread{z_1} \leq \spread{z_2}.$ Write $z_i = (n_i, x_i),$ for $i = 1, 2.$ Then,
			\begin{displaymath}\textstyle
				\sum\limits_{z \in \Z^{d + 1}} \spread{z}^{-(d - 2)} G(z_1, z) G(z_2, z) = \sum\limits_{\spread{z} > \frac{1}{2} \spread{z_1}} \spread{z}^{-(d - 2)} G(z_1, z) G(z_2, z) + I,
			\end{displaymath}
			and we may bound easily
			\begin{align*}
				\textstyle\sum\limits_{\spread{z} > \frac{1}{2} \spread{z_1}} \spread{z}^{-(d - 2)} G(z_1, z) G(z_2, z)
				&\textstyle\leq \spread{z_1}^{-(d - 2)} \sum\limits_{z \in \Z^{d + 1}} G(z_1, z) G(z_2, z) \\
				&\textstyle\leq c \spread{z_1}^{-(d - 2)} \spread{z_1 - z_2}^{-(d - 2)} \leq c \spread{0z_1z_2}^{-(d - 2)},
			\end{align*}
			where the penultimate inequality is obtained from (\ref{Lemma: Expectation of Kz is proportional to eta of z to the power of negative d minus two}) and the last inequality by definition of the spread. The rest of the calculation is to bound $I;$ this will be done in several cases. By definition
			\begin{displaymath}\textstyle
				I = \sum\limits_{\spread{z} \leq \frac{1}{2} \spread{z_1}} \spread{z}^{-(d - 2)} G(z_1, z) G(z_2, z).
			\end{displaymath}
			We will handle four cases, we will make use of the Green's function estimates (\ref{Theorem: Greens function bounds}). Notice that $G(z', z'') \leq c \norm{z' - z''}^{-\frac{d}{2}}$ in all cases and that $\japBracket{x} \leq \spread{z}$ for every $z \in \Z^{d + 1}$ (the left hand side is the spread relative to the Euclidean metric (\ref{Equation: Japanese bracket})). 
			\begin{description}
				\item[Case 1.] Here we assume that $\spread{z_1} = |n_1|^{\frac{1}{2}}$ and $\spread{z_2} = |n_2|^{\frac{1}{2}}.$ In this case, $\spread{z} \leq \frac{1}{2} \spread{z_1}$ implies that $|n| \leq \frac{1}{4} |n_1| \leq \frac{1}{4} |n_2|$ and $\norm{x - x_1} \leq \norm{x} + \norm{x_1} \leq \frac{3}{2} \norm{x_1} \leq \frac{3}{2} |n_1|^{\frac{1}{2}},$ and similarly $\norm{x - x_2} \leq \frac{3}{2} |n_2|^{\frac{1}{2}},$ obtaining $I \leq c |n_1|^{-\frac{d}{2}} |n_2|^{-\frac{d}{2}} \sum\limits_{\spread{z} \leq \frac{1}{2} \spread{z_1}} \spread{z}^{-(d - 2)}.$ We may apply (\ref{EqInside: Sum of the spread inside a ball}) to reach $	I \leq c \spread{z_1}^{-d} \spread{z_2}^{-d} \spread{z_1}^4 \leq c \spread{z_1}^{-(d - 2)} \spread{z_2}^{-(d - 2)} \leq c \spread{0z_1z_2}^{-(d - 2)}.$ This finishes \textbf{Case 1.}
		
				\item[Case 2.] Here we assume $\spread{z_1} = |n_1|^{\frac{1}{2}}$ and $\spread{z_2} = \norm{x_2}.$ Here we have $G(z_1, z) \leq c |n_1|^{-\frac{d}{2}}.$ We may split $I$ into the sum over $n$ and over $x.$ This gives $I \leq c |n_1|^{-\frac{d}{2}} \sum\limits_{\norm{x} \leq \frac{1}{2} |n_1|^{\frac{1}{2}}} \japBracket{x}^{-(d - 2)} \sum\limits_{|n| \leq \frac{|n_1|}{4}} G(z_2, z).$ For every $\norm{x} \leq \frac{1}{2} |n_1|^{\frac{1}{2}},$ we have $\norm{x - x_2} \leq \norm{x} + \norm{x_2} \leq \frac{3}{2} \norm{x_2}$ and $\norm{x - x_2} \geq \norm{x_2} - \norm{x} \geq \frac{1}{2} \norm{x_2}.$ Next, break the inner sum into two parts, the first one being over all $n$ such that $|n - n_2| \leq \frac{1}{2} \norm{x_2}$ and the second one being over all other $n.$ In the first part, we have a bound of the form $G(z_2, z) \leq c e^{-c' \norm{x_2}},$ thus $\sum\limits_{\substack{|n| \leq \frac{|n_1|}{4} \\ |n - n_2| \leq \frac{\norm{x_2}}{2}}} G(z_2, z) \leq c \norm{x_2}^2 e^{-c' \norm{x_2}} \leq c'' \norm{x_2}^{-(d - 2)}.$ We focus now in the second part of this inner sum. Here, we have $\frac{1}{2} \norm{x_2} \leq |n - n_2| \leq |n| + |n_2| \leq \frac{5}{4} \norm{x_2}^2.$ Thus, we can bound $G(z_2, z) \leq c e^{-c' \frac{\norm{x_2}^2}{k}} k^{-\frac{d}{2}}$ for some $1 \leq k \leq \frac{5}{4} \norm{x_2}^2.$ Then, we may bound the second part of the inner sum of $I$ as follows $\sum\limits_{\substack{|n| \leq \frac{|n_1|}{4} \\ |n - n_2| > \frac{\norm{x_2}}{2}}} G(z_2, z) \leq c \sum\limits_{1 \leq k \leq \frac{5}{4} \norm{x_2}^2} e^{-c' \frac{\norm{x_2}^2}{k}} k^{-\frac{d}{2}}.$ By virtue of (\ref{Proposition: Estimates in the p series with an exponential factor}), we have $\sum\limits_{1 \leq k \leq \frac{5}{4} \norm{x_2}^2} e^{-c' \frac{\norm{x_2}^2}{k}} k^{-\frac{d}{2}} \leq c \norm{x_2}^{2(-\frac{d}{2} + 1)} = c \norm{x_2}^{-(d - 2)}.$ Therefore, we may bound $I$ as follows
				\begin{displaymath}
					I \leq c |n_1|^{-\frac{d}{2}} \sum\limits_{\norm{x} \leq \frac{1}{2} |n_1|^{\frac{1}{2}}} \japBracket{x}^{-(d - 2)} \norm{x_2}^{-(d - 2)}.
				\end{displaymath}
				We have $\sum\limits_{\norm{x} \leq \frac{1}{2} |n_1|^{\frac{1}{2}}} \japBracket{x}^{-(d - 2)} \leq 1 + c \integral[1]^{\frac{1}{2} |n_1|^{\frac{1}{2}}} dt\ t^{-(d - 2)} t^{d - 1} = 1 + c \integral[1]^{\frac{1}{2} |n_1|^{\frac{1}{2}}} dt\ t \leq c' |n_1|.$ Substituting above yields $I \leq c |n_1|^{-\frac{d}{2} + 1} \norm{x_2}^{-(d - 2)} = c \spread{z_1}^{-(d - 2)} \spread{z_2}^{-(d - 2)} \leq c \spread{0z_1 z_2}^{-(d - 2)}.$
				This completed \textbf{Case 2.}
		
				\item[Case 3.] Here we assume $\spread{z_1} = \norm{x_1}$ and $\spread{z_2} = |n_2|^{\frac{1}{2}}.$ Then $|n - n_2| \asymp |n_2|$ since $|n| \leq \frac{\norm{x_1}}{2} \leq \frac{|n_2|}{2}$ and we may bound $G(z, z_2) \leq c \norm{z - z_2}^{-\frac{d}{2}} \leq c' |n_2|^{-\frac{d}{2}}.$ This gives $I \leq c |n_2|^{-\frac{d}{2}} \sum\limits_{\spread{z} \leq \frac{1}{2} \spread{z_1}} \spread{z}^{-(d - 2)} G(z_1, z).$ We may now employ the same split of this sum as in \textbf{Case 2}, deducing
				\begin{displaymath}\textstyle
					I \leq c |n_2|^{-\frac{d}{2}} \sum\limits_{\norm{x} \leq \frac{1}{2} \norm{x_1}} \japBracket{x}^{-(d -2)} \sum\limits_{|n| \leq \frac{1}{4} \norm{x_1}^2} G(z_1, z).
				\end{displaymath}
				We may break again the inner sum over those $n$ for which $|n - n_1| \leq \frac{1}{2} \norm{x_1}$ and over all other $n,$ and we will find that
				\begin{displaymath}\textstyle
					\sum\limits_{|n| \leq \frac{1}{2} \norm{x_1}} G(z_1, z) \leq c \left[ \sum\limits_{\substack{|n| \leq \frac{1}{2} \norm{x_1} \\ |n - n_1| \leq \frac{1}{2} \norm{x_1}}} e^{-c'\norm{x_1}} + \sum\limits_{1 \leq k \leq \frac{5}{4} \norm{x_1}^2} e^{-c' \frac{\norm{x_1}^2}{k}} k^{-\frac{d}{2}} \right].
				\end{displaymath}
				And again, we apply (\ref{Proposition: Asymptotic estimates on Eulerian sums}) and (\ref{Proposition: Estimates in the p series with an exponential factor}), obtaining
				\begin{align*}
					\textstyle\sum\limits_{\norm{x} \leq \frac{1}{2} \norm{x_1}} \japBracket{x}^{-(d -2)} &\leq c \norm{x_1}^2 \\
					\textstyle\sum\limits_{1 \leq k \leq \frac{5}{4} \norm{x_1}^2} e^{-c' \frac{\norm{x_1}^2}{k}} k^{-\frac{d}{2}} &\leq c \norm{x_1}^{2(-\frac{d}{2} + 1)} = c \norm{x_1}^{-(d - 2)}.
				\end{align*}
				It is also clear that $\norm{x_1} e^{-c' \norm{x_1}} \leq c \norm{x_1}^{-(d - 2)}.$ Substituting yields $I \leq c |n_2|^{-\frac{d}{2}} \norm{x_1}^2 \norm{x_1}^{-(d - 2)}.$ Since we assume $\spread{z_1} \leq \spread{z_2},$ we have $\norm{x_1}^2 \leq |n_2|,$ so that $	|n_2|^{-\frac{d}{2}} \norm{x_1}^2 \leq |n_2|^{-\frac{d - 2}{2}} = \spread{z_2}^{-(d - 2)}.$ The proof of \textbf{Case 3} is finished since $\spread{z_1} = \norm{x_1}$ and the definition of spread shows at once $\spread{z_1}^{-(d - 2)} \spread{z_2}^{-(d - 2)} \leq \spread{0z_1z_2}^{-(d - 2)}.$
		
				\item[Case 4.] Here we assume $\spread{z_1} = \norm{x_1}$ and $\spread{z_2} = \norm{x_2}.$
				Consider first the set of $z \in \Z^{d + 1}$ such that $\spread{z} \leq \frac{1}{2} \spread{z_1}$ and $|n - n_1| \leq \norm{x_1}.$ For $z$ in this set, $G(z_1, z) \leq c e^{-c' \norm{x_1}},$ this entails
				\begin{displaymath}\textstyle
					\sum\limits_{\substack{\spread{z} \leq \frac{1}{2} \norm{x_1} \\ |n - n_1| \leq \norm{x_1}}} \spread{z}^{-(d - 2)} G(z_1, z) G(z_2, z) \leq c e^{-c' \norm{x_1}} \sum\limits_{\norm{x} \leq \frac{1}{2} \norm{x_1}} \japBracket{x}^{-(d - 2)} \sum\limits_{|n| \leq \frac{1}{4} \norm{x_1}^2 } G(z_2, z)
				\end{displaymath}
				Recall we assume $\norm{x_1} \leq \norm{x_2},$ and, as done in \textbf{Case 2} and \textbf{Case 3} $\sum\limits_{|n| \leq \frac{1}{4} \norm{x_2}^2} G(z_2, z) \leq c \norm{x_2}^{-(d - 2)}.$ Therefore,
				\begin{align*}
					\textstyle\sum\limits_{\substack{\spread{z} \leq \frac{1}{2} \norm{z_1} \\ |n - n_1| \leq \norm{x_1}}} \spread{z}^{-(d - 2)} G(z_1, z) G(z_2, z) &\leq c \norm{x_1}^2 e^{-c' \norm{x_1}} \norm{x_2}^{-(d - 2)} \\
					&\leq c'' \norm{x_1}^{-(d - 2)} \norm{x_2}^{-(d - 2)} \leq c'' \spread{0z_1z_2}.
				\end{align*}
				There remains to handle the terms defining $I$ for which $\spread{z} \leq \frac{1}{2} \norm{x_1}$ and $|n - n_1| > \norm{x_1}.$ Let $J$ be the least integer for which $\frac{1}{2} \norm{x_1} \leq 2^J \norm{x_1}^{\frac{1}{2}}.$ In other words, $J$ is the least integer for which $\norm{x_1} \leq 4^{J + 1}.$ In particular, $4^J < \norm{x_1} \leq 4^{J + 1}.$ Next, observe that the relation $\spread{z} \leq \frac{1}{2} \norm{x_1}$ implies $|n - n_1| \leq |n| + |n_1| \leq \frac{5}{4} \norm{x_1}^2 \leq 5 \cdot 4^J \norm{x_1}.$ Define
				\begin{displaymath}\textstyle
					\begin{split}
						\strA_j &\textstyle= \left\{ z;\ 4^{j - 1} \norm{x_1} <  |n - n_1| \leq 4^j \norm{x_1}, \spread{z} \leq \frac{1}{2} \norm{x_1} \right\}, \quad 1 \leq j \leq J - 1 \\
						\strA_J &\textstyle= \left\{ z;\ 4^{J - 1} \norm{x_1} <  |n - n_1| \leq 5 \cdot 4^J \norm{x_1}, \spread{z} \leq \frac{1}{2} \norm{x_1} \right\}.
					\end{split}
				\end{displaymath}
				We can write $I \leq c \spread{0z_1z_2}^{-(d - 2)} + \sum\limits_{j = 1}^J \sum\limits_{z \in \strA_j} \spread{z}^{-(d - 2)} G(z_1, z) G(z_2, z).$ If $z \in \strA_j,$ then $G(z_1, z) \leq c e^{-c' \frac{\norm{x_1}}{4^j}} (4^j \norm{x_1})^{-\frac{d}{2}},$ and since $\norm{x_1} \asymp 4^J,$ $G(z_1, z) \leq c \norm{x_1}^{-\frac{d}{2}} e^{-c' 4^{J - j}} 4^{-j \frac{d}{2}}.$ Then,
				\begin{displaymath}\textstyle
					\sum\limits_{j = 1}^J \sum\limits_{z \in \strA_j} \spread{z}^{-(d - 2)} G(z_1, z) G(z_2, z) \leq c \norm{x_1}^{-\frac{d}{2}} \sum\limits_{j = 1}^J e^{-c' 4^{J - j}} 4^{-j \frac{d}{2}} \sum\limits_{z \in \strA_j} \spread{z}^{-(d - 2)} G(z_2, z).
				\end{displaymath}
				We may break the inner sum in terms of $n$ and $x$ (for details, see \textbf{Case 2} and \textbf{Case 3})  and this will give $\sum\limits_{z \in \strA_j} \spread{z}^{-(d - 2)} G(z_2, z) \leq \sum\limits_{\norm{x} \leq \frac{1}{2} \norm{x_1}} \japBracket{x}^{-(d - 2)} \sum\limits_{|n| \leq \frac{1}{4} \norm{x_1}^2} G(z_2, z) \leq c \norm{x_1}^2 \norm{x_2}^{-(d - 2)}.$ Whence,
				\begin{displaymath}\textstyle
					\sum\limits_{j = 1}^J \sum\limits_{z \in \strA_j} \spread{z}^{-(d - 2)} G(z_1, z) G(z_2, z) \leq c \norm{x_1}^{-\frac{d}{2}} \norm{x_1}^2 \norm{x_2}^{-(d - 2)} \sum\limits_{j = 1}^J e^{-c' 4^{J - j}} 4^{-j \frac{d}{2}}.
				\end{displaymath}
				Finally, in the inner sum reverse the order of summation, $\sum\limits_{j = 1}^J e^{-c' 4^{J - j}} 4^{-j \frac{d}{2}} = \sum\limits_{k = 0}^{J - 1} e^{-c'4^k} 4^{(k - J) \frac{d}{2}} \leq 4^{-J \frac{d}{2}} \sum\limits_{k = 0}^\infty 4^{k \frac{d}{2}} e^{-c' 4^k} \asymp \norm{x_1}^{-\frac{d}{2}}.$ Thus,
				\begin{align*}
					\textstyle\sum\limits_{j = 1}^J \sum\limits_{z \in \strA_j} \spread{z}^{-(d - 2)} G(z_1, z) G(z_2, z) &\leq c \norm{x_1}^{-\frac{d}{2}} \norm{x_1}^2 \norm{x_2}^{-(d - 2)} \norm{x_1}^{-\frac{d}{2}} \\
					&= c \norm{x_1}^{-(d - 2)} \norm{x_2}^{-(d - 2)} \\
					&= c \spread{z_1}^{-(d - 2)} \spread{z_2}^{-(d - 2)} \leq c \spread{0z_1z_2}^{-(d - 2)}.
				\end{align*}
				This completes \textbf{Case 4.}
			\end{description}
			With the four cases established, we have reached the conclusion of (\ref{EqInside: Weighted sum of Greens function by the spread}).
		\end{proof}

		A corollary of (\ref{EqInside: Weighted sum of Greens function by the spread}) is that for any vertex $v \in \Z^{d + 1},$
		\begin{displaymath}\textstyle
			\sum\limits_{z \in \Z^{d + 1}} \spread{vz}^{-(d - 2)} G(z_1, z) G(z_2, z) \leq c \spread{\{0, z_1 - v, z_2 - v\}}^{-(d - 2)} = c \spread{vz_1z_2}^{-(d - 2)}.
		\end{displaymath}
		Substituting this inequality in ($\ast$) allows deducing, with the aid of (\ref{Proposition: Lemma 2 point 6 of BKPS04}),
		\begin{displaymath}\textstyle
			\USF(\scrC_\strW(z_1, z_2)) \leq c \spread{\strV}^{-(d - 2)} \sum\limits_{v \in \strV} \spread{vz_1z_2}^{-(d - 2)} \leq c' \spread{\strW}^{-(d - 2)},
		\end{displaymath}
		where the last inequality follows form the fact that the union of a tree on $\{v, z_1, z_2\}$ and a tree on $\strV$ gives a tree on $\strW = \strV \cup \{z_1, z_2\}.$ We will show now that
		\begin{displaymath}\textstyle
			\scrC_\strW = \bigcup\limits_{z_1 \neq  z_2} \scrC_\strW(z_1, z_2),
		\end{displaymath}
		where the vertices $z_1$ and $z_2$ run on $\strW.$ For each vertex $w \in \strW,$ denote by $R^w$ the $\USF$-ray from $w.$ For every $z_1 \neq z_2,$ two vertices of $\strW,$ their corresponding rays $R^{z_1}$ and $R^{z_2}$ will meet at a vertex $v(z_1, z_2).$ Furthermore, on the event $\scrC_\strW,$ the set of vertices belonging to all $R^w,$ for $w \in \strW,$ is non empty and thus, there exists a vertex $z^*$ satisfying the following: $z^*$ is the closest vertex in the graph-distance of the $\USF$-tree, amongst all vertices belonging to all rays $R^w$ ($w \in \strW$), to each one of the vertices $w \in \strW.$ In other words, $z^*$ is the vertex where all the rays meet. Consider now a pair $z_1 \neq z_2$ such that the graph-distance of $v(z_1, z_2)$ and $z^*$ is maximal. Then, $\scrC_\strW(z_1, z_2)$ occurs. The proof of theorem (\ref{Theorem: The probability of finitely many vertices in the same component}) is complete.
	\end{proof}

	\section[Separation between components]{The separation between components: transitions in dimensions $d = 2, 6, 10,...$}\label{Section: Separation between components}
Given a vertex $z$ of $\Gamma_d(\lambda),$ denote by $\mathfrak{T}_z$ the component of $z$ in the $\USF$-forest, regarded as spanning subgraph. Construct now $\Gamma_d(\lambda)/\USF$ from $\Gamma_d(\lambda)$ by shorting all edges of all $\mathfrak{T}_z$ as $z$ runs on the vertices of $\Gamma_d(\lambda).$ Therefore, $\mathfrak{T}_z$ is a vertex of the new network $\Gamma_d(\lambda)/\USF.$ Define $D:\Z^d \to \Z_+$ by
\begin{displaymath}\textstyle
	D(z, z') = d_{\Gamma_d(\lambda)/\USF}\left( \mathfrak{T}_z, \mathfrak{T}_{z'} \right),
\end{displaymath}
the right hand side being the graph-distance of $\Gamma_d(\lambda)/\USF.$ We call this the \textbf{separation} between the components of $z$ and $z'$ relative to the chemical distance of the graph $\Z^{d + 1}.$ Another way of describing $D(z, z')$ is as follow: it is the minimal number of edges of $\Gamma_d(\lambda)$ that do not belong to any $\USF$-tree in a path connecting $z$ and $z'.$ In this section we are going to show that, for almost every realisation,
\begin{displaymath}\textstyle
	\max\limits_{z, z' \in \Z^{d + 1}} D(z, z') = \left\lceil \frac{d - 2}{4} \right\rceil.
\end{displaymath}
This section is substantially based on \cite{BKPS:GeomUSF}.

	\subsection{Random relations and stochastic dimension}
	We will be considering relations between vertices of $\Gamma_d(\lambda).$ If $R$ is one such relation, we will use the common notation $z_1Rz_2$ to mean $(z_1, z_2) \in R.$ Also, the \textbf{composition} of two relations $L$ and $R$ is, by definition, the set of pairs $(z_1, z_2)$ such that there exists a $z$ satisfying $z_1Lz$ \emph{and} $zRz_2.$ We write $LR$ to denote the composition of $L$ and $R$ (in this order). We will say that a random relation $R$ on the vertices of $\Gamma_d(\lambda)$ has \textbf{stochastic dimension} $\alpha \in [0, d + 1),$ relative to the metric $\eta$ (\ref{Equation: Metric in the probability of same component}), and write $\dim_\eta(R) = \alpha$ if there exists a constant $c = c(R) > 0$ such that for all vertices $z_1, z_2, z_3, z_4$ of $\Gamma_d(\lambda)$
	\begin{displaymath}\textstyle
		c\Probability{z_1Rz_2} \geq \spread{z_1z_2}^{-(d + 1 - \alpha)},
	\end{displaymath}
	and
	\begin{displaymath}\textstyle
		\Probability{z_1Rz_2, z_3Rz_4} \leq c \spread{z_1z_2}^{-(d + 1 - \alpha)} \spread{z_3 z_4}^{-(d + 1 - \alpha)} + c \spread{z_1z_2z_3z_4}^{-(d + 1 - \alpha)}.
	\end{displaymath}

	\begin{Remark}
		That a random relation $R$ has stochastic dimension $d + 1$ is the same as saying that $\ds \inf_{z_1, z_2} \Probability{z_1Rz_2} > 0.$
	\end{Remark}

	\begin{Remark}\label{Remark: The dimension of the USF relation}
		Let $d \geq 3.$ Denote by $U_{\Gamma_d(\lambda)}$ the random relation: ``$z_1$ and $z_2$ are in the same $\USF$-component.'' Then, $U_{\Gamma_d(\lambda)}$ has stochastic dimension \emph{three}, relative to $\eta,$ this is a consequence of the definition and of theorems (\ref{Theorem: Probability of two points belonging in the same component}) and (\ref{Theorem: The probability of finitely many vertices in the same component}). Indeed, write $z \sim z'$ to mean $z U_{\Gamma_d(\lambda)} z'.$ Then, (\ref{Theorem: Probability of two points belonging in the same component}) shows that $c \Probability{z_1 \sim z_2} \geq \spread{z_1 z_2}^{-(d + 1 - 3)}.$ Next, observe that
		\begin{displaymath}\textstyle
			\Probability{z_1 \sim z_2, z_3 \sim z_4} \leq \Probability{z_1 \sim z_2, z_3 \sim z_4, z_1 \sim z_3} + \Probability{z_1 \sim z_2, z_3 \sim z_4, z_1 \not \sim z_3}.
		\end{displaymath}
		The event $\{z_1 \sim z_2, z_3 \sim z_4, z_1 \sim z_3\}$ is the event where all the four vertices are in the same $\USF$-tree, this the probability of this event is bounded above by a universal multiple of $\spread{z_1z_2z_3z_4}^{-(d+1-3)}$ by (\ref{Theorem: The probability of finitely many vertices in the same component}). To estimate the probability of the event $\{z_1 \sim z_2, z_3 \sim z_4, z_1 \not \sim z_3\}$ we may construct the $\USF$-forest $\mathfrak{F}$ using Wilson's algorithm rooted at infinity with ordering of vertices $(z_1, z_2, z_3, z_4, \ldots).$ If $S^z$ denotes the random walk in this construction of $\mathfrak{F},$ then
		\begin{align*}
			\Probability{z_1 \sim z_2, z_3 \sim z_4, z_1 \not \sim z_3}& \leq \Probability{S^{z_1} \cap S^{z_2} \neq \varnothing, S^{z_3} \cap S^{z_4} \neq \varnothing} \\
			&\leq \Probability{S^{z_1} \cap S^{z_2} \neq \varnothing} \Probability{S^{z_3} \cap S^{z_4} \neq \varnothing} \\
			&\leq c \spread{z_1z_2}^{-(d + 1 - 3)} \spread{z_3 z_4}^{-(d + 1 - 3)},
		\end{align*}
		the last inequality by virtue of lemmas (\ref{Lemma: The probability of Kz is positive is asymptotic to its expectation}) and (\ref{Lemma: Expectation of Kz is proportional to eta of z to the power of negative d minus two}).
	\end{Remark}

	\begin{Remark}
		We are going to base the results of this section on those of the paper \cite{BKPS:GeomUSF}, making the necessary adaptations according to our conveniences. We will see that, in a sense, the stochastic dimension relative to $\eta$ is the same the stochastic dimension they defined (that is, relative to the Euclidean distance) \emph{minus} one.
	\end{Remark}

	We introduce the following annuli: for $n < N$ two integers and $z \in \Z^{d + 1},$
	\begin{displaymath}\textstyle
		\strA_n^N(z) = \{z' \in \Z^{d + 1} \mid 2^n \leq \spread{z - z'} < 2^N\}.
	\end{displaymath}
	Notice $\card{\strA_n^N(z)} \asymp 2^{N(d + 2)}.$

	\begin{Proposition}{\normalfont(\cite[Lemma 2.8]{BKPS:GeomUSF})}\label{Proposition: Lemma 2 point 8 of BKPS}
		Let $L$ and $R$ be two independent random relations of $\Gamma_d(\lambda).$ Suppose that $\dim_\eta(L) = d + 1 - \alpha$ and $\dim_\eta(R) = d + 1 - \beta,$ both exist. Denote  $\gamma = \alpha + \beta - d - 2.$ For $z_1, z_2 \in \Z^{d + 1},$ and $1 \leq n \leq N,$ let $S_{z_1, z_2} = S_{z_1, z_2}(n, N) \mathop{=}\limits^{\mathrm{def.}} \sum\limits_{z \in A_n^N(z_1)} \indic{\{z_1Lz\}} \indic{\{zRz_2\}}.$ If $\spread{z_1z_2} < 2^{n - 1}$ and $N \geq n,$ $\Probability{S_{z_1, z_2} > 0} \geq c \frac{\sum\limits_{k = n}^N 2^{-k \gamma}}{\sum\limits_{k = 0}^N 2^{-k \gamma}},$ where $c$ is a constant that may depend solely on $L,$ $R$ and dimension.
	\end{Proposition}
	\begin{proof}
		We may repeat word by word, with very minor notational modifications, the proof of Lemma 2.8 of \cite{BKPS:GeomUSF} bearing in mind that the only substantial thing that changes is $\card{\strA_k^{k + 1}(z_1)} \asymp 2^{k(d + 2)}.$ Also keep in mind that their Lemma 2.6 is our Proposition (\ref{Proposition: Lemma 2 point 6 of BKPS04}). Thus, their formula (2.6) changes to
		\begin{displaymath}\textstyle
			\Expectation{S_{z_1, z_2}} \asymp \sum\limits_{k = n}^N 2^{k(d + 2)} 2^{-k \alpha} 2^{-k \beta} = \sum\limits_{k = n}^N 2^{-k \gamma},
		\end{displaymath}
		the rest of the proof goes on without major modification except that we need to replace $d$ by $d + 2$ in their formula (2.7).
	\end{proof}

	\begin{Proposition}{\normalfont(\cite[Corollary 2.9]{BKPS:GeomUSF})}\label{Proposition: Corollary 2 point 9 of BKPS04}
		Let $L$ and $R$ be two independent random relations of $\Gamma_d(\lambda).$ Suppose that $\dim_\eta(L) = d + 1 - \alpha$ and $\dim_\eta(R) = d + 1 - \beta,$ both exist. Denote  $\gamma = \max(0, \alpha + \beta - d - 2).$ There exists a constant $c > 0$ such that for all $z_1, z_2 \in \Z^{d + 1},$
		\begin{displaymath}\textstyle
			\Probability{z_1 LR z_2} \geq c \spread{z_1 z_2}^{-\gamma}.
		\end{displaymath}
	\end{Proposition}
	\begin{proof}
		The proof of Corollary 2.9 of \cite{BKPS:GeomUSF} applies almost verbatim: in their second sentence replace ``lemma'' with ``Proposition (\ref{Proposition: Lemma 2 point 8 of BKPS}).''
	\end{proof}

	\begin{Proposition}{\normalfont(\cite[Lemma 2.10]{BKPS:GeomUSF})}\label{Proposition: Lemma 2 point 10 of BKPS04}
		Let $\alpha, \beta \in [0, d + 1)$ satisfy $\alpha + \beta > d + 2.$ Let $\gamma = \alpha + \beta - d - 2.$ Then,
		\begin{displaymath}\textstyle
			\sum\limits_{z \in \Z^{d + 1}} \spread{z_1z}^{-\alpha} \spread{z_2z}^{-\beta} \asymp \spread{z_1z_2}^{-\gamma}
		\end{displaymath}
		for all vertices $z_1, z_2$ of $\Z^{d + 1},$ and with any implicit constant being universal.
	\end{Proposition}
	\begin{proof}
		The proof of Lemma 2.10 of \cite{BKPS:GeomUSF} applies without changes except to notation.
	\end{proof}

	\begin{Proposition}{\normalfont(\cite[Lemma 2.11]{BKPS:GeomUSF})}\label{Proposition: Lemma 2 point 11 of BKPS04}
		Let $M$ be a positive integer and let $\alpha, \beta \in [0, d + 1)$ satisfy $\alpha + \beta > d + 2.$  Denote $\gamma = \alpha + \beta - d - 2.$ There exists a constant $c > 0$ such that for all subsets $\strV$ and $\strW$ of $\Z^{d + 1}$ with at most $M$ points, we have
		\begin{displaymath}\textstyle
			\sum\limits_{z \in \Z^{d + 1}} \spread{\strV \cup \{z\}}^{-\alpha} \spread{\strW \cup \{z\}}^{-\beta} \leq c \spread{\strV}^{-\alpha} \left( \min\limits_{(v, w) \in \strV \times \strW} \spread{vw} \right)^{-\gamma} \spread{\strW}^{-\beta} \leq c \spread{\strV \cup \strW}^{-\gamma}.
		\end{displaymath}
	\end{Proposition}
	\begin{proof}
		We may use the same proof as that of Lemma 2.11 of \cite{BKPS:GeomUSF} keeping in mind that their Lemma 2.6 is our Proposition (\ref{Proposition: Lemma 2 point 6 of BKPS04}) and, Lemma 2.10 is Proposition (\ref{Proposition: Lemma 2 point 10 of BKPS04}).
	\end{proof}

	\begin{Proposition}{\normalfont(\cite[Lemma 2.12]{BKPS:GeomUSF})}\label{Proposition: Lemma 2 point 12 of BKPS04}
		Let $\alpha, \beta \in [0, d + 1)$ satisfy $\alpha + \beta > d + 2.$ Set $\gamma = \alpha + \beta - d - 2.$
		\begin{enumerate}
			\item There exists a constant $c > 0$ such that for all vertices $z_1, z_2, z_3 \in \Z^{d + 1},$
			\begin{displaymath}\textstyle
				\sum\limits_{z \in \Z^{d + 1}} \spread{z_1z}^{-\alpha} \spread{z_2z}^{-\beta} \spread{z_3z}^{-\gamma} \leq c \spread{z_1z_2z_3}^{-\gamma}
			\end{displaymath}

			\item There exists a constant $c > 0$ such that for all vertices $z_1, z_2, z_3, z_4 \in \Z^{d + 1},$
			\begin{displaymath}\textstyle
				\sum\limits_{z \in \Z^{d + 1}} \spread{z_1z_2z}^{-\alpha} \spread{z_3z}^{-\beta} \spread{z_4z}^{-\gamma} \leq c \spread{z_1z_2z_3z_4}^{-\gamma}
			\end{displaymath}
		\end{enumerate}
		
	\end{Proposition} 
	\begin{proof}
		 The proof is essentially the same as that of Lemma 2.12 of \cite{BKPS:GeomUSF} and their formula (2.9), except we need to use (\ref{EqInside: Sum of the spread inside a ball}) in some step. For convenience of the reader, we repeat the proof.

		We begin by establishing the first item. By changing the indices should the need arise, we may assume that $\spread{z_1z_3} \leq \spread{z_2z_3}.$ Denote by $\strA$ the set of $z \in \Z^{d + 1}$ satisfying $\spread{z_3z} \leq \frac{1}{2} \spread{z_1z_3}.$ When $z \in \strA$ we have, $\spread{z_1z} \geq \spread{z_1z_3} - \spread{z_3z} \geq \frac{1}{2} \spread{z_1z_3}$ and, similarly, since $\spread{z_1z_3} \leq \spread{z_2z_3},$ we also have $\spread{z_2z} \geq \frac{1}{2} \spread{z_2z_3}.$ Therefore,
		\begin{displaymath}\textstyle
			\sum\limits_{z \in \strA} \spread{z_1z}^{-\alpha} \spread{z_2z}^{-\beta} \spread{z_3z}^{-\gamma} \leq 2^\alpha \spread{z_1z_3}^{-\alpha} 2^\beta \spread{z_2z_3}^{-\beta} \sum\limits_{z \in \strA} \spread{z_3z}^{-\gamma}.
		\end{displaymath}
		Apply (\ref{EqInside: Sum of the spread inside a ball}) to obtain the bound $\sum\limits_{z \in \strA} \spread{z_3z}^{-\gamma} \leq c \spread{z_1 z_3}^{d + 2 - \gamma}.$ Substitute above,
		\begin{align*}
			\textstyle\sum\limits_{z \in \strA} \spread{z_1z}^{-\alpha} \spread{z_2z}^{-\beta} \spread{z_3z}^{-\gamma}  &\leq c \spread{z_1z_3}^{-\alpha} \spread{z_2z_3}^{-\beta} \spread{z_1 z_3}^{d + 2 - \gamma} \\
			&= c \spread{z_1z_3}^{-\gamma} \spread{z_2z_3}^{-\gamma} \frac{\spread{z_1 z_3}^{d + 2 - \alpha}}{\spread{z_2z_3}^{d + 2 - \alpha}} \leq c \spread{z_1z_2z_3}^{-\gamma},
		\end{align*}
		because $\gamma > 0$ and $\spread{z_1z_3} \leq \spread{z_2z_3}.$ We now estimate the sum over the complement of $\strA.$ Here we apply (\ref{Proposition: Lemma 2 point 10 of BKPS04}), by definition of $\strA,$ the relation $z \notin \strA$ implies $\spread{z_3z}^{-\gamma} \leq 2^\gamma \spread{z_1z_3}^{-\gamma}.$ Thus,
		\begin{align*}
			\textstyle\sum\limits_{z \notin \strA} \spread{z_1z}^{-\alpha} \spread{z_2z}^{-\beta} \spread{z_3z}^{-\gamma}
			&\textstyle\leq 2^\gamma \spread{z_1z_3}^{-\gamma} \sum\limits_{z \in \Z^{d + 1}} \spread{z_1z}^{-\alpha} \spread{z_2z}^{-\beta} \\
			&\leq c \spread{z_1z_3}^{-\gamma} \spread{z_1z_2}^{-\gamma} \leq c \spread{z_1z_2z_3}^{-\gamma},
		\end{align*}
		by virtue of (\ref{Proposition: Lemma 2 point 10 of BKPS04}). Combining these two estimates we obtain the first item of the proposition.
	
		Now we prove the second item. By (\ref{Proposition: Lemma 2 point 6 of BKPS04}), we obtain $\spread{z_1z_2z} \geq c \spread{z_1z_2} \min\{\spread{z_1z}, \spread{z_2z}\}.$ Use this, the first item of the proposition, and recall that $\spread{z_1z_2}\spread{z_2z_3z_4} \geq \spread{z_1z_2z_3z_4}$ which holds by definition of the spread.
	\end{proof}	

	\begin{Theorem}{\normalfont(\cite[Theorem 2.4]{BKPS:GeomUSF})}
		Let $L$ and $R$ be two independent random relations on the vertices of $\Gamma_d(\lambda).$ Assume that their stochastic dimensions (relative to $\eta$) exist. Then, their composition $LR$ also has stochastic dimension relative to $\eta$ and it satisfies $\dim_\eta(LR) = \min \left( \dim_\eta(L) + \dim_\eta(R) + 1, d + 1 \right).$
	\end{Theorem}
	\begin{proof}
		Set $\alpha$ and $\beta$ to be defined as follows $\dim_\eta(L) = d + 1 - \alpha$ and $\dim_\eta(R) = d + 1 - \beta.$ Notice that $\gamma \mathop{=}\limits^{\mathrm{def.}} \alpha + \beta - d - 2 = d + 1 - \big( \dim_\eta(L) + \dim_\eta(R) + 1 \big).$ The arguments in the proof of Theorem 2.4 of \cite{BKPS:GeomUSF} can be modified very easily by recalling that Corollary 2.9 in their paper is (\ref{Proposition: Corollary 2 point 9 of BKPS04}), Lemma 2.10 is (\ref{Proposition: Lemma 2 point 10 of BKPS04}), Lemma 2.11 is (\ref{Proposition: Lemma 2 point 11 of BKPS04}) and equation (2.9) is part (b) of (\ref{Proposition: Lemma 2 point 12 of BKPS04}). Also, notice that their condition ``$\dim_S(\mathcal{L}) + \dim_S(\mathcal{R}) > d$'' changes to ``$\dim_\eta(L) + \dim_\eta(R) + 1 > d + 1$'', the latter being equivalent to $\gamma > 0,$ and similarly when the symbol $>$ is replaced by $\leq.$
	\end{proof}

	\subsection{Tail triviality of random relations}
	Recall that an event in any probability space is said to \textbf{trivial} if its probability is either zero or one, equivalently, if its probability $p$ satisfies $p = p^2.$ It is easy to see that the set of trivial events is a $\sigma$-field, called the \textbf{trivial} $\sigma$-field (relative to the given probability measure). A $\sigma$-algebra is said to be \textbf{trivial} if it is contained in the trivial $\sigma$-field. If $R$ is a random relation on the vertices of $\Gamma_d(\lambda),$ we denote by $\scrF^R_\strA,$ where $\strA \subset \Z^{d + 1} \times \Z^{d + 1},$ the $\sigma$-algebra generated by the events $\{zRz'\}$ as $(z, z') \in \strA.$ If $z$ is a vertex, we define the \textbf{tail $\sigma$-field based at $z$ on the left of $R,$} denoted by $\scrL_z^R,$ to be the $\sigma$-field generated by the events $\{zRz'\}$ as $z' \to \infty,$ in other words, $\scrL_z^R$ is the intersection of all $\scrF_{\{z\} \times \strK^\complement}^R$ as $\strK$ runs on all finite subsets of $\Z^{d + 1}.$ For simplicity, we will say ``left tail at $z$'' to refer to $\scrL_z^R.$ In a similar manner we define the \textbf{tail $\sigma$-field $\scrR_z^R$ based at $z$ on the right of $R,$} to be the intersection of all $\scrF_{\strK^\complement \times \{z\}}^R$ as $\strK$ runs on all finite subsets of $\Z^{d + 1};$ we will talk about ``right tail at $z$'' to refer to this $\sigma$-algebra. We also define the \textbf{(proper) tail} $\sigma$-algebra of $R$ to be the intersection of all $\scrF_{\strK_1^\complement \times \strK_2^\complement}^R$ as $\strK_1$ and $\strK_2$ run over all finite subsets of $\Z^{d + 1},$ we denote it as $\scrT^R.$ Finally, we define the \textbf{restricted composition (relative to $\eta$)} $L \diamond R$ of two relations $L$ and $R$ (in that order) to be the set of all pairs $(z_1, z_2)$ such that there exists a $z \in \Z^{d + 1}$ satisfying the relations $z_1Lz$ and $zRz_2$ and $z$ satisfies
	\begin{displaymath}\textstyle
		\spread{z_1z_2} \leq \min(\spread{z_1z}, \spread{z_2z}).
	\end{displaymath}

	\begin{Theorem}{\normalfont(\cite[Theorem 3.3]{BKPS:GeomUSF})}\label{Theorem: Tail triviality of composition}
		Let $L$ and $R$ be two independent random relations on the vertices of $\Gamma_d(\lambda).$
		\begin{enumerate}
			\item If all left tail $\sigma$-fields $\scrL^L_z$ of $L$ are trivial and $R$ has a trivial tail $\sigma$-field $\scrT^R,$ then the restricted composition $L \diamond R$ has trivial left tail $\sigma$-fields $\scrL^{L \diamond R}_z.$

			\item If $\dim_\eta(L)$ and $\dim_\eta(R)$ both exist, then $L \diamond R$ has also stochastic dimension (relative to $\eta$) and we have $\dim_\eta(LR) = \min \left( \dim_\eta(L) + \dim_\eta(R) + 1, d + 1 \right).$

			\item If $\dim_\eta(L)$ and $\dim_\eta(R)$ both exist, all left tail $\sigma$-fields $\scrL^L_z$ of $L$ are trivial, all right tail $\sigma$-algebras $\scrR^R_z$ of $R$ are trivial and $\dim_\eta(L) + \dim_\eta(R) \geq d,$ then $	\inf\limits_{z_1, z_2} \Probability{z_1 L \diamond R z_2} = 1.$
		\end{enumerate}
		Let now $m \geq 2$ and $(R_i)_{i = 1, \ldots, m}$ be independent random relations on the vertices of $\Gamma_d(\lambda)$ such that $\dim_\eta(R_i)$ exists for each $1 \leq i \leq m.$
		\begin{enumerate}\setcounter{enumi}{3}
			\item Suppose that $\sum\limits_{i = 1}^m \dim_\eta(R_i) + m \geq d + 2,$ all left tails $\sigma$-fields $\scrL^{R_1}_z$ are trivial, each of $\scrT^{R_2}, \ldots, \scrT^{R_{m - 1}}$ are trivial tail $\sigma$-algebras, and all right tail $\sigma$-algebras $\scrR^{R_m}_z$ are also trivial. Then $\Probability{z_1 R_1 \cdots R_m z_2} = 1$ for all $z_1$ and $z_2$ vertices of $\Gamma_d(\lambda).$
		\end{enumerate}
	\end{Theorem}
	\begin{proof}
		The proofs of Theorem 3.3 and Corollary 3.4 of \cite{BKPS:GeomUSF} carry over \emph{mutatis mutandis} to our context.
	\end{proof}

	\subsection{The $\WSF$ connectedness relation}
	Let $\Gamma$ be any network. We will consider the following relation $U_\Gamma$ on the vertices of $\Gamma:$ ``$v$ and $v'$ belong to the same $\WSF$-component.'' We have the following very general result regarding tail triviality of $U_\Gamma.$
	\begin{Theorem}\label{Theorem: Triviality for left, right and tail sigma fields of USF relation}
		Suppose $\Gamma$ is any network satisfying the usual assumptions. Assume further that
		\begin{enumerate}
			\item The network random walk of $\Gamma$ is transient.
			\item The network possesses the Liouville property {\normalfont($\mathsf{LP}$)} \S \ref{Section: Liouville property}.
		\end{enumerate}
		Then, all left and right tail $\sigma$-algebras of $U_\Gamma$ are trivial. If, additionally to the previous hypotheses,
		\begin{enumerate}\setcounter{enumi}{2}
			\item The $\WSF$ of $\Gamma$ is one ended, that is to say, for almost every realisation, every tree in the $\WSF$ has one end.
		\end{enumerate}
		Then, the tail $\sigma$-field of $U_\Gamma$ is also trivial.
	\end{Theorem}
	For the proof see Theorem 4.5 of \cite{BKPS:GeomUSF}, together with Remark 4.6.

	We need one last ingredient. Consider two random subsets of vertices, $A$ and $Z,$ of a network $\Gamma.$ It is said that $Z$ \textbf{stochastically dominates} $A$ if there exists a coupling $\mu$ of $A$ and $Z$ such that $\mu(A \subset Z) = 1.$

	\begin{Theorem}{\normalfont (\cite[Theorem 4.1]{BKPS:GeomUSF})}\label{Theorem: Domination of the USF relation}
		Let $\Gamma$ be any network. Let $\mathfrak{F},$ $\avec{\mathfrak{F}}{0}{m}$ be independent samples of the $\WSF$ of the network $\Gamma.$ If $x$ is a vertex of $\Gamma,$ denote by $\strV_x^{\mathfrak{F}}$ the vertex set of the component of $x$ in $\mathfrak{F}.$ Fix any vertex $o$ of $\Gamma$ and set $\strV_0 = \strV_o^{\mathfrak{F}}.$ For $j \geq 1,$ define inductively $\strV_j$ to be the union of all vertex components of $\mathfrak{F}$ that are contained in, or adjacent to, $\strV_{j - 1};$ in other words, $\strV_j$ is the union of $\strV_x^{\mathfrak{F}}$ as $x$ runs on $\closure{\strV}_{j - 1}$ (closure in $\Gamma$). Let $\strQ_0 = \strV_o^{\mathfrak{F}_0}.$ For $j \geq 1,$ define inductively $\strQ_j$ to be the union of all vertex components of $\mathfrak{F}_j$ that intersect $\strQ_{j - 1};$ in other words, $\strQ_j$ is the union of all $\strV_x^{\mathfrak{F}_j}$ such that $\strV_x^{\mathfrak{F}_j} \cap \strQ_{j - 1} \neq \varnothing.$ Then, $\strV_m$ stochastically dominates $\strQ_m.$
	\end{Theorem}
	\begin{proof}
		It is a simple adaptation of Theorem 4.1 of \cite{BKPS:GeomUSF}, indeed in their first sentence change $B_R$ and replace it with $B_R = \{v; d(v, o) < R\},$ where $d$ is the graph-distance of $\Gamma.$ The remainder of their proof proceeds verbatim.
	\end{proof}

	\subsection{The main theorem}
	\begin{Theorem}\label{Theorem: The separation between components}
		With the notation at the beginning of the section, $\max\limits_{z, z' \in \Z^{d + 1}} D(z, z') = \left\lceil \frac{d - 2}{4} \right\rceil, \quad a.s.$
	\end{Theorem}
	\begin{proof}
		The proof follows the lines of that of Theorem 1.1 and Proposition 5.1 of \cite{BKPS:GeomUSF}. For convenience of the reader, we provide full details.

		When $d = 1, 2,$ we know that $\USF$ of $\Gamma_d(\lambda)$ is a tree and thus the separation between components is zero, this agrees with the formula given. Assume then that $d \geq 3.$

		Let $m = \left\lceil \frac{d - 2}{4} \right\rceil.$ Consider $m + 1$ independent copies of the $U_{\Gamma_d(\lambda)}$ relation. By (\ref{Remark: The dimension of the USF relation}), we know that $\dim_\eta\left( U_{\Gamma_d(\lambda)} \right) = 3.$ Next, by (\ref{Theorem: Triviality for left, right and tail sigma fields of USF relation}), bearing in mind that $\USF$ is one ended (see sect. \ref{Section: One endedness}) this random relation possesses trivial tail $\sigma$-field as well as trivial left tail and right tail $\sigma$-algebras. We may apply item (d) of (\ref{Theorem: Tail triviality of composition}) to conclude that, if $R$ is the composition of the $m + 1$ copies of $U_{\Gamma_d(\lambda)}$ then $\inf\limits_{z_1, z_2 \in \Z^{d + 1}} \Probability{z_1Rz_2} = 1.$ Finally, define $S$ to be the relation $D(z, z') \leq m.$ By virtue of (\ref{Theorem: Domination of the USF relation}), we know that $\Probability{z_1Sz_1} \geq \Probability{z_1Rz_2}$ and thus, $\Probability{D(z, z') \leq m, \quad \forall z,z'} = 1.$

		The proof of the cases $3 \leq d \leq 6$ is complete, so it only remains to prove the lower bound for the cases $d \geq 7.$ To prove the lower bounds we will establish first
		\begin{equation}\tag{$\ast$}
			\Probability{D(z, z') \leq k} \leq c \spread{zz'}^{4k - (d - 2)}
		\end{equation}
		for all $z, z' \in \Z^{d + 1}$ and $k \in \N,$ the constant $c$ depending on $k$ and $d,$ but not on $z$ or $z'.$ When $k \geq m,$ the criterion holds since then $4k \geq d - 2 \geq 5$ and $\spread{zz'} \geq 1.$ Assume $k < m,$ so that $4k - (d - 2) \leq -1.$ Consider a finite sequence $(z_j, z_j')_{j = 0, \ldots, k}$ of pairs of vertices of $\Gamma_d(\lambda)$ with $z_0 = z,$ $z_k' = z'$ and $z_j' \sim z_{j + 1}$ for $j = 0, \ldots, k - 1.$ Let $\EuScript{A} = \EuScript{A}(z_j, z_j')$ be the event in which $z_j'$ is a vertex of $\mathfrak{T}_{z_j}$ for $j = 0, \ldots, k$ and for each of these $j,$ $z_j'$ is not a vertex of $\mathfrak{T}_{z_i}$ for $i \neq j.$ Then,
		\begin{displaymath}\textstyle
			\{D(z, z') = k\} \subset \bigcup\limits_{(z_j, z_j')} \EuScript{A}(z_j, z_j')
		\end{displaymath}
		where the pairs $(z_j, z_j')$ run over all sequences $(z_j, z_j')$ as stated. To estimate the probability of a particular $\EuScript{A}(z_j, z_j')$ we may construct the $\USF$-forest using Wilson's rooted at infinity, beginning with the random walks started at $z_0, z_0', \ldots, z_k, z_k'.$ For the event $\EuScript{A}(z_j, z_j')$ to hold, it is necessary that for each $j = 0, \ldots, k,$ the random walk started at $z_j$ must intersect the path of the random walk started at $z_j'.$ Hence, if $S^z$ denotes the random walk started at $z$ (with the different random walks mutually independent),
		\begin{displaymath}\textstyle
			\Probability{\EuScript{A}} \leq \prod\limits_{j = 0}^k \Probability{S^{z_j} \cap S^{z_j'} \neq \varnothing} = \prod\limits_{j = 0}^k \sum\limits_{z \in \Z^{d + 1}} G(z_j, z) G(z_j', z) \asymp \prod\limits_{j = 0}^k \spread{z_jz_j'}^{-(d - 2)},
		\end{displaymath}
		where $\asymp$ follows from (\ref{Lemma: Expectation of Kz is proportional to eta of z to the power of negative d minus two}). Since $z_j' \sim z_{j + 1},$ we can bound
		\begin{displaymath}\textstyle
			\Probability{D(z, z') = k} \leq c \sum\limits_{(z_j)_{j = 0}^{k + 1}} \prod\limits_{j = 0}^k \spread{z_jz_{j + 1}}^{-(d - 2)},
		\end{displaymath}
		the sum extending over all finite paths $(z_j)$ starting at $z$ and ending at $z'.$ We may finally break up this sum, starting with the sum over all $z_k$ and applying (\ref{Proposition: Lemma 2 point 10 of BKPS04}) with $\alpha = \beta = d - 2,$ so that $\gamma = \alpha + \beta - d - 2 = d - 6 > 0,$ and then continue with the sum over $z_{k - 1},$ and so on until $z_1,$ this gives the desired bound ($\ast$). Having establish ($\ast$), we can terminate Theorem (\ref{Theorem: The separation between components}). Indeed, consider $z \in \Z^{d + 1},$ the event $\ds \left\{ \max_{z' \in \Z^{d + 1}} D(z, z') = k \right\}$ belongs to the left-tail at $z$ relative to the relation $U_{\Gamma_d(\lambda)},$ as such it is a trivial event (\ref{Theorem: Triviality for left, right and tail sigma fields of USF relation}). It is clear from ($\ast$) that this event has probability zero, therefore the event $\ds \left\{\max_{z, z' \in \Z^{d + 1}} D(z, z') = k \right\}$ also have probability zero.
	\end{proof}

	\bibliographystyle{alpha}
	\nocite{*}
	\bibliography{PhDBibliography}
\end{document}